\documentclass[11pt]{report}
\usepackage{amsmath,amsfonts,amssymb,amsthm,mathrsfs,mathtools}
\usepackage[dvipsnames]{xcolor}
\usepackage[colorlinks,allcolors=Green]{hyperref}
\usepackage{array}
\usepackage{tikz}
\usepackage{marginnote}

\title{Small resolutions of \\ special three-dimensional varieties}
\author{J\"{u}rgen Werner \\ translated by Simon Venter and Nicolas Addington\footnote{simonv@uoregon.edu, adding@uoregon.edu \newline
Dept.\ of Mathematics, University of Oregon, Eugene, OR 97403, United States  \newline
N.A.\ was partially supported by NSF grant no.\ DMS-1902213.}}
\date{} 

\usepackage{enumitem}
\setlist{  
  listparindent=\parindent,
  parsep=0pt,
}

\theoremstyle{definition}
\newtheorem{thm}{Theorem}[chapter]
\newtheorem{dfn-thm}[thm]{Definition and Theorem}
\newtheorem{remark}[thm]{Remark}
\newtheorem{lem}[thm]{Lemma}

\newtheorem{conj}[thm]{Conjecture}
\newtheorem{dfn}[thm]{Definition}
\newtheorem{exa}[thm]{Example}
\newtheorem{cor}[thm]{Corollary}
\newtheorem{prop}[thm]{Proposition}

\counterwithout{footnote}{chapter}

\reversemarginpar

\newcommand \wernerpage[1] {\marginnote{#1}}

\makeatletter
\def\@makechapterhead#1{%
  {\parindent \z@ \raggedright \normalfont
    \ifnum \c@secnumdepth >\m@ne
        \LARGE\bfseries \Roman{chapter}.
    \fi
   #1\bigskip
  }}
\makeatother

\begin{document}

\pagenumbering{roman}

\renewcommand \thefootnote {\fnsymbol{footnote}}
\maketitle

\newcommand \httpurl [1] {\href{https://#1}{\nolinkurl{#1}}}

\paragraph{Translators' note.} \ 

J\"urgen Werner's dissertation, written in 1987 under the supervision of Friedrich Hirzebruch, explores nodal threefolds and their small resolutions, which continue to interest algebraic geometers and string theorists a generation later.  We hope this translation will make it accessible to a wider audience.  We thank Werner for permission to translate his dissertation, Rachel Rodman for advice on several passages, Manfred Lehn for advice on the word ``speziell,'' and David Levin for help with the end of chapter \ref{chapter7}.

The original was published as Bonner Mathematische Schriften No.\ 186, and is available online at \httpurl{archive.mpim-bonn.mpg.de/id/eprint/2784}.  The original page numbers appear in the margins. The theorem numbers in this version do not appear in the original.  We have corrected typos, and in footnotes commented on a few minor errors and added a few clarifying remarks, but we have not attempted to check every detail.  We hope that we have not introduced errors, but we welcome corrections from readers, and will update the arXiv file with any corrections that we receive.

\newpage
\tableofcontents
\newpage
\renewcommand \thefootnote {\arabic{footnote}}
\chapter*{Introduction}
\addcontentsline{toc}{chapter}{Introduction}
\wernerpage{i}

The present work deals with some 
complex three-dimensional singular varieties and the various ways of resolving their singularities. 
Our main interest is in so-called ``small" resolutions: the singularities are not replaced by a surface, but by a curve, thus by a set of codimension 2.  
This is possible for the first time in dimension 3.

For a long time, mathematicians have mainly been concerned with the local properties of small resolutions (cf.\ Pinkham’s article \cite{pinkham83} and Laufer’s article \cite{laufer81}). 
Brieskorn also describes small resolutions of three-dimensional singularities in his articles \cite{brieskorn66} and \cite{brieskorn68}.

In this work, we will investigate the global properties of small resolutions, especially the question of projectivity. \bigskip

Why are we interested in the global properties of small resolutions?

For some time, physicists such as E. Witten (Princeton) have sought examples of three-dimensional complex K\"{a}hler manifolds which have $c_1^{\mathbb R} = 0$ and absolute value of the Euler number as small as possible, but not equal to zero.
Some canonical examples with $c_1^{\mathbb R} = 0$ include the complete intersection of $k$ hypersurfaces of degrees $d_1, \dotsc, d_k$ in $\mathbb P^{3 + k}(\mathbb C)$ with $\sum_{i = 1}^k d_i = k + 4$ and double covers of $\mathbb P^3(\mathbb C)$ branched over surfaces of degree 8. 
These examples have the disadvantage that their Euler numbers are between $-300$ and $-100$.

\wernerpage{ii}

One way to obtain examples with smaller Euler numbers is to let a group act freely 
on the examples above and take the quotient (see \cite{candelas85} and \cite{strominger85}).

Another means for altering Euler numbers was found by F.~Hirzebruch: he deforms the examples above to singular varieties and then resolves the singularities (cf. \cite{hirzebruch85} and \cite{hirzebruch86}). 
In some cases where a ``big'' resolution would change the first Chern class, one must resolve the singularities with curves.

This raises the question: are these small resolutions K\"ahler?

Although this was the starting point of the present work, we will not confine ourselves to examples with $c_1^{\mathbb R} = 0$, but will take this opportunity to study global properties of small resolutions in general.

More examples with $c_1^{\mathbb R} = 0$ can be found in \cite{schoenIV} and \cite{yau85}. \bigskip

The mathematical basis of my considerations was the work ``Double Solids'' by H.~Clemens, in which, among other things, the homology groups of various resolutions of double solids are described.
Clemens' results can also be transferred to nodal hypersurfaces in $\mathbb P^4(\mathbb C)$ (and to some extent to arbitrary nodal threefolds and even to varieties with higher singularities).

\wernerpage{iii}

The overview 
of Clemens’ article with regard to our problem (see Chapter \ref{chapter1}) takes place in chapter \ref{chapter2} of this work. 
Of particular interest is the ``defect'' defined by Clemens, which can best be introduced as the difference between the fourth and second Betti numbers of the singular variety.
Furthermore, we work out how changing the small resolution affects the homology groups.
Altogether, a nodal variety $V$ allows exactly $2^s$ different small resolutions, where $s = \# V_{\text{sing}}$. \bigskip

Chapter \ref{chapter3} contains our main theoretical results for double solids and nodal hypersurfaces in $\mathbb P^4$: we characterize the projective algebraic small resolutions, and derive a necessary and sufficient criterion for the existence of at least one projective algebraic small resolution for a given nodal variety $V$. 
It turns out that projective algebraic small resolutions exist if and only if all irreducible, exceptional curves over $\mathbb Q$ are not nullhomologous (this property is independent of the particular 
small resolution). \bigskip

General considerations of Moishezon \cite{moisezon67}, Harvey-Lawson \cite{harvey83}, Siu \cite{siu74}, and Peternell \cite{peternell86} enter here. 
We give several versions of the proof, each time trying to build less on general results and instead making direct inferences from our particular situation. \bigskip

\wernerpage{iv}

In Chapter \ref{chapter4}, we investigate how our newfound necessary and sufficient criterion for existence of at least one projective algebraic small resolution can be applied in concrete situations. 
Here the defect plays an important role. 
One easily sees that a threefold with defect 0 admits no projective algebraic small resolution 
(the converse is unfortunately false: there exist varieties with arbitrarily large defect, but with no projective algebraic small resolution; see Chapter \ref{chapter9}). \bigskip

A method for computing defects, stated by Clemens for double solids and derived in Chapter \ref{chapter4} for nodal hypersurfaces in $\mathbb P^4 (\mathbb C)$, is therefore worth mentioning in this context. 
It turns out that for generic nodal varieties of the kinds we consider, the defect vanishes: only double solids and hypersurfaces with many ordinary double points in very special position have positive defect, and only in these exceptional cases is there a chance of getting projective algebraic small resolutions.

As it will transpire in the course of the work, however, there are enough of these ``exceptional cases'' to warrant a general study: there are infinitely many examples with positive defect, and the defect can be arbitrarily large.\bigskip

Chapters \ref{chapter5} and \ref{chapter6} focus on a large class of interesting examples: Chmutov hypersurfaces in $\mathbb P^4$ and double solids branched over Chmutov surfaces in $\mathbb P^3$.

\wernerpage{v}

In \cite{varchenko832}, Varchenko describes how Chmutov uses Chebyshev polynomials to construct hypersurfaces in $\mathbb P^n(\mathbb C)$ with many ordinary double points. 
This construction can be further generalized through sign changes; in the case of double solids, the choice of signs is crucial for answering the question of whether projective algebraic small resolutions exist.

Chapter \ref{chapter5} examines the various families in general with regard to this question, while Chapter \ref{chapter6} singles out some special cases. 
The defects are calculated (in part with the help of a computer), and in various examples by knowing the defect one can show the non-existence of a projective algebraic small resolution. \bigskip

Chapter \ref{chapter7} is only indirectly related to small resolutions, as it deals with the maximal number $\mu_n(d)$ of nodes that a hypersurface of degree $d$ in $\mathbb P^n(\mathbb C)$ can have.
Among other things, we observe the asymptotic behavior as $d \to \infty$ and $n \to \infty$, and also give a current table of estimates for $\mu_n(d)$ for $n = 3, 4$ and small values of $d$ (cf.\ also \cite{gallarati83}). \bigskip

Chapters \ref{chapter8}, \ref{chapter9}, and \ref{chapter10} feature further examples of nodal threefolds.

Chapter \ref{chapter8} deals with cubics in $\mathbb P^4(\mathbb C)$, with results developed in collaboration with mathematicians from Leiden University;
\wernerpage{vi}
a thorough treatment of cubics (especially their geometric properties) will be included in a joint project with Hans Finkelnberg. \bigskip

An important tool for studying cubics is the associated curve, from which one can read off everything that interests us in the context of small resolutions. 
Chapter \ref{chapter9} describes to what extent the principle of the associated curve can also be carried over to quartics and special hypersurfaces of higher degree.

A quartic with 45 nodes, studied in detail by Todd in \cite{todd47}, admits projective algebraic small resolutions. 
A proof that $\mu_4(4) = 45$ can be found in \cite{friedman86}; thus 45 is the maximum number of nodes that a quartic in $\mathbb P^4(\mathbb C)$ can have, as can also be seen from Varchenko's general considerations in \cite{varchenko831}.
It remains an open question whether there are any quartics with 45 ordinary double points other than Todd’s example. \bigskip

Chapter \ref{chapter10} contains further examples of nodal hypersurfaces and double solids. 
Especially interesting is how geometric properties of some quartic surfaces, studied by Kummer over 100 years ago, can be recovered from the corresponding double solids or small resolutions.
Finally, we investigate the minimum number of ordinary double points that a hypersurface or a double solid must have in order to obtain positive defect. \bigskip

\wernerpage{vii}

The last two chapters attempt to generalize the previous results further. 
First, arbitrary nodal threefolds will be allowed, then varieties with certain higher singularities.
Local information about small resolutions, like that contained in Brieskorn's articles \cite{brieskorn66} and \cite{brieskorn68}, will be put together to give global information.
The aim is not to set out a comprehensive theory, but rather to give a feeling for the extent to which the double point case is an expression of a general principle.

It turns out this principle is not so easy to pin down using nullhomologous curves, but works better with the help of smooth divisors containing certain singularities. This is the method used by Schoen in the case of ordinary double points (see \cite{schoen86} and \cite{schoenIV}). 
Brieskorn's consideration 
point the way toward a generalization.
This also makes clear the close connection with nodal threefolds arising as deformations of varieties with higher singularities.

Our first examples are again cubics in $\mathbb P^4 (\mathbb C)$.
The general principle can here be read off from the the associated curve, which contains important information about small resolutions.
We also study Kummer surfaces in $\mathbb P^3 (\mathbb C)$ with $D_4$-singularities and the corresponding double solids. \bigskip

\wernerpage{viii}

I would like to sincerely thank Prof.\ F.~Hirzebruch for suggesting this project and for the opportunity to work at the Max Planck Institute for Mathematics in Bonn. 
My thanks go to the Max Planck Society for financial support and to the many guests and staff of the Institute for professional support, especially Ms.\ H.~Konrad for her help with computer calculations. 
I learned a great deal from the hospitable mathematicians of Leiden University, and am indebted to them for their important advice. 
Most notably, I would like to sincerely thank Prof.\ W.~Meyer for his interest in my work and for many valuable conversations. 
Without his help this work (at least in its present form) would not have come into being.

\chapter{Problem statement} \label{chapter1}
\pagenumbering{arabic}
\wernerpage{1}

Our investigations focus on two classes of complex three-dimensional singular varieties:
\begin{enumerate}[label=\arabic*.]
     \item Hypersurfaces in $\mathbb P^{4}( \mathbb C )$;
     \item Double solids, which are double covers of $\mathbb P^{3}( \mathbb C )$ branched over singular surfaces, as studied by Clemens in \cite{clemens83}.
\end{enumerate}

For singularities, we will only have isolated ordinary double points, which are locally described by the equation
\[
\sum_{ i = 1 }^{ 4 } z_i^2 = 0. 
\]
In the last two chapters of this work, higher singularities will also be allowed.\bigskip

Let $V$ be a nodal threefold of either of the types described above.  In the first case, let $n$ denote the degree of $V$ and, in the second, let $n$ denote the degree of the nodal branching surface $B \subset \mathbb P^{3}( \mathbb C )$, which must be even.
The singularities of $B$ are given locally by an equation of the form
\[
\sum_{ i = 1 }^{ 3 } z_i^{2} = 0 .
\]
In what follows the singularities of $B$ will be identified with those of the corresponding double solid.

\wernerpage{2}

The set of singularities will be denoted by $\mathscr{S}$, and their number by $s$.
There are two ways to resolve a singularity $P \in \mathscr{S}$.
The first, which we call a \emph{big resolution}, can be described as follows:

\begin{enumerate}[label=\arabic*.]
    \item If $P \in V \subset \mathbb P^{4}$, blow up $\mathbb P^{4}$ at $P$; the proper transformation of $V$ contains, instead of $P$, a smooth quadric in $\mathbb P^{3}$, which is isomorphic to $\mathbb P^{1} \times \mathbb P^{1}$.
    
    \item If $P \in B \subset \mathbb P^{3}$ is a singularity of a double solid, blow up $\mathbb P^3$ at $P$; the double cover of this [blow-up of] $\mathbb P^3$ branched over the proper transform of $B$ contains, instead of $P$, a double cover of $\mathbb P^{2}$ branched over a smooth conic.
    The exceptional surface is again isomorphic to $\mathbb P^{1} \times \mathbb P^{1}$.
\end{enumerate}

In both cases, the two rulings of the $\mathbb P^1 \times \mathbb P^1$ have intersection number $-1$ with the exceptional surface. 
Thus one can blow down either fibration of the $\mathbb P^1 \times \mathbb P^{1}$ to the corresponding base without introducing singularities. \bigskip

Starting from the singular variety, one obtains these \emph{small resolutions}, where
$P \in \mathscr{S}$ is only replaced by a curve isomorphic to $\mathbb P^{1}$, as follows.
By a change of coordinates, the local equation of a singularity
\[
\sum_{ i  = 1 }^{ 4 } z_i^{2} = 0
\]
can be put in the form
\[
\Phi_1 \Phi_2 = \Phi_3 \Phi_4.
\]
\wernerpage{3}
The local meromorphic function
\[
\frac{\Phi_1}{\Phi_3} = \frac{\Phi_2}{\Phi_4}
\]
has a point of indeterminacy at $(0, 0, 0, 0)$.
The closure of its graph is smooth and contains a copy of $\mathbb P^1$ in place of $P$.
The graph of the function
\[
\frac{\Phi_1}{\Phi_4} = \frac{\Phi_3}{\Phi_2}
\]
yields the other small resolution, which arises by blowing down the other ruling of the $\mathbb P^1 \times \mathbb P^1$. \bigskip

Thus one has two distinct small resolutions for each singular point.  They can be exchanged by transposing the coordinates $\Phi_3$ and $\Phi_4$.
In $z$-coordinates, this means
\[ z_2 \longrightarrow -z_2 \]
and $z_i$ is unchanged for $i \neq 2$.
The general principle will be clarified in Chapter \ref{chapter11}; see also \cite{brieskorn66} and \cite{hirzebruch85}. \bigskip

A manifold obtained from $V$ by simultaneous small resolutions of each singularity $P \in \mathscr{S}$ will be denoted $\hat V$.
Altogether, this gives $2^s$ such \emph{small} resolutions of $V$.
By a \emph{partial} resolution $\check{V}$ of $V$, we will mean a variety in which only the singularities $P \in \mathscr{S}_1 \subset \mathscr{S}$ are given a small resolution.

\wernerpage{4}

The manifold obtained by a simultaneous big resolution of all the singularities of $V$ is called $\tilde{V}$.
Given a small resolution $\hat{V}$, one can construct the \emph{big} resolution $\tilde{V}$ by blowing up along all the exceptional curves.
If one performs a big resolution of all $P \in \mathscr{S}_1 \subset \mathscr{S}$ and a small resolution of the others, one gets a \emph{mixed} resolution $\overline{V}$.

The singular variety $V$ is embedded in a projective space by assumption.
The big resolution $\tilde{V}$ is also always projective algebraic.
But the projectivity of the small resolutions $\hat{V}$ is an interesting question.
As we shall see, this question is closely tied to the projectivity of certain 
partial and mixed resolutions. \bigskip

The central problem is to find a necessary and sufficient criterion for a given singular variety $V$ that guarantees the existence of at least one projective algebraic small resolution $\hat{V}$. \bigskip

The following study of the homology groups of various $V$, $\check{V}$, $\hat{V}$, $\overline{V}$, and $\tilde{V}$ represents a first step toward such a criterion.
\chapter{Homology of \texorpdfstring{$V$}{V} and related varieties} \label{chapter2}
\wernerpage{5}

The considerations in this chapter are primarily based on the first sections of H.~Clemens' paper ``Double Solids" \cite{clemens83}.
We assume his results without proof. \bigskip

To determine the homology groups of $V$, it is useful to regard $V$ as a singular member of a family of threefolds $\{V_{t}\}$ ($|t| < \epsilon$), where $V_t$ is smooth for $t \neq 0$ and $\deg V_t = \deg V$, or $V_t$ is a double solid branched over a smooth surface $B_t$ such that $\deg B_t = \deg B$, respectively.

The transition from $V_t$ for $t \neq 0$ to $V_0 = V$ corresponds in an affine neighborhood of $P \in \mathscr{S}$ with the contraction of the real $3$-sphere
\[
\alpha( P, t ) = \left\{ \sum_{i = 1}^{4} z_i^{2} = t, \frac{z_i}{\sqrt{t}} \in \mathbb{R}  \text{ for } i \in \{1, \dotsc, 4\} \right\}.
\]
Clemens describes this transition as follows:

By attaching a 4-cell along every vanishing 3-cycle $\alpha( P, t )$, one obtains from $V_{t}$ a space $X$ whose homology groups are isomorphic to those of $V$.
From this, Clemens deduces the exact sequence
\[
0 \to H_4( V_{t}, \mathbb{Z} ) \to H_4( V, \mathbb{Z} ) \xrightarrow{k} \mathscr M \to H_3( V_{t}, \mathbb{Z} ) \to H_3( V, \mathbb{Z} ) \to 0
\]
as well as the relation
\[
H_2( V_{t}, \mathbb{Z} ) \cong H_2( V, \mathbb{Z} ).
\]
\wernerpage{6}Here $\mathscr M$ denotes the free $\mathbb Z$-module generated by $\mathscr S$, and the map
\[ \mathscr M  \to H_3( V_{t}, \mathbb{Z} ) \]
is given by
\[ P \mapsto \alpha( P, t ). \]

The exactness of the sequence above implies that the 4-cycles on $V$ that are not images of 4-cycles on $V_{t}$ correspond exactly to global linear relations
\[
\sum_{ i = 1 }^s m_i \alpha( P_i, t ) \simeq 0
\]
between the vanishing cycles on $V_{t}$.
The submodule $\mathscr A$ of $\mathscr M$ spanned by the $s$-tuples $(m_1, \dots, m_s)$ is the image of $H_4( V, \mathbb{Z} )$ under the map $k$.
Let $d$ be the rank of the module $\mathscr A$, which Clemens calls the \emph{defect}.
We calculate the Betti numbers of $V$ to be
\begin{align*}
    \beta_2( V ) &= 1 \\
    \beta_4( V ) &= 1 + d \\
    \beta_3( V ) &= \beta_3( V_{t} ) - s + d,
\end{align*}
where
\[
\beta_3( V_{t} ) =
\begin{cases}
n^{4} - 5 n^{3} + 10 n^{2} - 10 n + 4 & \text{if } n = \deg V_{t} \\
n^{3} - 4 n^{2} + 6 n - 4 & \text{if } n = \deg B_{t}
\end{cases}
\]
For the Euler number, we have
\[
e( V ) = e( V_{t} ) + s.
\]

\wernerpage{7}

We proceed to the small resolution $\hat{V}$, so singular points $P \in \mathscr S$ are replaced exceptional curves $L_P$ isomorphic to $\mathbb P^{1}$.
Clemens describes the dual process as follows:

Let $Y$ be the space obtained from $\hat V$ by attaching a 3-cell along each curve $L_P$; then $Y$ is homotopy equivalent to $V$, so for all $q \ge 0$ we have
\[
H_q( Y, \mathbb Z ) \cong H_q( V, \mathbb Z ).
\]
We obtain the exact sequence
\begin{align*}
    0 \to H_3( \hat V, \mathbb Z ) \to H_3( V, \mathbb{Z} ) \xrightarrow{k'} \mathscr M &\to H_2( \hat V, \mathbb Z ) \to H_2( V, \mathbb Z ) \to 0 \\
    P &\mapsto \{ L_P \}
\end{align*}
and the relation
\[
    H_4( \hat V, \mathbb Z ) \cong H_4( V, \mathbb Z ).
\]
We define a bilinear form on $\mathscr M$ by $( P_i, P_j ) := \delta_{ij}$, with respect to which $\mathscr M$ is self-dual.
For $\hat V$ we get the following Betti numbers:
\begin{align*}
    \beta_2( \hat V ) &= 1 + d \\
    \beta_4( \hat V ) &= 1 + d \\
    \beta_3( \hat V ) &= \beta_3( V ) - s + d \\
                        &= \beta_3( V_{t} ) - 2s + 2d.
\end{align*}
The Euler number of $\hat V$ is thus equal to
\begin{align*}
    e(\hat V) &= e( V ) + s \\
         &= e( V_{t} ) + 2s.
\end{align*}

\wernerpage{8}

Between the exceptional curves on $\hat V$ there must be $(s - d)$ linear relations
\[
    \sum_{ i = 1 }^s n_i L_{P_i} \simeq 0.
\]
Define $\mathscr B$ as the rank-$(s-d)$ submodule of $\mathscr M$ spanned by the $s$-tuples $(n_1, \dotsc, n_s)$. \bigskip

This raises the following questions:
What is the relationship between $\mathscr A$ and $\mathscr B$?
How do $\mathscr A$ and $\mathscr B$ depend on the particular choice of small resolution?
These questions will now be answered. \bigskip

As Clemens proves, we can obtain $\hat V$ directly from $V_{t}$ by surgery along the vanishing cycle $\alpha( P, t )$.
The particular small resolution of a point $P \in \mathscr S$ is determined by the orientation of the vanishing cycle $\alpha( P, t )$.
In the local $z$-coordinates, a transformation
\[
    z_2 \mapsto -z_2
\]
causes a change in the small resolution of the point $P$ as well as in the orientation of $\alpha( P, t )$.
Thus the choices of small resolution of a point $P \in \mathscr S$ are in bijection with the choices of orientation of the associated vanishing cycle $\alpha( P, t )$ on $V_{t}$. \bigskip

\wernerpage{9}

Clemens deduces that the intersection number of a $\gamma \in H_4( \hat V, \mathbb Z )$ with a curve $L_P$ is
\[
    ( \gamma . L_P ) = ( k( \gamma ), P ),
\]
where
\[
    k \colon H_4( \hat V, \mathbb Z ) \cong H_4( V, \mathbb Z ) \to \mathscr M
\]
is the map given in the exact sequence above, and $(.,.)$ is the bilinear form on $\mathscr M$ already defined.
For 4-cycle $\gamma$ with
\[
    k( \gamma ) = \sum_{i = 1}^s m_i P_i,
\]
we have
\[
    (\gamma . L_{P_j}) = m_j.
\]

In particular, for every $s$-tuple $( n_1, \dotsc, n_s ) \in \mathscr B$ we have
\[
    0 = \gamma . \left( \sum_{i = 1}^s n_i L_{P_i} \right)  = \sum_{i = 1}^s m_i n_i.
\]
With regard to the bilinear form on $\mathscr M$ determined by $(P_i, P_j) = \delta_{ij}$, we see that $\mathscr A$ and $\mathscr B$ are orthogonal to one another.
Note that $\operatorname{rank} \mathscr M = \operatorname{rank} \mathscr A + \operatorname{rank} \mathscr B$. \bigskip

In general there are elements of $\mathscr M$ that cannot be expressed as a linear combination of vectors from $\mathscr A$ and $\mathscr B$ with integer coefficients.
Also, torsion groups can occur in the homology groups over $\mathbb Z$.
Therefore we tensor $\mathscr A$ and $\mathscr B$ with $\mathbb Q$, and for simplicity we let these $\mathbb Q$-vector spaces be denoted again by $\mathscr A$ and $\mathscr B$.
We have
\[
    \mathbb Q^s = \mathscr A \perp \mathscr B.
\]

In what follows, torsion is irrelevent to our considerations.
\wernerpage{10}
Every statement about homology --- for example, that $C \simeq 0$ for some curve $C$ --- is henceforth understood to be over $\mathbb Q$. \bigskip

We return to the small resolution $\hat V$ and recall that a change of small resolution at a point $P_q \in \mathscr S$ corresponds to a change in the orientation of the cycle $\alpha( P_q, t )$ on $V_{t}$.
In homology, this means
\[
    \alpha( P_q, t ) \longrightarrow -\alpha( P_q, t ),
\]
and thus
\[
    m_q \longrightarrow -m_q
\]
in all tuples $( m_1, \dots, m_s ) \in \mathscr A$.

Since
\[
    \sum_{i = 1}^s m_i n_i = 0
\]
for all $( m_1, \dotsc, m_s ) \in \mathscr A$ and all $( n_1, \dotsc, n_s ) \in \mathscr B$, it follows that
\[
    n_q \longrightarrow -n_q
\]
for all $n_q$ in the tuples $( n_1, \dots, n_s ) \in \mathscr B$. \bigskip

Thus the decomposition of $\mathbb Q^s$ into $\mathscr A$ and $\mathscr B$ depends in general on the particular 
choice of small resolution. 
Unfortunately, there is no intrinsic description of $\mathscr A$ and $\mathscr B$ independent of the small resolution, so in examples, one must specify the small resolution used to decompose $\mathbb Q^s$ into $\mathscr A$ and $\mathscr B$.
\wernerpage{11}
A change of small resolution at a singularity causes a change in the signs of the associated coefficients in all the relations
\[
    \sum_{i = 1}^s m_i \alpha( P_i, t ) \simeq 0
\]
and
\[
    \sum_{i = 1}^s n_i L_{P_i} \simeq 0.
\]
However, if a singularity $P_q$ is replaced with a nullhomologous curve, this is independent of the particular 
small resolution; in this case, $m_q = 0$ in all the $s$-tuples in $\mathscr A$, and $\mathscr B$ contains of the vectors $\pm e_q$.
The defect and hence the Betti numbers are also independent of the particular 
small resolution. \bigskip

One can get the big resolution $\tilde V$ from any small resolution $\hat V$ by blowing up along all the curves $L_P$.
In $\tilde V$, the exceptional surface lying over $P$ is called
\[
Q_P \cong A_P \times B_P:
\]
$A_P$ corresponds to the ruling that is contracted to a point in the map to $\hat V$, and $B_P$ to the ruling that is mapped to $L_P$.
If one starts from another small resolution, then the relevant $A_P$ and $B_P$ are interchanged accordingly.

\wernerpage{12}

Clemens describes the relationship between the homology groups of $\hat V$ and $\tilde V$ with the following exact sequences:
\begin{align*}
    0 \to \mathscr M &\to H_2( \tilde V, \mathbb Z ) \to H_2( \hat V, \mathbb Z ) \to 0 \\
                   P &\mapsto [ A_P ]
\end{align*}
\begin{align*}
    0 \to \mathscr M &\to H_4( \tilde V, \mathbb Z ) \to H_4( \hat V, \mathbb Z ) \to 0 \\
                   P &\mapsto [ Q_P ]
\end{align*}
\[ H_3( \tilde V, \mathbb Z ) \cong H_3( \hat V, \mathbb Z ) \]
From this we calculate the Betti numbers of the big resolution to be
\begin{align*}
    \beta_2( \tilde V ) &= 1 + d + s \\
    \beta_4( \tilde V ) &= 1 + d + s \\
    \beta_3( \tilde V ) & = \beta_3( \hat V ).
\end{align*}
The Euler number of the big resolution is equal to
\begin{align*}
    e( \tilde V ) &= e( \hat V ) + 2s \\
                  &= e( V_{t} ) + 4s.    
\end{align*}

Between the curves $A_P$ and $B_P$ (or more precisely, between their homology classes) there are again $(s - d)$ linear relations.
Since
\[
    Q_{ P_i } . A_{ P_j } = Q_{ P_i } . B_{ P_j } = -\delta_{ij},
\]
these relations are of the form
\[
    \sum_{i = 1}^s n_i ( B_{ P_i } - A_{ P_i } ) \simeq 0,
\]
where the tuple $( n_1, \dotsc, n_s )$ again ranges over $\mathscr B$.

\wernerpage{13}

Here it once again becomes clear: a change of the small resolution at a point $P_q \in \mathscr S$ --- appearing on the big resolution an exchange of $A_{ P_q }$ and $B_{ P_q }$ --- causes a sign change
\[
    n_q \mapsto -n_q
\]
in all tuples $( n_1, \dotsc, n_s ) \in \mathscr B$.
A relation of the form
\[
    L_P \simeq 0
\]
exists on a small resolution $\hat V$, and thus on all $\hat V$, precisely when
\[
    A_P \simeq B_P
\]
on $\tilde V$. \bigskip

How do the partial and mixed resolutions fit into this scheme?

Let $\mathscr S_1$ be a subset of $\mathscr S$, let $s_1 = \#\mathscr S_1$, and let $\overline V$ be a mixed resolution where we resolve the singularities in $\mathscr S_1$ with a big resolution and the rest with a small resolution.
Let $\mathscr M_1$ and $\mathscr M_2$ be the free $\mathbb Z$-modules spanned by $\mathscr S_1$ and $\mathscr S - \mathscr S_1$, respectively.
Using the exact sequences
\[
    0 \to \mathscr M_2 \to H_2( \tilde V, \mathbb Z ) \to H_2( \overline V, \mathbb Z ) \to 0
\]
and
\[
    0 \to \mathscr M_1 \to H_2( \overline V, \mathbb Z ) \to H_2( \hat V, \mathbb Z ) \to 0,
\]
\wernerpage{14}
we calculate the Betti numbers of $\overline V$ as
\begin{align*}
    \beta_2( \overline V ) &= 1 + d + s_1 \\
    \beta_4( \overline V ) &= 1 + d + s_1 \\
    \beta_3( \overline V ) &= \beta_3( \hat V ).
\end{align*}
The Euler number of $\overline V$ is equal to
\begin{align*}
    e( \overline V ) &= e( \hat V ) + 2 s_1 \\
                     &= e( V_{t} ) + 2s + 2s_1.
\end{align*}
The $(s - d)$ relations between the $L_P$, $A_P$, and $B_P$ are
\[
    \sum_{ P_i \in \mathscr S_1 } n_i \left( A_{ P_i } - B_{ P_i } \right) + \sum_{ P_i \in \mathscr S - \mathscr S_1 } n_i L_{ P_i } \simeq 0,
\]
where $( n_1, \dotsc, n_s ) \in \mathscr B$. \bigskip

The Betti numbers of a partial resolution $\check V$ are unfortunately not so easy to calculate; $\beta_2( \check V )$ depends in general on which singularities are resolved.
But by considering the map on second homology induced by the natural map
\[
\hat V \to \check V,
\]
one sees that between the exceptional curves on $\check V$, the relations
\[
    \sum_{ P_i } n_i L_{ P_i } \simeq 0
\]
hold, where the sum is only taken over the $P_i \in \mathscr S$ that are resolved in $\check V$.
Here $\hat V$ is a suitable small resolution, and $( n_1, \dotsc, n_s ) \in \mathscr B$.

\wernerpage{15}
This implies that if $L_P \simeq 0$ on $\hat V$ for some $P \in \mathscr S$, then $L_P \simeq 0$ on any partial resolution $\check V$ that resolves $P$. \bigskip

Building on the results of this chapter, it is possible to find a necessary and sufficient condition for the existence of projective algebraic small resolutions.
Moreover, we can characterize general projective algebraic small resolutions.
\chapter{Characterization of projective algebraic small resolutions} \label{chapter3}
\wernerpage{16}

The manifolds $\hat V$, $\overline V$, and $\tilde V$ are all Moishezon manifolds: the transcendence degree of their function field is maximal, in our case, equal to three.
According to Moishezon \cite{moisezon67}, for such manifolds the properties \emph{K\"ahler} and \emph{projective algebraic} are equivalent.

The big resolution $\tilde V$ is projective algebraic in any case.  For what reasons might a small resolution not have this property?

One possibility is the existence of an effective nullhomologous curve on $\hat V$.
Such a curve could not exist on a K\"ahler manifold, since the integral of the K\"ahler form on a closed curve must be positive.

For a singular variety $V$ with $d < s$, there is always a small resolution which is not K\"ahler:
one resolves its singularities such that a non-trivial vector in $\mathscr B$ has no negative components; this vector describes an effective nullhomologous curve on some 
small resolution.

As the examples will show later, $d$ is usually very small compared to $s$, and thus the number of non-projective small resolutions $\hat V$ is quite large.
In the case when $d = 0$, every exceptional curve on $\hat V$ is nullhomologous, and there are no projective algebraic small resolutions.

\wernerpage{17}

In fat, this is already the case if even one irreducible exceptional curve on $\hat V$ is nullhomologous, as this means a relation
\[
L_P \simeq 0
\]
holds in $H_2( \hat V, \mathbb Z )$.
As mentioned earlier, this is independent of the particular 
small resolution.

But in the course of this chapter we will show that this is also the only obstruction to the existence of at least one projective algebraic small resolution for a given $V$.
But first we will prove a general theorem characterizing projective algebraic small resolutions. \bigskip

\begin{thm}\label{thm1}
Let $V$ be a nodal complex three-dimensional variety, either a hypersurface in $\mathbb P^4$ or a double solid.
Then for any small resolution $\hat V$, the following are equivalent:
\begin{enumerate}[label=\arabic*.]
    \item $\hat V$ is projective algebraic.
    
    \item There exists a divisor $D$ on $\hat V$ such that
    \[
        D . L_P > 0
    \]
    for all $P \in \mathscr S$.
    
    \item On $\hat V$ there is no nontrivial relation
    \[
        \sum_{i = 1}^s n_i L_{P_i} \simeq 0
    \]
    such that $n_i \geq 0$ for all $i = 1, \ldots, s$.
\end{enumerate}
\end{thm}

\begin{proof}
\wernerpage{18}
``$1 \Rightarrow 2$'' : An ample divisor satisfies condition 2.

``$2 \Rightarrow 3$'' : Let
\[
    m_i := D . L_{ P_i } > 0
\]
be the intersection number of the $i$-th exceptional curve with the given divisor $D$.
Then for all $( n_1, \ldots, n_s ) \in \mathscr B$ we have
\[
    \sum_{i = 1}^s m_i n_i = 0.
\]
Since each $m_i > 0$, there can be no nontrivial element of $\mathscr B$ with $n_i \ge 0$ for all $i = 1, \dotsc, s$, from which property 3 follows immediately.

``$3 \Rightarrow 1$'' : This implication can be proved using a theorem of Peternell (\cite{peternell86}, Theorem 2.5).
In order to apply his results, we must first introduce some new notation. \bigskip

For a compact complex manifold $X$, define
\begin{align*}
\mathscr E^{p, q}( X ) &:= \left\{ \mathscr C^\infty \text{-forms on } X \text{ of type } (p, q) \right\} \\
\mathscr D^{p, q}( X ) &:= \mathscr E^{p, q}( X )'
\end{align*}
Here $\mathscr D^{p, q}( X )$ is the dual space of $\mathscr E^{p, q}( X )$ with respect to the $\mathscr C^\infty$-topology and is called the space of $(p, q)$-currents.
The definition does not agree with the one given in \cite{griffiths78}!
We have $\mathscr D^{n}( X ) = \bigoplus_{p + q = n} \mathscr D^{p, q}( X )$.

A form $\varphi \in \mathscr E^{p, p}( X )$ is called \emph{positive} if for all $x \in X$, the vector $\varphi(x)$ lies in the cone generated by $(i u_1 \wedge \overline u_1 ) \wedge \cdots \wedge (i u_p \wedge \overline u_p )$, where $u_i \in \bigwedge^{1, 0} \, T^*_x(X)$.
A current $T \in \mathscr D^{p, p}( X )$ is called \emph{positive} if $T(\varphi) \geq 0$ for all positive $\varphi \in \mathscr E^{p, p}( X )$.

\wernerpage{19}

Let $\mathscr P_a^1( X )$ be the set of positive $(1,1)$-currents that are weak limits of the form
$\lim_j \sum_{i = 1}^{ n_j } \lambda_{i j} T_{ C_{i j} }$,
where $\lambda_{i j} \in \mathbb{R}$ and $C_{i j}$ are irreducible curves in $X$.
Here $T_{C}( \varphi ) = \int_{C} \varphi$.
Futhermore, define
\[
    d \mathscr D^3( X ) = \left\{ d S : S \in \mathscr D^3( X ) \right\}.
\]

Theorem 2.5 of Peternell's paper now says that a Moishezon threefold $X$ is projective algebraic if there is no current $T \in \mathscr P_a^1( X )$ and no effective curve $C$ with $C + T \simeq 0$.
We want to apply this to our small resolutions. \bigskip

Let $\hat V$ be a small resolution that satisfies condition 3 of our theorem;
thus there are no effective nullhomologous curves on $\hat V$.
Suppose there were an effective curve $C$ and a positive current $T$ such that $C + T \simeq 0$.
Let $\mathscr L$ be an ample line bundle on $V$ and $\pi^* \mathscr L$ the pullback on $\hat V$.
Then
\[
    0 = (C + T) . \pi^* \mathscr L = C . \pi^* \mathscr L + T.\pi^* \mathscr L;
\]
both terms $C . \pi^* \mathscr L$ and $T . \pi^* \mathscr L$ are non-negative, so
\[
    T . \pi^* \mathscr L = 0.
\]
\wernerpage{20}
Hence the support of $T$ must be contained in the fibers of the map $\pi \colon \hat V \to V$, meaning that
\[
    \operatorname{supp} T  \subset \bigcup_{i = 1}^s L_{ P_i }.
\]
A theorem of Siu (\cite{siu74}, Proposition 12.3) yields:
\[
    T = \sum_{i = 1}^s \lambda_i T_{ L_{ P_i } }
\]
for appropriate $\lambda_i \geq 0$.
Thus
\[
C = - \sum_{i = 1}^s \lambda_i T_{ L_{ P_i } },
\]
contradicting the effectiveness of the curve $C$.
Thus the hypothesis of Peternell's theorem is satisfied, so $\hat V$ is projective algebraic.
\end{proof}

The following alternative proof will show how one can achieve the same goal without using Peternell's general theorem.
We do not go directly from 3 to 1, but first prove ``$3 \Rightarrow 2$.''
Transferred to the the vector space $\mathbb Q^s = \mathscr A \perp \mathscr B$, this amounts to
\begin{prop}
If
\[
\mathscr B \cap \mathscr P = \{ 0 \},
\]
where
\[
\mathscr P := \left\{ ( x_1, \ldots, x_s ) \in \mathbb R^s \colon x_i \geq 0 \text{ for all } i = 1, \ldots, s  \right\},
\]
then there is a vector $( m_1, \ldots, m_s ) \in \mathscr A$ such that $m_i > 0$ for all $i = 1, \ldots, s$.
\end{prop}

\begin{proof}
\wernerpage{21}
For technical reasons we must first consider $\mathbb R^s = ( \mathscr A \otimes \mathbb R ) \perp ( \mathscr B \otimes \mathbb R )$.
Furthermore, define
\[
\mathscr K := \left\{ x \in \mathscr P \colon \tfrac12 \leq \sum_{ i = 1 }^s x_i \leq 2 \right\}.
\]
$\mathscr K$ is compact.
By assumption $( \mathscr B \otimes \mathbb R ) \cap \mathscr K = \varnothing$.
The Hahn--Banach theorem gives a hyperplane $\mathscr B_0 \subset \mathbb R^s$ with $\mathscr B_0 \supset \mathscr B \otimes \mathbb R$ and $\mathscr B_0 \cap \mathscr K = \varnothing$.
Let $H_i := \mathbb R e_i$.
For dimension reasons we have
\[ \mathscr B_0 \cap (H_i + H_j) \ne 0 \]
for $i \neq j$.
Moreover, $e_i$ and $e_j$ lie in $\mathscr K$, and thus not in $\mathscr B_0$, so for $i \ne j$ we have
\[
\dim \big( \mathscr B_0 \cap ( H_i + H_j ) \big) = 1.
\]
For each $i = 2, \ldots, s$ there is a vector in $\mathscr B_0 \cap ( H_1 + H_i )$ of the form
\[
y_i := a_i e_1 - e_i \qquad \text{with }a_i > 0.
\]
The vectors $y_2, \ldots, y_s$ are linearly independent and therefore span $\mathscr B_0$.
The vector $( 1, a_2, \ldots, a_s )$ lies in $\mathscr B_0^\perp \subset ( \mathscr B \otimes \mathbb R )^\perp = \mathscr A \otimes \mathbb R$.
This vector can be approximated by vectors in $\mathscr A$, so there is a vector
\[
( \bar m_1, \ldots, \bar m_s ) \in \mathscr A
\]
with $\bar m_i > 0$ for all $i = 1, \ldots, s$.
Clearing denominators, we get a vector $( m_1, \ldots, m_s ) \in \mathscr A$ where $m_i > 0$ and $m_i \in \mathbb Z$ for all $i = 1, \ldots, s$.
Let $\gamma$ be a 4-cycle on $\hat V$ with $k( \gamma ) = ( m_1, \ldots, m_s )$.
By a remark of Clemens in \cite{clemens83},\footnote{Clemens' remark is at the bottom of page 120 and the top of page 121 of his paper.} all of $H_4( \hat V, \mathbb Z )$ is algebraic, so $\gamma$ can be represented by a divisor $D$ satisfying
\[
D . L_P > 0
\]
for all $P \in \mathscr S$.
\end{proof}

\wernerpage{22}

The implication ``$2 \Rightarrow 1$'' can again be proved using a result of Peternell from his article \cite{peternell86}.
This time, instead of the general Theorem 2.5, we require only Corollary 2.3, a generalization of a result of Harvey--Lawson (\cite{harvey83}, Theorem 38). \bigskip

Peternell's result states that a Moishezon manifold $X$ is projective algebraic if and only if the condition
\[\label{stareq}
    \mathscr P_a^1( X ) \cap d \mathscr D^3( X ) = \{ 0 \} \tag{$*$}
\]
is satisfied.
In our situation, we have

\begin{prop}
A small resolution $\hat V$ on which there is a divisor $D$ such that $D . L_P > 0$ for all $P \in \mathscr S$ satisfies condition \eqref{stareq}.
\end{prop}

\begin{proof}
Let $\hat V$ be a small resolution, let $D$ be a divisor such that $D . L_P > 0$ for all $P \in \mathscr S$, and let $T \in \mathscr P_a^1( X ) \cap d \mathscr D^3( X )$.
Then
\[
    T( c_1 (\pi^* \mathscr L) ) = d S( c_1 (\pi^* \mathscr L) )  = S( d c_1 (\pi^* \mathscr L) ) = 0,
\]
where $\mathscr L$ is again an ample line bundle on $V$.
As above, the theorem of Siu implies that there are $\lambda_i \geq 0$ such that
\[
    T = \sum_{i = 1}^s \lambda_i T_{ L_{ P_i } }.
\]
Then
\[
    0 = S( d c_1( [ D ] ) ) = T( c_1( [D] ) ) = \sum_{i = 1}^s \lambda_i D . L_{ P_i }.
\]
Therefore each $\lambda_i = 0$, so $T = 0$.
\end{proof}

\wernerpage{23}

Thus we have also proved the implication ``$2 \Rightarrow 1$'' without using the general theorem of Peternell.
The only disadvantage of this proof is that it fails to make clear how one finds an ample divisor of $\hat V$ given a divisor $D$ that intersects all exceptional curves positively.
The proof of Harvey--Lawson/Peternell is not constructive in this regard, as it relies on the Hahn-Banach theorem. \bigskip

Therefore we give another proof of the implication ``$2 \Rightarrow 1$,'' which shows more precisely how to get an ample divisor on a small resolution $\hat V$ from a divisor $D$ on $\hat V$ such that $D . L_P > 0$ for all $P \in \mathscr S$.

\begin{prop}
Let $\hat V$ be small resolution, and let $D$ be a divisor on $\hat V$ such that $D . L_P > 0$ for all $P \in \mathscr S$.
Furthermore, let $\mathscr L$ be an ample line bundle on $V$, and let $H$ be a divisor on $\hat V$ belonging to $\pi^* \mathscr L$.
Then $D + m H$ is an ample divisor on $\hat V$ for $m \gg 0$.
\end{prop}

\begin{proof}
The line bundle $\mathscr L$ carries a metric whose curvature form is
\[
    \Xi_{ \mathscr L } = \frac{ i }2\, \Omega_{ \mathscr L },
\]
where $\Omega_{ \mathscr L }$ describes a positive definite bilinear form on $T'( V_{ \text{reg} } )$.
If
\[
    \Xi_{ \pi^* \mathscr L } = \frac{ i }2\, \Omega_{ \pi^* \mathscr L }
\]
is the curvature form of the induced metric on the line bundle $\pi^* \mathscr L$ on $\hat V$, then
\[
    \Omega_{ \pi^* \mathscr L } = \pi^* \Omega_{ \mathscr L }.
\]
\wernerpage{24}
Hence $\Omega_{ \pi^* \mathscr L }$ is positive on $\hat V - \bigcup_{P \in \mathscr S} L_P$.
Moreover, for every point $x \in \bigcup_{P \in \mathscr S} L_P$ and ever vector $v \in T_x'( \hat V )$, we have
\[
    \Omega_{ \pi^* \mathscr L }( v , v ) = \Omega_{ \mathscr L }( \pi_* v , \pi_* v ) \geq 0.
\]
Equality holds exactly when $\pi_* v = 0$, meaning that $v$ is tangent to the exceptional curve $L_x$ containing $x$.
On the other hand, by hypothesis we can choose $\Omega_{ [D] }$ so that
\[
    \Omega_{[D]}( v , v ) > 0
\]
for $0 \ne v \in T_x'( L_x )$ and each point $x \in \bigcup_{P \in \mathscr S} L_P$.

Thus for each point $x \in \hat V$ there is an $m_x$ so that $\Omega_{ [D] + m_x \pi^* \mathscr L}$ is positive definite on $T_x'( \hat V )$.
By continuity, this bilinear form remains positive definite on $T_y'( \hat V )$ for $y \in U_x$, where $U_x$ is an affine open neighborhood of $x$.
The $U_x$s form an open cover of $\hat V$.
Since $\hat V$ is compact, there is a finite subcover.
This guaranttes the existence of an $m_1 > 0$ such that $\Omega_{ [D] + m \pi^* \mathscr L }$ is positive definite for all $x \in \hat V$ and all $m \geq m_1$.

Hence $D + m H$ is an ample divisor on $\hat V$ for $m \gg 0$.
\end{proof}

\wernerpage{25}

Now it is easy to state the promised criterion for the existence of projective algebraic small resolutions:

\begin{thm}
Let $V$ be a nodal complex three-dimensional variety, either a hypersurface in $\mathbb P^4$ or a double solid.
Then there exists a projective algebraic small resolution if and only if none of the exceptional curves $L_P$ is nullhomologous in $H_2( \hat V, \mathbb Q )$.
\end{thm}

\begin{proof}
Suppose that no curve $L_P$ is nullhomologous.
Then for each $P_i \in \mathscr S$ there is a vector $( m_1, \ldots, m_s ) \in \mathscr A$ with $m_i \neq 0$; 
otherwise there would be a unit vector with $n_{ i } = 1$ 
in the orthogonal complement of $\mathscr A$, thus in $\mathscr B$, contradicting our assumption.
Then $\mathscr A$ also contains a vector $(m_1, \ldots, m_s)$ such that $m_i \ne 0$ for all $i = 1, \ldots, s$.
Now choose a small resolution such that
\[
m_i > 0 \quad\quad \text{ for all } i = 1, \ldots, s.
\]
On this small resolution there exists a divisor $D$ satisfying
\[
D . L_P > 0 \quad\quad \text{ for all } i = 1, \ldots, s.
\]
Therefore this small resolution, at least, is projective algebraic.
\end{proof}

\chapter{Defect calculations} \label{chapter4}
\wernerpage{26}

In this chapter, we shall describe several ways of checking the criterion we have just proved in concrete examples.
As we have already noted, varieties without defect have no projective algebraic small resolutions.
Thus a method for calculating the defect is certainly useful in this context; we could then immediately exclude uninteresting varieties with $d = 0$ from further investigation. \bigskip

One such method for calculating the defect is given by Clemens in the case of double solids.
Consider the homogeneous polynomials of degree $(3n/2-4)$ in $\mathbb P^3$ that vanish on $\mathscr S$.
Points in general position impose $s$ conditions on these polynomials, so the dimension of the subspace of polynomials vanishing on $\mathscr S$ is only $\binom{3n/2-1}{3} - s$.
But if the singularities are in special position, then there are more polynomials that vanish on $\mathscr S$.

\begin{thm}[Clemens]
Let $V$ be a double solid with branching surface $B \subset \mathbb P^3$, where $\deg B = n \geq 4$.
Then the dimension of the space of homogeneous polynomials of degree $(3n/2-4)$ that vanish on $\mathscr S \subset \mathbb P^3$ is equal to
\[
    \binom{3n/2 - 1}{3} - s + d.
\]
\end{thm}

\wernerpage{27}

Alternatively, one can express this dimension as the rank of a matrix
\[
    \begin{pmatrix}
    f_i( P_j )
    \end{pmatrix}_{\substack{1 \leq j \leq s \\ 1 \leq i \leq t}},
\]
where the $( f_i )$ runs through a basis for the space of homogeneous polynomials of degree $(3n/2 - 4$) and $t = \binom{3n / 2 - 1}{3}$.
In this way, the defect can be calculated as \emph{$s$ minus the rank of the matrix above}.
An analogous method can also be derived for hypersurfaces in $\mathbb P^4$. \pagebreak 

\begin{thm}
Let $V \subset \mathbb P^4$ be a nodal hypersurface with $\deg V = n \geq 3$.
Consider the homogeneous polynomials of degree $(2n - 5)$ that vanish on $\mathscr S$.
These polynomials form a space of dimension
\[
    \binom{2n - 1}{4} - s + d.
\]
\end{thm}

\begin{proof}
Tensor the exact sequence
\[
    0 \longrightarrow \mathscr I_{ \mathscr S } \longrightarrow \mathscr O_{ \mathbb P^4 } \longrightarrow \mathscr O_{ \mathscr S } \longrightarrow 0
\]
with $\Omega_{ \mathbb P^4 }^4( 2V ) \cong \mathscr O_{ \mathbb P^4 }( 2n - 5 )$.
The long exact sequence in cohomology is
\begin{centering}
\[
\setlength{\arraycolsep}{2pt}
\begin{array}{rcccccl}
    0 & \longrightarrow & H^0( \mathbb P^4 ; \mathscr I_{ \mathscr S } \otimes \mathscr O( 2n - 5 ) ) & \longrightarrow & H^0( \mathbb P^4 ; \mathscr O( 2n - 5 ) ) & \longrightarrow & H^0( \mathscr S ; \mathscr O( 2n - 5 ) |_{ \mathscr S } ) \\
    & \longrightarrow & H^1( \mathbb P^4 ; \mathscr I_{ \mathscr S } \otimes \mathscr O( 2n - 5 ) ) & \longrightarrow & H^1( \mathbb P^4 ; \mathscr O( 2n - 5 ) ) & \longrightarrow & \dotsb.
\end{array} \]
\end{centering}

We have
\[
    h^i( \mathscr S ; \mathscr O( 2n - 5 )|_{ \mathscr S } ) =
    \begin{cases}
    s & \text{for } i = 0 \\
    0 & \text{otherwise}
    \end{cases}
\]
and
\begin{align*}
    h^i( \mathbb P^4 ; \mathscr O( 2n - 5 ) ) = 0 && \text{for $i > 0$}.
\end{align*}
\wernerpage{28}
Thus the cohomology sequence above gives
\begin{align*}
    h^i( \mathbb P^4 ; \mathscr I_{ \mathscr S } \otimes \mathscr O( 2n - 5 ) ) = 0
    && \text{for all $i > 1$,}
\end{align*}
and
\[
    h^0( \mathbb P^4 ; \mathscr I_{ \mathscr S } \otimes \mathscr O( 2n - 5 ) ) = \binom{2n - 1}{4} - s + h^1( \mathbb P^4 ; \mathscr I_{ \mathscr S } \otimes \mathscr O( 2n - 5 ) ).
\]
By a result of Schoen (\cite{schoen85}, Proposition 1.3),
\[
    h^1( \mathbb P^4; \mathscr I_{ \mathscr S } \otimes \Omega_{ \mathbb P^4 }^4( 2V ) ) = d.
\]
Now according to Griffiths \cite{griffiths69} we have\footnote{This seems to confuse the ideal sheaf of $\mathscr S$ in $\mathbb P^5$ with the ideal sheaf of $\mathscr S$ in $V$; a calculation with the former should not involve the polynomial $F$ that defines $V$. Nonetheless it is true that $H^0( \mathbb P^4 ; \mathscr I_{ \mathscr S } \otimes \mathscr O( 2n - 5 ) )$ is the space of polynomials of degree $2n-5$ that vanish on $\mathscr S$, and the theorem is valid.}
\[
    H^0( \mathbb P^4 ; \mathscr I_{ \mathscr S } \otimes \mathscr O( 2n - 5 ) ) = \left\{ \frac{f \Omega}{F^2} : \deg f = 2n - 5 \text{ and } f |_{ \mathscr S } \equiv 0 \right\},
\]
where $\Omega = \sum_{i = 0}^{4} ( -1 )^i x_i \, dx_0 \wedge \dotsb \wedge \widehat{dx_i} \wedge \dotsb \wedge dx_4$, $V = \{ F = 0 \}$, and $\deg F = n$, which gives the claim.
\end{proof}

A special case is illustrated by the quadric
\[
    \sum_{i = 1}^{4} x_i^2 = 0
\]
with one ordinary double point.
Here $\beta_3( V_t ) = 0$, so the vanishing cycle is globally nullhomologous on $V_t$.
Hence $d = 1$, and both small resolutions are projective algebraic.
They are given globally by the closures of the graphs of the functions
\[
    - \frac{ x_1 + ix_2 }{ x_3 - ix_4 } = \frac{ x_3 + ix_4 }{ x_1 - ix_2 }
\]
and
\[
    - \frac{ x_1 + ix_2 }{ x_3 + ix_4 } = \frac{ x_3 - ix_4 }{ x_1 - ix_2 }
\]
\wernerpage{29}
This the only example without non-projective small resolutions!

\begin{thm}
Let $V \subset \mathbb P^4$ be a hypersurface with $\deg V \geq 3$, or let $V$ be a double solid branched over a surface $B$ with $\deg B \geq 4$.
Then there exist small resolutions that are not projective algebraic.
\end{thm}

\begin{proof}
Since the singular points impose at least one condition on the homogeneous polynomials of degree $2n - 5$ in $\mathbb P^4$ or of degree $(3n / 2) - 4$ in $\mathbb P^3$, we always have $d < s$.
\end{proof}

By what we have shown above, it is in principle possible to calculate the defect of any given singular hypersurface or double solid, but in practice it is only feasible for values of $n$ that are not too large.
As later examples will show, the condition $d > 0$ alone is not enough to guarantee the existence of a projective algebraic small resolution.
But we have the following

\begin{thm}
Let $V$ be a nodal threefold as above with $d > 0$, and suppose furthermore that there is a group of automorphisms of $V$ which acts transitively on $\mathscr S$.
Then [in any small resolution $\hat V$,] every singularity $P \in \mathscr S$ is resolved by a curve that is not nullhomologous in $H_2( \hat V, \mathbb Q )$.
\end{thm}

\begin{proof}
\wernerpage{30}
Extend the automorphisms of $V$ to automorphisms of the big resolution $\tilde V$, and of $H_2( \tilde V, \mathbb Q )$.
If $L_P \simeq 0$ for some $P \in \mathscr S$, then on $\tilde V$ we would have
\[
    A_P \simeq B_P.
\]
Due to the transitive group operation, this would hold for all $P \in \mathscr S$, contradicting our assumption that $d > 0$.
\end{proof}

Partial and mixed resolutions offer more ways of testing the condition $L_P \simeq 0$.

\begin{thm}
Let $V$ be as above, and let $P \in \mathscr S$.
\begin{enumerate}[label=\alph*)]
    \item If a projective mixed resolution $\overline V$ exists in which the resolution of $P$ is small, then $L_P$ is not nullhomologous on [any small resolution] $\hat V$.
    \item If a projective partial resolution $\check V$ exists in which $P$ is resolved, then $L_P$ is not nullhomologous on [any small resolution] $\hat V$.
    \item If a partial resolution $\check V$ exists with $\beta_2( \check V ) = 1 + d$ that does not resolve $P$, then $L_P$ is nullhomologous on [any small resolution] $\hat V$.
\end{enumerate}
\end{thm}

\begin{proof} \ 
\begin{enumerate}[label=\alph*)]
    \item There are no effective nullhomologous curves on a K\"ahler manifold, so $L_P$ is not nullhomologous on $\overline V$.
    As shown in Chapter \ref{chapter2}, this implies that $L_P$ is not nullhomologous on $\hat V$ either.
    \item There are also no effective nullhomologous curves on a singular projective variety.
    Consider an irreducible curve $C$ on a threefold $X \subset \mathbb P^N$ as a curve in $\mathbb P^N$, $[ C ] \in H_2( \mathbb P^N, \mathbb Z )$.
    \wernerpage{31}
    Let $H$ be a hyperplane in $\mathbb P^N$ in general position with respect to $C$ and $X$, and let $[H] \in H_{2N - 2}( \mathbb P^N, \mathbb Z )$.
    Then we have
    \[
        [H] . [C] > 0.
    \]
    The intersection $X \cap H$ represents a homology class $\gamma \in H_4( X, \mathbb Z )$, and
    \[
        \gamma . [C] = [H] . [C] > 0.
    \]
    This holds in general for every effective curve.
    Hence $L_P$ is not nullhomologous on $\check V$.
    As we saw in Chapter \ref{chapter2}, this implies that $L_P$ is not nullhomologous on $\hat V$.
    
    \item Since $\beta_2( \check V ) = 1 + d$, the second Betti number cannot increase further when we pass to a small resolution $\hat V$.
    Thus all remaining singularities must be replaced by nullhomologous curves. \qedhere
\end{enumerate}
\end{proof}

\begin{remark}
There are two natural ways to get projective partial resolutions. One is to take the closure of the graph of a meromorphic function on $V$ whose indeterminacy set contains $\mathscr S_1 \subset \mathscr S$.
The other is to blow up along a smooth Weil divisor that contains certain singularities of $V$. \bigskip

Suppose we can write the defining equation of $V$ as
\[
    AB + CD = 0,
\]
where $A$, $B$, $C$, and $D$ are smooth.
Let
\[
    \mathscr S_1 := \bigl\{\ P \in V \bigm| 0 = A( P ) = B( P ) = C( P ) = D( P )\  \bigr\} \subset \mathscr S.
\]
\wernerpage{32}
The Cartier divisor $\{ A = 0 \}$ decomposes into two Weil divisors \linebreak 
$\{ A = 0,\, C = 0 \}$ and $\{ A = 0,\, D = 0 \}$, both of which are smooth at all points of $\mathscr S_1$.
Blowing up along one divisor or the other yields the two small resolutions of the singularities in $\mathscr S_1$. \bigskip

In general, every decomposition of the defining equation of $V$ into $AB+CD$ leads to a vector $( m_1, \ldots, m_s ) \in \mathscr A$ such that $m_i \neq 0$ for exactly those $P_i \in \mathscr S$ where $A = 0$, $B = 0$, $C = 0$, and $D = 0$.

It follows immediately that if there is a global Weil divisor on $V$ that is smooth at a point $P \in \mathscr S$, then $L_P$ is not nullhomologous on [any small resolution] $\hat V$.
\end{remark}

\chapter{Chmutov threefolds (1)} \label{chapter5}
\wernerpage{33}

The initial motivation for all our considerations was to study Chmutov hypersurfaces in $\mathbb P^4$ and double solids branched over Chmutov surfaces in $\mathbb P^3$, which together we will refer to as \textit{Chmutov threefolds}.
These classes are to be discussed as first examples; here the theoretical methods developed in the last chapter to prove the existence or non-existence of projective algebraic small resolutions can be applied directly to concrete examples. \bigskip

Let $\mathscr T_n$ denote the Chebyshev polynomial in one variable of degree $n$,
\[
    \mathscr T_n( x ) = \sum_{j = 0}^{\lfloor n / 2 \rfloor} ( -1 )^j \binom{n}{2j} x^{n - 2j} ( 1 - x^2 )^j.
\]
The Chebyshev polynomials satisfy the equation
\[
    \mathscr T_n( \cos \alpha ) = \cos( n \alpha ).
\]
The derivative $\mathscr T_n'( x )$ has simple zeroes at the points $\alpha_k := \cos( k \pi / n )$, where $1 \leq k \leq n - 1$; these points are maxima of $\mathscr T_n$, considered as functions of a \emph{real} variable, if $k$ is even, or minima if $k$ is odd.
The value of $\mathscr T_n$ at these points is
\[
    \mathscr T_n( \alpha_k ) =
    \begin{cases}
    +1 & \text{if $k$ is even} \\
    -1 & \text{if $k$ is odd.}
    \end{cases}
\]


The following lemma will be useful later:

\begin{lem}\label{decomp_lem}
Consider the curves
\[
    C_1 := \{ \mathscr T_n( x ) + \mathscr T_n( y ) = 0 \} \hspace{2em} \text{and} \hspace{2em} C_2 := \{ \mathscr T_n( x ) - \mathscr T_n( y ) = 0 \}
\]
in $\mathbb C^2$.
If $n$ is even, then $C_1$ decomposes into $\frac{n}{2}$ irreducible conics, and $C_2$ decomposes into $\frac{n - 2}{2}$ irreducible conics and two lines.
If $n$ is odd, then $C_1$ and $C_2$ each decompose into $\frac{n - 1}{2}$ irreducible conics and one line.
\end{lem}

\begin{proof}
The curve $C( \mu )$ parametrized by $x = \cos \alpha$ and $y = \cos( \alpha + \mu \pi/n )$, where $0 \le \alpha \le 2\pi$, lies in $C_1$ for $\mu = 1, 3, 5, \ldots, 2 \lfloor(n - 1)/2 \rfloor + 1$, and lies in $C_2$ for $\mu = 0, 2, 4, \ldots, 2 \lfloor n/2 \rfloor$.
Since
\[
    \cos( \alpha + \mu \pi / n ) = \cos( \alpha ) \cos( \mu \pi / n ) - \sin( \alpha ) \sin( \mu \pi / n ),
\]
we can describe $C( \mu )$ by the equation
\[
y = x\,\cos( \mu \pi /n ) - \sqrt{1 - x^2}\,\sin( \mu \pi / 2 ).
\]
This means:
\begin{align*}
    C( 0 ) &= \{ y = x \}
    \\
    C( \mu ) &= \left\{ y^2 + x^2  - 2 [ \cos( \mu \pi / n ) ] x y - [ \sin( \mu \pi / n ) ]^2 = 0 \right\} \quad \text{for } 0 < \mu < n
    \\
    C( n ) &= \{ y = -x \}.
\end{align*}
If $n$ is even then $C_1$ decomposes into the irreducible conics $C( 1 )$, $C( 3 )$, \dots, $C( n - 1 )$ and $C_2$ decomposes into the irreducible conics $C( 2 )$, $C( 4 )$, \dots, $C( n - 2 )$ and the lines $C( 0 )$ and $C( n )$.
If $n$ is odd then $C_1$ decomposes into the conics $C( 1 )$, $C( 3 )$, \dots, $C( n - 2 )$ as well as the line $C( n )$, and $C_2$ decomposes into the conics $C( 2 )$, $C( 4 )$, \dots, $C( n - 1 )$ as well as the line $C( 0 )$.
\end{proof}

\begin{dfn-thm}
\wernerpage{35}
The Chmutov hypersurface of degree $n$ in $\mathbb P^m( \mathbb C )$ is defined in $m$ affine coordinates by the equation
\[
    \sum_{j = 1}^m \mathscr T_n( x_j ) =
    \begin{cases}
    0 & \text{if $m$ is even} \\
    +1 & \text{if $m$ is odd.}
    \end{cases}
\]
The singular points are $( \alpha_{ k_1 }, \ldots, \alpha_{ k_m } )$, where $1 \leq k_i \leq n - 1$ and $\lfloor \frac{m}{2} \rfloor$ of the indices $k_j$ are odd and the rest are even; all singularities are ordinary double points.  Homogenizing the equation above gives no new singularities: all nodes of a projective Chmutov hypersurface lie in the specified affine patch.
Thus in what follows it is enough to calculate in affine coordinates.

More generally, define hypersurfaces by the equations
\[
    \sum_{j = 1}^m ( -1 )^{ \beta_j } \mathscr T_n( x_j ) =
    \begin{cases}
    0 & \text{ if $m$ is even} \\
    ( -1 )^{ \beta_0 } & \text{ if $m$ is odd},
    \end{cases}
\]
where $\beta_0, \beta_1, \ldots, \beta_m \in \{-1, +1\}$. \\
All of these Chmutov hypersurfaces contain many double points, the precise number depending on the choice of the $\beta_i$s. \bigskip
\end{dfn-thm}

We begin with three-dimensional Chmutov hypersurfaces in $\mathbb P^4$.  In the first case, let
\[
    V = \left\{ \sum_{j = 1}^{4} \mathscr T_n( x_j ) = 0\right\}, \quad n \geq 3.
\]
\wernerpage{36}
The singularities lie at the points $( \alpha_{ k_1 }, \alpha_{ k_2 }, \alpha_{ k_3 }, \alpha_{ k_4 })$ where $1 \leq k_i \leq n - 1$ and each point has two even and two odd indices.
The number of nodes is equal to
\[
    6 \bigg( \frac{n}{2} \bigg)^2 \bigg( \frac{n - 2}{2} \bigg)^2 = \frac{3}{8} n^2 ( n - 2 )^2
\]
if $n$ is even and
\[
    6 \bigg( \frac{n - 1}{2} \bigg)^4 = \frac{3}{8} ( n - 1 )^4
\]
if $n$ is odd.

Let $P \in \mathscr S$.
Without loss of generality, we may assume $P$ is of the form $( \alpha_{ u_1 }, \alpha_{ u_2 }, \alpha_{ g_1 }, \alpha_{ g_2 } )$, where $u_1$ and $u_2$ are odd and $g_1$ and $g_2$ are even.\footnote{In English the letters $u$ and $g$ may seem like strange choices, but they stand for \emph{ungerade} (odd) and \emph{gerade} (even).}
Both $\mathscr T_n( x_1 ) + \mathscr T_n( x_3 )$ and $\mathscr T_n( x_2 ) + \mathscr T_n( x_4 )$ vanish at $P$.
By Lemma \ref{decomp_lem} we can factor 
\begin{align*}
    \mathscr T_n( x_1 ) + \mathscr T_n( x_3 ) &= \prod_{i} C_i
    \\
    \mathscr T_n( x_2 ) + \mathscr T_n( x_4 ) &= \prod_{j} D_j,
\end{align*}
where $C_1, C_2, D_1,$ and $D_2$ all vanish at $P$.
The defining equation of $V$ can then be rearranged as
\[
    \mathscr T_n( x_1 ) + \mathscr T_n( x_3 ) = -[ \mathscr T_n( x_2 ) + \mathscr T_n( x_4 )].
\]
The global meromorphic function
\[
    \frac{ C_1 }{ D_1 } = -\, \frac{ \prod_{j > 1} D_j }{ \prod_{i > 1} C_i }
\]
has a point of indeterminacy at $P$; the closure of the graph of this function describes a projective partial resolution which resolves $P$.
Hence the exceptional curve $L_P$ is not nullhomologous on [any small resolution] $\hat V$.
For each $P \in \mathscr S$, such a meromorphic function with a point of indeterminacy at $P$ can be found, thus ensuring the existence of a projective algebraic small resolution.

\wernerpage{37}

The same argument can be used in the general case
\[
    V = \left\{ \sum_{j = 4}^{m} ( -1 )^{ \beta_j } \mathscr T_n( x_j ) = 0 \right\} \hspace{1em} ( n \geq 3 ),
\]
where $\beta_j \in \{+1, -1\}$ for all $j = 1, \ldots, 4$.
Again, all singularities are of the form $( \alpha_{ k_1 }, \alpha_{ k_2 }, \alpha_{ k_3 }, \alpha_{ k_4 } )$ for $1 \leq k_i \leq n - 1$; the distribution of even and odd indices must be such that for each point $P \in \mathscr S$, two of the terms $( -1 )^{ \beta_j } \mathscr T_n( x_j )$ take the a value $+1$ and two take the value $-1$. \bigskip

Let $P \in \mathscr S$, and assume without loss of generality that at this point,
\[
    ( -1 )^{ \beta_1 } \mathscr T_n( x_1 ) = +1 = ( -1 )^{ \beta_2 } \mathscr T_n( x_2 ).
\]
Then the polynomials
\begin{align*}
    ( -1 )^{ \beta_1 } \mathscr T_n( x_1 ) + ( -1 )^{ \beta_3 } \mathscr T_n( x_3 )
    \hspace{2em} \text{and} \hspace{2em}
    ( -1 )^{ \beta_2 } \mathscr T_n( x_2 ) + ( -1 )^{ \beta_4 } \mathscr T_n( x_4 )
\end{align*}
vanish at $P$.
By Lemma \ref{decomp_lem}, both factor as products of non-singular polynomials of degrees one and two, of which two vanish at $P$.
As in the first case, we get a global meromorphic function with a point of indeterminacy at $P$.

Thus for any choice of $\beta_j$ and all $n \geq 3$, the existence of a projective algebraic small resolution is guaranteed. \bigskip

\wernerpage{38}

The situation with the Chmutov double solids is more interesting.
Here there are four possible choices of sign up to isomorphism.
Unlike with Chmutov hypersurfaces, with Chmutov double solids the existence of projective algebraic small resolutions depends on the particular 
choice of $\beta_i$s. \bigskip

We begin with the family of double solids branched over the surface
\[
    B = \left\{ \mathscr T_n( x_1 ) + \mathscr T_n( x_2 ) + \mathscr T_n( x_3 ) + 1 = 0 \right\},
\]
where $n \geq 4$ is an even integer.
The singular points are $( \alpha_{ k_1 }, \alpha_{ k_2 }, \alpha_{ k_3 } )$, where $1 \leq k_i \leq n - 1$ and exactly one of the indices $k_i$ is even.
The nodes of $V$ are identified with those of $B$, and there are $\frac{3}{8} n^2 (n - 2)$ of them.
Let $P \in \mathscr S$, and assume without loss of generality that $k_1$ is even for $P$.
We write the equation of the double solid $V$ as
\[
    \mathscr T_n( x_1 ) + \mathscr T_n( x_2 ) = \omega^2 - [ \mathscr T_n( x_3 ) + 1 ].
\]
The left-hand side factors into the product of $n / 2$ non-singular polynomials $C_1, \ldots, C_{n / 2}$ of degree 2, two of which vanish at $P$.
Let these be $C_1$ and $C_2$.
The polynomial $\mathscr T_n( x ) + 1$ has $\frac{n}{2}$ double roots, so there is a polynomial $F_n( x )$ of degree $\frac{n}{2}$ with
\[
    \mathscr T_n( x ) + 1 = F_n( x )^2.
\]
Thus the right-hand side of the equation above factors as
\[
    [ \omega + F_n( x_3 ) ] [ \omega - F_n( x_3 ) ].
\]
\wernerpage{39}
The meromorphic function
\[
    \frac{ C_1 }{ \omega + F_n( x_3 )} = \frac{ \omega - F_n( x_3 )}{\prod_{i = 2}^{n / 2} C_i}
\]
on $V$ has a point of indeterminacy at $P$.
Thus in these examples as well, for all $n$ there exists a projective algebraic small resolution. \bigskip

Now let
\[
    B = \left\{ \mathscr T_n( x_1 ) + \mathscr T_n( x_2 ) - \mathscr T_n( x_3 ) + 1 = 0 \right\}.
\]
In this case the singularities lie at the points $( \alpha_{ u_1 }, \alpha_{ u_2 }, \alpha_{ u_3 } )$, $( \alpha_{ u_4 }, \alpha_{ g_1 }, \alpha_{ g_2 } )$, and $( \alpha_{ g_3 }, \alpha_{ u_5 }, \alpha_{ g_4 } )$, where $1 \leq u_i, g_i \leq n - 1$, the $g_i$ are even, and the $u_i$ are odd.
Let $P \in \mathscr S$; then either the index of the first or the second coordinate is odd.
Without loss of generality, assume that the second coordinate has an odd index.
Write the equation of $V$ as
\[
    \mathscr T_n( x_2 ) - \mathscr T_n( x_3 ) = \omega^2 - [ \mathscr T_n( x_1 ) + 1 ].
\]
As in the first case, both sides factor into two polynomials, each of which vanishes at $P$.
Again there exists a meromorphic function with a point of indeterminacy at $P$, and in each case there are projective algebraic small resolutions. \bigskip

\wernerpage{40}

This does not hold for the double solid branched over
\[
    B = \{ \mathscr T_n( x_1 ) - \mathscr T_n( x_2 ) - \mathscr T_n( x_3 ) + 1 = 0 \},
\]
whose singularities lie at the points $( \alpha_{ u_1 }, \alpha_{ g_1 }, \alpha_{ u_2 } )$, $( \alpha_{ u_3 }, \alpha_{ u_4 }, \alpha_{ g_2 } )$, and \linebreak $( \alpha_{ g_3 }, \alpha_{ g_4 }, \alpha_{ g_5 } )$, where $1 \leq u_i, g_i \leq n - 1$, the $g_i$ are even, and the $u_i$ are odd.
The first two types of singularities can be resolved by curves that are not nullhomologous: if we write the equation for the double solid as
\[
    - [ \mathscr T_n( x_2 ) + \mathscr T_n( x_3 ) ] = \omega^2 - [ \mathscr T_n( x_1 ) + 1 ].
\]
then the closure of the graph of the meromorphic function
\[
    \frac{ \prod_{i = 1}^{j} C_i }{ \omega + F_n( x_1 ) } = \frac{ \omega - F_n( x_1 ) }{ \prod_{i = j + 1}^{n / 2} C_i },
\]
where $j = 1,\, \ldots,\, n/2 - 1$, provides a projective algebraic partial resolution $\check V$ with $\beta_2( \check V ) \geq n/2$, which only contains singularities of the third type.
As will be shown in the next chapter, at least for $n \in \{ 4, 6, 8 \}$ the defect is equal to $n/2 - 1$.
Thus in these cases,
\[
    \beta_2( \check V ) = \beta_2( \hat V ) = n/2,
\]
so the singularities of the third type are all resolved by nullhomologous curves.
Despite the positive defect, there are no projective algebraic small resolutions.

\begin{conj}
This holds for all even integers $n \geq 10$.
\end{conj}

\wernerpage{41}

To close this chapter, we briefly consider the double solid branched over
\[
    B = \{ \mathscr T_n( x_1 ) + \mathscr T_n( x_2 ) + \mathscr T_n( x_3 ) - 1 = 0 \}.
\]
The number of singularities is only $\frac{3}{8} n (n - 2)^2$.
If $n \in \{4, 6, 8\}$, then the defect is 0, as we will show in the next chapter; all nodes are resolved by nullhomologous curves.

\begin{conj}
For all even integers $n \geq 10$, the defect of the double solid branched over $B$ vanishes.
\end{conj}

\chapter{Chmutov threefolds (2)} \label{chapter6}
\wernerpage{42}

Having focused on general properties of the various families of Chmutov threefolds in the previous chapter, we shall now examine particular 
examples in greater detail.
This chapter is mainly about calculating the defect.
In several cases we need to know the defect precisely in order to prove the non-existence of projective algebraic small resolutions. \bigskip

As we found in chapter \ref{chapter4}, the defect of a hypersurface (or a double solid) can be calculated as
\[
    s - \operatorname{rank}\big( f_i ( P_j ) \big),
\]
where $f_1, \ldots, f_N$ runs through a basis for the homogeneous polynomials of degree $(2n - 5)$ in $\mathbb P^4$ (or of degree $(\frac{3n}{2} - 4)$ in $\mathbb P^3$) and $P_1, \ldots, P_s$ denote the singularities of $V \subset \mathbb P^4$ (or of $B \subset \mathbb P^3$). \bigskip

When $\deg( V ) = 3$ and $\deg( B ) = 4$, these calculations can be carried out without any problems.

The singularities of the Chmutov cubic
\[
    V = \Bigg\{ \sum_{i = 1}^{4} \mathscr T_3( x_i ) = 0 \Bigg\}
\]
lie at the points
\begin{align*}
    P_{--++} &:= \left( -\tfrac{1}{2}, -\tfrac{1}{2}, +\tfrac{1}{2}, +\tfrac{1}{2} \right),
    \\
    P_{-+-+} &:= \left( -\tfrac{1}{2}, +\tfrac{1}{2}, -\tfrac{1}{2}, +\tfrac{1}{2} \right)
\end{align*}
and so on.

\wernerpage{43}

There are 6 different ways to arrange two plus signs and two minus signs, so $s = 6$.
All singularities lie in the hyperplane $\{ \sum_{i = 1}^{4} x_i = 0 \}$, hence $d = 2$.
With the help of partial resolutions, one determines the linear relations between the exceptional curves (up to signs) as
\begin{align*}
    L_{++--} &\simeq L_{--++},
    \\
    L_{+-+-} &\simeq L_{-+-+},
    \\
    L_{+--+} &\simeq L_{-++-}
    \\
    \text{and} \hspace{2em}
    L_{++--} - L_{+-+-} &- L_{+--+} \simeq 0.
\end{align*}
There are 6 projective algebraic small resolutions.

The general case
\[
    V = \Bigg\{ \sum_{j = 1}^{4} ( -1 )^{ \beta_j } \mathscr T_3( x_j ) = 0 \Bigg\}
\]
yields nothing new; all of these hypersurfaces are isomorphic, as $\mathscr T_n( x ) = -\mathscr T_n( -x )$ for all odd natural numbers $n$. \bigskip

In the case of double solids branched over surfaces of degree 4, one arrives at the following results:

\begin{enumerate}[label=\arabic*.]
    \item Let
    \[
    B = \big\{ \mathscr T_4( x_1 ) + \mathscr T_4( x_2 ) + \mathscr T_4( x_3 ) + 1 = 0 \big\}.
    \]
\wernerpage{44}
    The singular points are
    \begin{align*}
        P_{\pm \pm 0} &:= \left( \pm \tfrac{1}{\sqrt{2}},\, \pm \tfrac{1}{\sqrt{2}},\, 0 \right),
        \\
        P_{\pm 0 \pm} &:= \left( \pm \tfrac{1}{\sqrt{2}},\, 0,\, \pm \tfrac{1}{\sqrt{2}}\right),
        \\
        \text{and} \hspace{2em} P_{0 \pm \pm} &:= \left( 0,\, \pm \tfrac{1}{\sqrt{2}},\, \pm \tfrac{1}{\sqrt{2}} \right).
    \end{align*}
    We have $s = 12$.
    Since all singularities lie on the quadric \linebreak $\{ \sum_{i = 1}^{3} x_{i}^{2} = 1 \}$, we have $d = 3$.
    Using this quadric, we can write the defining equation of $B$ as
    \[
        \left( \sum_{i = 1}^{3} x_i^2 - 1 \right)^2 + ( x_1 - x_2 - x_3 ) ( x_1 + x_2 - x_3 ) ( x_1 - x_2 + x_3 ) ( x_1 + x_2 + x_3 ) = 0.
    \]
    The four planes $\{ x_1 \pm x_2 \pm x_3 = 0 \}$ intersect pairwise in 6 lines, the intersection of these 6 lines with the quadric $\{ \sum_{i = 1}^{3} x_i^2 = 1 \}$ being exactly the 12 singularities of $B$. \bigskip
    
Quartics of the form $\{ q^2 + \prod_{i = 1}^{4} \ell_{i} = 0 \}$, $\deg q = 2$, 
$\deg \ell_i = 2$, and their double solids are generally examined in chapter \ref{chapter10}.
    Writing the quartic in this form gives new global meromorphic functions on $V$:\footnote{Here $\omega$ comes from writing the double solid as $\omega^2 = $ quartic.}
        \[
        \frac{ \omega - q }{ \prod_{i \in I} \ell_i } = \frac{ \prod_{i \notin I} \ell_i }{ \omega + q },
    \]
where $I$ is a nonempty subset of $\{ 1, 2, 3, 4 \}$.
    The locus of indeterminacy is contained in $\mathscr S$, and the graphs of these functions describe certain 
    partial resolutions.
    
\wernerpage{45}
    
    The linear relations between the exceptional curves read as follows on a particular 
    small resolution:
    \begin{align*}
        L_{+0+} &\simeq L_{-0-},
        \\
        L_{+0-} &\simeq L_{-0+},
        \\
        L_{++0} &\simeq L_{--0},
        \\
        L_{+-0} &\simeq L_{-+0},
        \\
        L_{0++} &\simeq L_{0--},
        \\
        L_{0+-} &\simeq L_{0-+},
        \\
        L_{+0+} - L_{0++} & + L_{-+0} \simeq 0,
        \\
        L_{+0+} - L_{0+-} & - L_{++0} \simeq 0
       \\
       \text{and} \hspace{2em} L_{+0-} - L_{0++} & + L_{++0} \simeq 0.
    \end{align*}
    There are 24 projective algebraic small resolutions.
    
    \item Things are similar in the case
    \[
        B = \{ \mathscr T_4( x_1 ) + \mathscr T_4( x_2 ) - \mathscr T_4( x_3 ) + 1 = 0 \}.
    \]
    The singular points are
    \begin{align*}
        P_{\pm \pm \pm} &:= \left( \pm \tfrac{1}{\sqrt{2}},\, \pm \tfrac{1}{\sqrt{2}},\, \pm \tfrac{1}{\sqrt{2}} \right),
        \\
        P_{\pm 0 0} &:= \left( \pm \tfrac{1}{\sqrt{2}},\, 0,\, 0 \right)
        \\
        \text{and} \hspace{2em} P_{0 \pm 0} &:= \left( 0,\, \pm \tfrac{1}{\sqrt{2}},\, 0 \right).
    \end{align*}
    We have $s = 12$.
    All singularities lies on the quadric
    \[
        \{ 2x_1^2 + 2x_2^2 - 2x_3^2 = 1 \},
    \]
    so $d = 3$.
    
\wernerpage{46}
    
    The defining equation of $B$ can be written as
    \[
        ( 2x_1^2 + 2x_2^2 + 2x_3^2 - 1)^2 - 8( x_1 - x_3 ) ( x_1 + x_3 ) ( x_2 - x_3 ) ( x_2 + x_3 ) = 0,
    \]
    which is again a special case of an example that will be discussed in more detail in chapter \ref{chapter10}.
    The linear relations are
    \begin{align*}
        L_{+++} &\simeq L_{---},
        \\
        L_{+-+} &\simeq L_{-+-},
        \\
        L_{+--} &\simeq L_{-++},
        \\
        L_{++-} &\simeq L_{--+},
        \\
        L_{+00} &\simeq L_{-00},
        \\
        L_{0+0} &\simeq L_{0-0},
        \\
        L_{+++} - L_{+-+} & - L_{+00} \simeq 0,
        \\
        L_{+++} - L_{0+0} & - L_{+--} \simeq 0
        \\
        \text{and} \hspace{2em} L_{+-+} - L_{0+0} & - L_{++-} \simeq 0.
    \end{align*}
    Exactly 24 small resolutions are projective algebraic.
    
    \item Especially interesting is the case
    \[
        B = \{ \mathscr T_4( x_1 ) - \mathscr T_4( x_2 ) - \mathscr T_4( x_3 ) + 1 = 0 \}
    \]
    with $s = 9$; the singular points are $(0, 0, 0)$, $( \pm \frac{1}{\sqrt{2}}, 0, \pm \frac{1}{\sqrt{2}} )$ and \linebreak $(\pm \frac{1}{\sqrt{2}}, \pm \frac{1}{\sqrt{2}}, 0 )$.
    The singular set $\mathscr S$ lies on the two quadrics $\{ x_1^2 - x_2^2 - x_3^2 = 0 \}$ and $\{ x_2 x_3 = 0 \}$ ; observe that no other quadrics contain $\mathscr S$.
    Hence $d = 1$.
    As we saw in the last chapter, there exists a partial resolution $\check V \subset \mathbb P^N$ of the last $8$ singularities.
    The ninth double point is thus resolved nullhomologously.
    
\wernerpage{47}
    
    \item In the case
    \[
        B = \{ \mathscr T_4( x_1 ) + \mathscr T_4( x_2 ) + \mathscr T_4( x_3 ) - 1 = 0 \},
    \]
    the singular points are $( 0, 0 \pm \frac{1}{\sqrt{2}} )$, $(0, \pm \frac{1}{\sqrt{2}}, 0 )$ and $( \pm \frac{1}{\sqrt{2}}, 0, 0 )$. \linebreak 
    Precisely the quadrics $\{ x_1 x_2 = 0 \}$, $\{ x_1 x_3 = 0 \}$, $\{ x_2 x_3 = 0 \}$, and $\{ x_1^2 + x_2^2 + x_3^2 = \frac{1}{2} \}$ vanish on all of $\mathscr S$, so $d = 0$.
    All nodes are resolved by nullhomologous curves.
    
\end{enumerate}

To obtain a positive defect, the singularities must all lie on a plane conic.
Later we will encounter such an example with $\deg B = 4$, $s = 6$, and $d = 1$.
Thus it is clear that the defect really depends on the specific position of the singularities. \bigskip

When $\deg V = 4, 5$ and $\deg B = 6, 8$, the help of a computer is needed to calculate the defect.
The coordinates of the singularities and thus the entries of the matrix $(f_i( P_j ) )$ each lie in some $\mathbb Z[ \alpha ]$, where $\alpha$ is integral over $\mathbb Z$.
Except for one special case $(\deg B = 8)$, $\alpha$ is always the square root of a prime number.

For technical reasons it was necessary in the computer calculations to reduce all matrix entries modulo a prime.
If this prime is inert --- meaning that $p \mathbb Z[ \alpha ]$ is a prime ideal in $\mathbb Z[ \alpha ]$ --- then $\mathbb Z[ \alpha ] / p\mathbb Z[ \alpha ]$ has no zero-divisors.
Thus the rank of the matrix $(f_{i}( P_j ) \text{ mod } p )$ is less than or equal to the rank of the matrix $( f_i( P_j ) )$.
Define
\[
    d'( p ) := s - \text{Rank}\big( f_i( P_j) \text{ mod } p \big).
\]
Then $d \leq d'( p )$ for all inert primes $p$.

\wernerpage{48}
    
The defect can be bounded below using divisors on $V$, so its exact value can be given in all the cases treated here. \bigskip
    
How does one get suitable primes for reduction?
From number theory, we know that a prime number $p$ is inert in $\mathbb Z[ \alpha ]$ exactly when the minimal polynomial of $\alpha$ over $\mathbb Z$ is also irreducible over $\mathbb Z_p$.
If $\alpha = \sqrt{q}$ for a prime $q$, this means that if $x^2-q$ has no roots in $\mathbb Z_p$, then $p$ is a suitable prime for reduction. \bigskip

Consider the Chmutov quartics in $\mathbb P^4$ whose singularities have coordinates lying in $\mathbb Z[ \sqrt{2} ]$ after homogenization and suitable normalization.
The hypersurface
\[
    V = \{ \mathscr T_4( x_1 ) + \mathscr T_4( x_2 ) + \mathscr T_4( x_3 ) + \mathscr T_4( x_4 ) = 0 \}
\]
has 24 singularities; we divide them into three classes:
\[ \begin{array}{llllll}
    \text{Type I} & : & & ( \pm, 0, \pm, 0 ) & \text{ and } & ( 0, \pm, 0, \pm ) \\
    \text{Type II} & : & & ( \pm, 0, 0, \pm ) & \text{ and } & ( 0, \pm, \pm, 0 ) \\
    \text{Type III} & : & & ( \pm, \pm, 0, 0 ) & \text{ and } & ( 0, 0, \pm, \pm ). \\
\end{array} \]
Here ``$+$'' means $1 / \sqrt{2}$ and ``$-$'' means $-1 / \sqrt{2}$.

\wernerpage{49}

With the computer we calculated that $d'( 181 ) = 2$.
The global divisors
\[
    \left\{ x_1^2 + x_2^2 + x_1 x_2 - \tfrac{1}{2} = 0,\ x_3^2 + x_4^2 + x_3 x_4 - \tfrac{1}{2} = 0 \right\}
\]
and
\[
    \left\{ x_1^2 + x_3^2 + x_1 x_3 - \tfrac{1}{2} = 0,\ x_2^2 + x_4^2 + x_2 x_4 - \tfrac{1}{2} = 0 \right\}
\]
yield two linearly independent vectors in $\mathscr A$, so we have
\[
    d = 2.
\]
On a particular 
small resolution, we have the following relations between the irreducible exceptional curves: all curves over singularities of the same type are homologous, and furthermore
\[
    L_\text{I} - L_\text{II} + L_\text{III} \simeq 0.
\]
There are 6 projective algebraic small resolutions. \bigskip

The quartic
\[
    V = \{ \mathscr T_4( x_1 ) + \mathscr T_4( x_2 ) - \mathscr T_4( x_3 ) - \mathscr T_4( x_4 ) = 0 \}
\]
has 33 nodes, abbreviated as $( 0, 0, 0, 0 )$, $( \pm, \pm, \pm, \pm )$, $( 0, \pm, \pm, 0 )$, $( \pm, 0, \pm, 0 )$, $( 0, \pm, 0, \pm )$ and $( \pm, 0, 0, \pm )$.
Again, a ``$+$'' stands for $1 / \sqrt{2}$ and a ``$-$'' for $- 1 / \sqrt{2}$.

The computer calculates that $d'( 181 ) = 7$.

The quartic contains the 8 planes $\{ x_1 \pm x_i = 0, x_j \pm x_k = 0 \}$, $(i, j, k) \in \{ (3, 2, 4), (4, 2, 3) \}$, each of which contains 9 of the nodes.
The vectors of the intersection numbers of the proper transformations of these planes with the exceptional curves on a small resolution span a subspace of dimension greater than or equal to 7.

\wernerpage{50}

To see this, choose 7 planes $E_i$, one after another, so that each $E_i$ contains a node $P_i$ not contained in the the previously chosen planes, for example:
\begin{align*}
    E_1 &= \{ x_1 - x_3 = 0, x_2 - x_4 = 0 \} & P_1 &= ( +, +, +, + ) \\
    E_2 &= \{ x_1 + x_3 = 0, x_2 - x_4 = 0 \} & P_1 &= ( +, -, -, - ) \\
    E_3 &= \{ x_1 - x_3 = 0, x_2 + x_4 = 0 \} & P_1 &= ( +, +, +, - ) \\
    E_4 &= \{ x_1 + x_3 = 0, x_2 + x_4 = 0 \} & P_1 &= ( -, +, -, + ) \\
    E_5 &= \{ x_1 + x_4 = 0, x_2 + x_3 = 0 \} & P_1 &= ( 0, +, -, 0 ) \\
    E_6 &= \{ x_1 + x_4 = 0, x_2 - x_3 = 0 \} & P_1 &= ( 0, +, +, 0 ) \\
    E_7 &= \{ x_1 - x_4 = 0, x_2 - x_3 = 0 \} & P_1 &= ( +, 0, 0, + ) \\
\end{align*}
Thus the defect $d = 7$.
The decomposition of $\mathbb Q^s$ into $\mathscr A$ and $\mathscr B$, and thus also the number of projective algebraic small resolutions, can be determined in theory, although it would involve considerable computing time. \bigskip

Let us turn instead to the quartic
\[
    V = \{ \mathscr T_4( x_1 ) + \mathscr T_4( x_2 ) + \mathscr T_4( x_3 ) - \mathscr T_4( x_4 ) = 0 \},
\]
with $s = 30$.
The computer calculates that $d'( 181 ) = 8$.
Write the equation of the quartic as
\[
    \mathscr T_4( x_1 ) + \mathscr T_4( x_2 ) 
    = \mathscr T_4( x_4 ) - \mathscr T_4( x_3 ).
\]
\wernerpage{51}
The term on the left-hand side of the equation factors
\[
    2 ( 2x_1^2 + 2x_2^2 + 2x_1 x_2 - 1 ) ( 2x_1^2 + 2x_2^2 - 2x_1 x_2 - 1 ),
\]
and the term on the right-hand side as
\[
    ( x_4 + x_3 ) ( x_4 - x_3 ) ( x_3^2 + x_4^2 - 1 ).
\]
We see that the divisors
\begin{align*}
    D_1 &:= \{ x_3 + x_4 = 0,\ 2x_1^2 + 2x_2^2 + 2x_1 x_2 - 1 = 0 \}
    \\
    \text{and} \qquad
    D_2 &:= \{ x_3 - x_4 = 0,\ 2x_1^2 + 2x_2^2 + 2x_1 x_2 - 1 = 0 \}
\end{align*}
are smooth Weil divisors on $V$.
By transposing $x_3$ with $x_1$ or $x_2$ one obtains more divisors
\begin{align*}
    D_3 &:= \{ x_2 + x_4 = 0, 2x_1^2 + 2x_3^2 + 2x_1 x_3 - 1 = 0 \}
    \\
    D_4 &:= \{ x_2 - x_4 = 0, 2x_1^2 + 2x_3^2 + 2x_1 x_3 - 1 = 0 \}
    \\
    \text{and} \qquad
    D_5 &:= \{ x_1 + x_4 = 0, 2x_2^2 + 2x_3^2 + 2x_2 x_3 - 1 = 0 \}.
\end{align*}
Additional divisors can be seen by rewriting the equation of the quartic as
\[
    \mathscr T_4( x_1 ) + \mathscr T_4( x_2 ) + \mathscr T_4( x_3 ) + 1
    = \mathscr T_4( x_4 ) + 1.
\]
As we have already seen with double solids, the term on the right-hand side of the equation equals $q_1^2 := 2 ( 2x_4^2 - 1 )^2$ and the term on the left-hand side equals $q_2^2 + 4 \prod_{i = 1}^4 \ell_i$, where $\ell_i$ runs through the 4 linear polynomials $x_1 \pm x_2 \pm x_3$ and $q_2$ is defined as $2 \big( \sum_{i = 1}^3 x_i^2 - 1 \big)$.
The equation of the quartic can be transformed into
\[
    \prod_{i = 1}^4 \ell_i = q_1^2 - q_2^2;
\]
\wernerpage{52}
the divisors
\begin{align*}
    D_6 &:= \{ x_1 + x_2 + x_3 = 0,\ q_1 + q_2 = 0 \}
    \\
    D_7 &:= \{ x_1 + x_2 - x_3 = 0,\ q_1 + q_2 = 0 \}
    \\
    \text{and} \qquad
    D_8 &:= \{ x_1 - x_2 + x_3 = 0,\ q_1 + q_2 = 0 \}.
\end{align*}
are again smooth Weil divisors on $V$.
The 33-tuples $(m_1, \ldots, m_{33}) \in \mathscr A$ associated to the divisors $D_1$ through $D_8$ are linearly independent.
I would rather not present the computationally laborious proof here; the method of proof will be explained in more detail later with a more important example (a double solid branched over a Chmutov octic with 144 nodes). \bigskip

We now come to the Chmutov quintic
\[
    V = \bigg\{ \sum_{i = 1}^4 \mathscr T_5( x_i ) \bigg\}
\]
with 96 nodes.
The coordinates of its singularities lie in $\mathbb Z[ \sqrt{5} ]$ after suitable normalization, and the computer calculates that
\[
    d'( 173 ) = 10.
\]
To also bound the defect from below, we should now again construct a sequence $D_i$ of 10 global Weil divisors, such that each divisor $D_i$ contains at least one point $P_i \in \mathscr S$ with $P_i \notin D_j$ for all $j < i$.

\wernerpage{53}

As we know, $\mathscr T_5( x ) + \mathscr T_5( y )$ decomposes into 3 polynomials of degrees 1, 2, and 2.
The factorization is
\[
    \mathscr T_5( x ) + \mathscr T_5( y ) = ( x + y ) f( x, y ) g( x, y ),
\]
\begin{align*}
    &\text{with} & f( x, y ) &= 4x^2 + 4y^2 + 2( \sqrt{5} - 1 ) x y - \tfrac{( 5 + \sqrt{5} )}{2}
    \\
    &\text{and} &
    g( x, y ) &= 4x^2 + 4y^2 - 2( \sqrt{5} - 1 ) x y - \tfrac{( 5 + \sqrt{5} )}{2}.
\end{align*}

We now choose:
\begin{align*}
    D_1 &= \{ x_1 + x_2 = 0, x_3 + x_4 = 0 \}, & P_1 &= ( \alpha_1, \alpha_4, \alpha_1, \alpha_4 ) ;
    \\
    D_2 &= \{ x_1 + x_3 = 0, x_2 + x_4 = 0 \}, & P_2 &= ( \alpha_1, \alpha_1, \alpha_4, \alpha_4 ) ;
    \\
    D_3 &= \{ x_1 + x_4 = 0, x_2 + x_3 = 0 \}, & P_3 &= ( \alpha_1, \alpha_2, \alpha_3, \alpha_4 ) ;
    \\
    D_4 &= \{ x_1 + x_2 = 0, g( x_3, x_4 ) = 0 \}, & P_4 &= ( \alpha_1, \alpha_4, \alpha_1, \alpha_2 ) ;
    \\
    D_5 &= \{ x_1 + x_3 = 0, g( x_2, x_4 ) = 0 \}, & P_5 &= ( \alpha_1, \alpha_1, \alpha_4, \alpha_2 ) ;
    \\
    D_6 &= \{ x_1 + x_4 = 0, g( x_2, x_3 ) = 0 \}, & P_6 &= ( \alpha_1, \alpha_1, \alpha_2, \alpha_4 ) ;
    \\
    D_7 &= \{ f( x_1, x_2 ) = 0, g( x_3, x_4 ) = 0 \}, & P_7 &= ( \alpha_1, \alpha_2, \alpha_1, \alpha_2 ) ;
    \\
    D_8 &= \{ f( x_1, x_3 ) = 0, g( x_2, x_4 ) = 0 \}, & P_8 &= ( \alpha_1, \alpha_1, \alpha_2, \alpha_2 ) ;
    \\
    D_9 &= \{ f( x_1, x_2 ) = 0, x_3 + x_4 = 0 \}, & P_9 &= ( \alpha_2, \alpha_1, \alpha_4, \alpha_1 ) ;
    \\
    D_{10} &= \{ f( x_1, x_3 ) = 0, x_2 + x_4 = 0 \}, & P_{10} &= ( \alpha_2, \alpha_4, \alpha_1, \alpha_1 ) .
\end{align*}
Here we recall that
\begin{align*}
    \alpha_1 &= \frac{\sqrt{5} + 1}{4} = - \alpha_4, \\
    \alpha_2 &= \frac{\sqrt{5} - 1}{4} = -\alpha_3.
\end{align*}
Hence the vector space $\mathscr A$ has at least dimension 10, and the equation
\[ d = 10 \]
is proved!

\wernerpage{54}

The projective algebraic small resolutions are K\"ahler manifolds with trivial canonical bundle and Euler characteristic $-8$.
The interesting Hodge numbers are $h^{1, 1}( \hat V ) = 11$ and $h^{2, 1}( \hat V ) = 15$. \bigskip

We come now to the remaining double solids.
The coordinates of the singularities of the double solids branched over Chmutov sextics lie in in $\mathbb Z[ \sqrt{3} ]$, for which 173 is a suitable prime for reduction.
With the octics, the coordinates of the singularities lie in $\mathbb Z[ \sqrt{2 + \sqrt{2}} ]$.
The minimal polynomial of $\sqrt{2 + \sqrt{2}}$ over $\mathbb Z$ is $f( x ) = x^4 - 4 x^2 + 2$. 
For which prime numbers $p$ is this polynomial irreducible over $\mathbb Z_p$?

Suppose that $f$ factors into two quadratic terms over $\mathbb Z_p$.
Then there are $a, b \in \mathbb Z_p$ such that $\overline f( x ) = ( x^2 - ax + b )( x^2 + ax + b )$.\footnote{One should also consider the possibility that $f(x)$ factors as $(x^2+c)(x^2+d)$.  But in that case $x^2-4x+2$ would have a root in $\mathbb Z_p$, so $2$ would be a square modulo $p$, which Werner excludes by the end of the paragraph.}
Thus $b^2 \equiv 2 \pod p$.
If $p \equiv 3 \pod 8$ or $5 \pod 8$, then this is impossible; in this case 2 is not a square modulo $p$.
A prime $p \equiv 5 \pod 8$ with $\overline f( a ) \not\equiv 0 \pod p$ for all $a \in \mathbb Z_p$ is therefore suitable for reduction; 181 is such a prime.

\wernerpage{55}

We have the following results:
\[
\begin{tabular}{m{5.5cm}<{\centering} m{1.5cm}<{\centering} c m{2.6cm}<{\centering}}
    Branching surface & Number \text{of nodes} & Defect & \text{Number of} nullhomologous curves on $\hat V$ \\ \hline
    $\mathscr T_6( x_1 ) + \mathscr T_6( x_2 ) + \mathscr T_6( x_3 ) + 1 = 0$ & 54 & 6 & 0 \\
    $\mathscr T_6( x_1 ) + \mathscr T_6( x_2 ) + \mathscr T_6( x_3 ) - 1 = 0$ & 36 & 0 & 36 \\
    $\mathscr T_6( x_1 ) + \mathscr T_6( x_2 ) - \mathscr T_6( x_3 ) + 1 = 0$ & 51 & 5 & 0 \\
    $\mathscr T_6( x_1 ) - \mathscr T_6( x_2 ) - \mathscr T_6( x_3 ) + 1 = 0$ & 44 & 2 & 8 \\
    $\mathscr T_8( x_1 ) + \mathscr T_8( x_2 ) + \mathscr T_8( x_3 ) + 1 = 0$ & 144 & 9 & 0 \\
    $\mathscr T_8( x_1 ) + \mathscr T_8( x_2 ) + \mathscr T_8( x_3 ) - 1 = 0$ & 108 & 0 & 108 \\
    $\mathscr T_8( x_1 ) + \mathscr T_8( x_2 ) - \mathscr T_8( x_3 ) + 1 = 0$ & 136 & 7 & 0 \\
    $\mathscr T_8( x_1 ) - \mathscr T_8( x_2 ) - \mathscr T_8( x_3 ) + 1 = 0$ & 123 & 3 & 27
\end{tabular}
\]

Upper bounds on the defect were supplied by a computer, and lower bounds were obtained by giving $d$ linearly independent vectors in $\mathscr A$.
The only difficulty was in the case of the octic with 144 double points, where one must find 9 global divisors on $V$ whose associated 144-tuples in $\mathscr A$ are linearly independent! \bigskip \pagebreak 

Recall that $\mathscr T_8( x ) + \mathscr T_8( y )$ decomposes into 4 irreducible polynomials of degree 2.
These are (up to a constant factor)
\begin{align*}
    f_1( x, y ) &= x^2 + y^2 + \sqrt{ 2 + \sqrt{2} }\, x y + \frac{\sqrt{2} - 2}{4}
    \\
    f_2( x, y ) &= x^2 + y^2 - \sqrt{ 2 + \sqrt{2} }\, x y + \frac{\sqrt{2} - 2}{4}
    \\
    f_3( x, y ) &= x^2 + y^2 - \sqrt{ 2 + \sqrt{2} }\, x y - \frac{\sqrt{2} - 2}{4}
    \\
    f_4( x, y ) &= x^2 + y^2 + \sqrt{ 2 + \sqrt{2} }\, x y - \frac{\sqrt{2} - 2}{4}.
\end{align*}

\wernerpage{56}

Furthermore we need that $\mathscr T_8( z ) + 1$ is the square of a polynomial $g( z )$ of degree 4.
Up to a constant factor,
\[
    g( z ) = z^4 - z^2 + \frac{1}{8}.
\]
The equation of the double solid can be written as\footnote{Again $\omega$ comes from writing the double solid as $\omega^2 = $ octic.}
\[
    \prod_{i = 1}^4 f_i( x_1, x_2 ) = \big[ \omega - g( x_3 ) \big] \big[ \omega + g( x_3 ) \big].
\]
Choose as special 
divisors on $V$
\begin{align*}
    D_1 &= \{ \omega - g( x_3 ) = 0,\ f_1( x_1, x_2 ) = 0 \}
    \\
    D_2 &= \{ \omega - g( x_3 ) = 0,\ f_2( x_1, x_2 ) = 0 \}
    \\
    D_3 &= \{ \omega - g( x_3 ) = 0,\ f_3( x_1, x_2 ) = 0 \}
    \\
    D_4 &= \{ \omega - g( x_2 ) = 0,\ f_1( x_1, x_3 ) = 0 \}
    \\
    D_5 &= \{ \omega - g( x_2 ) = 0,\ f_2( x_1, x_3 ) = 0 \}
    \\
    D_6 &= \{ \omega - g( x_2 ) = 0,\ f_3( x_1, x_3 ) = 0 \}
    \\
    D_7 &= \{ \omega - g( x_1 ) = 0,\ f_1( x_2, x_3 ) = 0 \}
    \\
    D_8 &= \{ \omega - g( x_1 ) = 0,\ f_2( x_2, x_3 ) = 0 \}
    \\
    D_9 &= \{ \omega - g( x_1 ) = 0,\ f_3( x_2, x_3 ) = 0 \}.
\end{align*}
Let
\[
    \alpha_1 = \frac{\sqrt{ 2 + \sqrt{2} }}{2} = -\alpha_7
\]
and
\[
    \alpha_3 = \frac{\sqrt{ 2 - \sqrt{2} }}{2} = -\alpha_5.
\]
We now examine how these divisors behave in a neighborhood of specific 
singularities. \bigskip

\wernerpage{57}

Let $P_1 := ( \alpha_3, 0, \alpha_3 )$.
The divisors $D_1$, $D_2$, $D_7$, and $D_8$ go through this point: $D_1$ and $D_2$ represent the two local divisors at $P_1$, 
as do $D_7$ and $D_8$.
The coordinate transformation $x_3 \leftrightarrow x_1$ leaves $P_1$ invariant and exchanges $D_1$ and $D_7$.
Locally, this coordinate transformation exchanges the two small resolutions, so $D_1$ and $D_7$ are different locally at $P_1$.
We write it thus:
\[
    D_1( 1 ),\ D_2( -1 ),\ D_7( -1 ),\ D_8( 1 ).
\]
In parentheses are the intersection numbers with the exceptional curve of a small resolution of $P_1$. \bigskip

Let $P_2 := ( \alpha_1, 0, \alpha_1 )$.
Here we have $D_3( 1 )$ and $D_9( -1 )$. \bigskip

Choose $P_3 := ( \alpha_3, 0, \alpha_5 )$.
The coordinate transformation $x_3 \leftrightarrow -x_1$ leaves $P_3$ fixed and maps $D_1$ to $\{ \omega + g( x_1 ) = 0, f_1( x_2, x_3 ) = 0 \}$.
Thus here we have
\[
    D_1( 1 ), D_2( -1 ), D_7( 1 ), D_8( -1 ).
\]
Analogously we get
\[
    \begin{array}{l l l l l l}
    P_4 := ( 0, \alpha_3, \alpha_3 ) & : & D_1( 1 ), & D_2( -1 ), & D_4( -1 ), & D_5( 1 )
    \\
    P_5 := ( 0, \alpha_1, \alpha_1 ) & : & D_3( 1 ), & D_6( -1 ) & &
    \\
    P_6 := ( \alpha_3, -\frac{1}{\sqrt{2}}, \alpha_3 ) & : & D_1( 1 ), & D_3( -1 ), & D_7( -1 ), & D_9( 1 )
    \\
    P_7 := ( -\frac{1}{\sqrt{2}}, \alpha_3, \alpha_3 ) & : & D_1( 1 ), & D_3( -1 ), & D_4( -1 ), & D_6( 1 )
    \\
    P_8 := ( \alpha_3, \frac{1}{\sqrt{2}}, \alpha_3 ) & : & D_2( 1 ), & D_7( -1 ), & D_9( 1 ) &
    \end{array}
\]

\wernerpage{58}

Choose the last point to be $P_9 := ( \alpha_1, 1 / \sqrt{2}, \alpha_3 )$.
The divisors $D_2$, $D_3$, and $D_8$ contain $P_9$, and $D_2$ and $D_3$ represent the two local divisors at $P_9$.
Which of these two divisors does $D_8$ agree with locally?
    
Consider the local meromorphic functions
\[
    F_1 := \frac{ g( x_3 ) }{ f_2( x_1, x_2 ) }, \qquad F_2 := \frac{ g( x_3 ) }{ f_3( x_1, x_2 ) }, \quad\text{and}\quad F_3 := \frac{ g( x_1 ) }{ f_2( x_2, x_3 ) }
\]
in a small affine neighborhood of $P_1$.
Choose sequences $0 < \epsilon_n \to 0$ and $0 < \delta_n \to 0$, such that $( \alpha_1 + \epsilon_n, 1 / \sqrt{2} + \delta_n, \alpha_3 + \epsilon_n )$ lies in $B$ [the neighborhood].
Then
\begin{align*}
    F_1( \alpha_1 + \epsilon_n, 1 / \sqrt{2} + \delta_n, \alpha_3 + \epsilon_n ) &> 0,
    \\
    F_2( \alpha_1 + \epsilon_n, 1 / \sqrt{2} + \delta_n, \alpha_3 + \epsilon_n ) &< 0,
    \\
    \text{and} \qquad    
    F_3( \alpha_1 + \epsilon_n, 1 / \sqrt{2} + \delta_n, \alpha_3 + \epsilon_n ) &> 0.
\end{align*}
Therefore $D_8$ locally agrees with $D_2$, which is to say: $D_2( 1 )$, $D_3( -1 )$, $D_8( 1 )$.

\wernerpage{59}

Assemble 
the calculated intersection numbers on a particular 
small resolution in the following matrix $( \hat D_i . L_{P_j} )_{\substack{1 \leq i \leq 9 \\ 1 \leq j \leq 9}}$:
\[
    \begin{pmatrix*}[r]
    1 & 0 & 1 & 1 & 0 & 1 & 1 & 0 & 0
    \\
    -1 & 0 & -1 & -1 & 0 & 0 & 0 & 1 & 1
    \\
    0 & 1 & 0 & 0 & 1 & -1 & -1 & 0 & -1
    \\
    0 & 0 & 0 & -1 & 0 & 0 & -1 & 0 & 0
    \\
    0 & 0 & 0 & 1 & 0 & 0 & 0 & 0 & 0
    \\
    0 & 0 & 0 & 0 & -1 & 0 & 1 & 0 & 0
    \\
    -1 & 0 & 1 & 0 & 0 & -1 & 0 & -1 & 0
    \\
    1 & 0 & -1 & 0 & 0 & 0 & 0 & 0 & 1
    \\
    0 & -1 & 0 & 0 & 0 & 1 & 0 & 1 & 0
    \end{pmatrix*}
\]
This matrix has rank 9, so the 144-tuples in $\mathscr A$ belonging to the divisors $D_1, \ldots, D_9$ are linearly independent.
The defect is equal to 9, and the projective algebraic small resolutions are K\"ahler manifolds with trivial canonical bundle and Euler characteristic $-8$.
The interesting Hodge numbers are $h^{1, 1}( \hat V ) = 10$ and $h^{2, 1}( \hat V ) = 14$.
\chapter{Estimates for \texorpdfstring{$\mu_n(d)$}{mu\_n(d)}} \label{chapter7}
\wernerpage{60}

In this chapter we will briefly leave small resolutions, and occupy ourselves with the question of how many ordinary double points a hypersurface of degree $d$ in $\mathbb P^n$ can have.
Since the chances of finding a projective algebraic small resolution for a given nodal threefold become greater as the number of singularities increases, such a question is definitely interesting in this context. \bigskip

About notation:
following Varchenko, in this chapter we denote the dimension of a projective space by $n$ and the degree of a hypersurface by $d$.
Since defects will not reappear until the next chapter, this should not cause confusion.
The maximal number of nodes of a hypersurface of degree $d$ in $\mathbb P^n( \mathbb C )$ is denoted by $\mu_n( d )$. \bigskip

A general estimate for $\mu_n( d )$ comes from Bruce \cite{bruce81}.
He shows that
\begin{align*}
    \mu_n( d ) &\le \frac{ (d - 1)^n (d + 1) + d - 1 }{ 2 d }
    && \text{if $n$ is even,}
    \\
    \mu_n( d ) &\le \frac{ (d - 1)^n }{ 2 }
    && \text{if $d$ and $n$ are odd,}
    \\
    \mu_n( d ) &\le \frac{ (d - 1)^n (d + 1) + 1 }{ 2 d }
    && \text{if $d$ is even and $n$ is odd.}
\end{align*}
Varchenko's estimate in \cite{varchenko831} is sharper:
\[
    \mu_n( d ) \leq A_n( d ),
\]
where $A_n( d )$ is the number of lattice points
\wernerpage{61}
$(k_1, \ldots, k_n)$ inside the cube $(0, d)^n$ 
such that
\[
    (n - 2)\, d / 2 + 1 < \sum_{i = 1}^n k_i \leq n d / 2 .
\]

Kalker shows in his dissertation \cite{kalker86} that
\[
    \mu_n( 3 ) = A_n( 3 ) \qquad \text{for all }n \in \mathbb N,
\]
so Varchenko's estimate is sharp for cubics.
For surfaces in $\mathbb P^3$, Miyaoka proved in \cite{miyaoka84} that
\[
    \mu_3( d ) \leq \frac{4}{9} \, d ( d - 1 )^2.
\]
The following table contains all known results for $n = 3$ and $2 \leq d \leq 12$:\footnote{The upper bounds 256 and 488 should be 246 and 480, following Varchenko rather than Miyaoka.  The lower bound 66 is incorrect: Jaffe and Ruberman \cite{jaffe97} later proved that $\mu_3(6) = 65$.  The lower bounds for $d \ge 7$ have since been improved; at this writing the Wikipedia article for ``Nodal surface'' has a good list. \label{65_vs_66_footnote}}
\begin{center}
\begin{tabular}{c c c c c c c c c c c c}
    $d$ & 2 & 3 & 4 & 5 & 6 & 7 & 8 & 9 & 10 & 11 & 12
    \\ \hline
    $\mu_3( d ) \leq$ & 1 & 4 & 16 & 31 & $66^*$ & 104 & 174 & $256^*$ & 360 & $488^*$ & 645
    \\ \hline
    $\mu_3( d ) \geq$ & 1 & 4 & 16 & 31 & $66^*$ & 90 & 160 & 192 & 300 & 375 & 540
\end{tabular}
\end{center}
Beauville \cite{beauville80} proved that $\mu_3( 5 ) = 31$, and Stagnaro \cite{stagnaro85} proved that $\mu_3( 6 ) = 66$.
The surface of degree 7 with 90 nodes was also discovered by Stagnaro.
The octic with 160 double points comes from Kreiss \cite{kreiss56}. \bigskip

For hypersurfaces in $\mathbb P^4$ we have the following individual results:
\begin{align*}
    \text{Bruce:} && \mu_4( 2 ) &\leq 1 & \mu_4( 3 ) &\leq 11 & \mu_4( 4 ) &\leq 51 & \mu_4( 5 ) &\leq 154 \\
    \text{Varchenko:} && \mu_4( 2 ) &\leq 1 & \mu_4( 3 ) &\leq 10 & \mu_4( 4 ) &\leq 45 & \mu_4( 5 ) &\leq 135  
\end{align*}
\wernerpage{62}
Examples that will be discussed in more detail later show us that
\begin{align*}
    \mu_4( 2 ) &= 1  &  \mu_4( 3 ) &= 10  &  \mu_4( 4 ) &= 45.
\end{align*}
Quintics with 125 and 126 nodes were constructed by Schoen \cite{schoen86} and Hirzebruch \cite{hirzebruch85}, respectively.
Thus $126 \leq \mu_4( 5 ) \leq 135$. \bigskip

In what follows we are interested in the asymptotic behavior of $\mu_n( d )$.

\begin{dfn}
Let
\begin{align*}
    a_n &:= \lim_{d \to \infty} \frac{ A_n( d ) }{ d^n }, \qquad b_n^1 := \varlimsup_{d \to \infty} \frac{ \mu_n( d ) }{ d^n }
    \\[1ex]
    b_n^2 &:= \varliminf_{d \to \infty} \frac{ \mu_n( d ) }{ d^n }, \qquad c_n := \binom{ n }{ \lfloor n / 2 \rfloor } 2^{-n}.
\end{align*}
\end{dfn}

The asymptotic number of singularities in a Chmutov hypersurface as $d \to \infty$ is given by $c_n$. 
Then we have
\[
    c_n \leq b_n^2 \leq b_n^1 \leq a_n.
\]
From Bruce's results it follows immediately that $b_n^1 \leq \frac{1}{2}$ for all $n \in \mathbb N$.
In the cases $n = 3, 4$, we have the following more precise estimates:
\begin{align*}
    3/8 &\leq b_3^2 \leq b_3^1 \leq 4/9
    \\
    3/8 &\leq b_4^2 \leq b_4^1 \leq 11/24.
\end{align*}
The lower bounds come from Chmutov hypersurfaces, the upper bound for $b_3^1$ from Miyaoka, and the upper bound for $b_4^1$ from Varchenko. 

\begin{lem}
\wernerpage{63}
\[
    a_n = \text{Volume} \left\{ ( x_1, \ldots, x_n ) \in [0, 1]^n : (n - 2)/2 \leq \sum_{i = 1}^{n} x_i \leq n/2 \right\}
\]
\end{lem}

The proof is immediate from the definition of the so-called Arnold numbers $A_n( d )$. \bigskip

By direct integration one calculates $a_3 = 23/48$ and $a_4 = 11/24$.
The sequence $a_n$ has yet another description: the condition
\[
    (n - 2)/2 \leq \sum_{i = 1}^{n} x_i \leq n/2
\]
is equivalent to
\[
    -1 \leq \sum_{i = 1}^{n} ( x_i - \tfrac12 ) \leq 0.
\]
Interpret $( x_i - \tfrac12 )$ as a random variable that is uniformly distributed in the interval $[-\tfrac12,\,\tfrac12]$.
The Fourier transform of this density is
\[
    \int_{-\frac{1}{2}}^{\frac{1}{2}} e^{i \lambda x} \, dx = \frac{\sin( \lambda / 2 )}{\lambda / 2}.
\]
Then $\sum_{i = 1}^{n} ( x_i - \frac{1}{2})$ has a distribution function with density
\[
    \frac{1}{2 \pi} \int_{\mathbb R} e^{-i \lambda x} \bigg( \frac{\sin( \lambda / 2 )}{\lambda / 2} \bigg)^n \, d\lambda.
\]
Therefore
\begin{align*}
    a_n 
    &= \frac{1}{2 \pi} \int_{-1}^{0} \int_{\mathbb R} e^{-i \lambda x} \bigg( \frac{\sin( \lambda / 2 )}{\lambda / 2} \bigg)^n \, d\lambda \, dx
    \\
    &= \frac{1}{2 \pi} \int_0^1 \int_{\mathbb R} \cos( \lambda x ) \bigg( \frac{\sin( \lambda / 2 )}{\lambda / 2} \bigg)^n \, d\lambda \, dx.
\end{align*}
\wernerpage{64}
For $n \geq 2$, we follow with Fubini's Theorem:
\begin{align*}
    a_n 
    &= \frac{1}{2 \pi} \int_{\mathbb R} \bigg( \int_{0}^{1} \cos( \lambda x ) \, dx \bigg) \bigg( \frac{\sin( \lambda / 2 )}{\lambda / 2} \bigg)^n \, d\lambda
    \\
    &= \frac{1}{2 \pi} \int_{\mathbb R} \frac{\sin \lambda }{\lambda} \bigg( \frac{\sin( \lambda / 2 )}{\lambda / 2} \bigg)^n \, d\lambda
    \\
    &= \frac{1}{\pi} \int_{0}^{\infty}  \frac{\sin 2 x }{ x }  \bigg( \frac{\sin x }{x} \bigg)^n \, dx.
\end{align*}
This is equivalent to
\begin{align*}
    a_n 
    &= \frac{2}{\pi} \int_{0}^{\infty} \bigg( \frac{\sin x }{x} \bigg)^{n + 1} \cos x  \, dx.
    \\
    &= \frac{2}{\pi} \bigg( \frac{n + 1}{n + 2} \bigg) \int_{0}^{\infty} \bigg( \frac{\sin x }{x} \bigg)^{n + 2} \, dx.
\end{align*}
The variances of the independent random variables $(x_i - \frac{1}{2})$ are equal to $\frac{1}{12}$.
Thus the sequence
\[
    S_n^* := \frac{\sum_{i = 1}^n ( x_i - \frac{1}{2} )}{\sqrt{ n / 12}}
\]
converges to the normal distribution by the central limit theorem.
Let $F_n$ be the [density of the] 
distribution function of $S_n^*$.
Then
\[
    a_n = \int_{- \sqrt{12 / n}}^0 F_n( x ) \, dx.
\]
Since $F_n$ is uniformly bounded, $a_n$ converges to $0$.
Therefore $b_n^1$ and $b_n^2$ also converge to 0.

\wernerpage{65}

Moreover we know that the sequence of functions $F_n$ converges uniformly on every compact interval to the density $\phi$ of the normal distribution.
By continuity of $\phi$, we have
\begin{gather*}
    \lim_{n \to \infty} \sqrt{n} \, a_n 
    = \lim_{n \to \infty} \sqrt{n} \int_{- \sqrt{12 / n}}^0 F_n( x ) \, dx
    \\
    = \lim_{n \to \infty} \sqrt{n}\, \sqrt{12 / n}\ \phi( 0 )
    = \sqrt{12} / \sqrt{2\pi}
    = \sqrt{6 / \pi}.
\end{gather*}
Using Stirling's formula we easily calculate that
\[
    \lim_{n \to \infty} \sqrt{n} \, c_n = \sqrt{2 / \pi}.
\]
Both of these formulas were already cited by Varchenko in \cite{varchenko832}, albeit without proof.

In summary, we have the inequality
\[
    \sqrt{2 / \pi} \leq \varliminf_{n \to \infty} b_n^1 \leq \varlimsup_{n \to \infty} b_n^2 \leq \sqrt{6 / \pi}.
\]

\chapter{Cubics in \texorpdfstring{$\mathbb P^4(\mathbb C)$}{P\^{}4(C)}} \label{chapter8}
\wernerpage{66}

This chapter essentially builds on results from Ton Kalker's dissertation \cite{kalker86}, and was developed in close collaboration with Hans Fin\-kelnberg (Leiden); see also \cite{finkelnbergIV}.
We will examine nodal cubic hypersurfaces in $\mathbb P^4( \mathbb C )$ more closely. \bigskip

Let $K$ be a singular cubic in $\mathbb P^4$ with only ordinary double points and a singularity at the point $P_1 := ( 0 : 0 : 0 : 0 : 1 )$.
Then
\[
    K = \{ x_4 Q + R = 0 \},
\]
where $Q, R \in \mathbb C[ x_0, \dotsc, x_3 ]$ are homogeneous polynomials of degrees 2 and 3, respectively.
We let $K( P_1 )$ denote $K$ blown up at $P_1$; the quadric $\{ Q = 0 \}$ corresponds to the exceptional surface $\mathbb P^1 \times \mathbb P^1$.
The restriction of the map
\[
    \begin{array}{ r r c l}
    \Pi : & \mathbb P^4 - \{ P_1 \} & \longrightarrow & \mathbb P^3  \\
          & ( x_0 : \dotsc : x_4 ) & \longmapsto & ( x_0 : \dotsc : x_3 )  
    \end{array}
\]
to $K - \{ P_1 \}$ can be extended to a map $\tilde\Pi : K( P_1 ) \to \mathbb P^3$.
The curve
\[
    S := \{ Q = 0 \} \cap \{ R = 0 \} \subset \mathbb P^3
\]
is called the \textit{associated curve}.
is mapped bijectively onto the singularities of $S$ via $\Pi$.
The preimage $\Pi^{-1}( S )$ is contained in $K$.
Away from $S$, the map $\tilde \Pi$ is bijective.
Proofs of these facts can be found in \cite{kalker86}.

\wernerpage{67}

The blow-up $K( P_1 )$ is isomorphic to $\mathbb P^3$ blown up along $S$. 
Thus we have
\[
    H_4( K( P_1 ) , \mathbb Z ) \cong \mathbb Z \oplus \bigoplus_{i} \mathbb Z \, 
\tilde \Pi^{-1} S_i,
\]
where the $S_i$ run through the irreducible components of $S$. 
Together with the homology class of the generic hyperplane section $H$, the divisors $\Pi^{-1} S_i $ generate $H_4( K, \mathbb Z )$.

The associated curve $S$, thought of as a curve in $Q \cong \mathbb P^1 \times \mathbb P^1$, is of type $(3, 3)$.
Using the associated curve we can give a new characterization of the cubics that admit projective algebraic small resolutions.

\begin{thm}
Let $K$ be a nodal cubic in $\mathbb P^4( \mathbb C )$.
A projective algebraic small resolution exists if and only if all irreducible components of $S$ are smooth and at least one is of type $(a, b)$ with $a \neq b$.

The number of irreducible components of $S$ is equal to $d + 1$. 
\end{thm}

\begin{proof}
Let $P \in \mathscr S - \{ P_1 \}$.
Then $\Pi( P )$ is a singularity of $S$.
There are thus two possibilities:
\begin{enumerate}[label=\arabic*.]
    \item Two locally smooth irreducible components of $S$ intersect at $\Pi( P )$.
    
    \item An irreducible component of $S$ intersects itself at $\Pi( P )$.
\end{enumerate}
In the first case, let $S_1$ be a component that is smooth at $\Pi( P )$.
Then $\Pi^{-1} S_1 $ is a divisor on $K$ which is smooth at $P$.
\wernerpage{68}
On a small resolution $\hat V$, the exceptional line $L_P$ intersects the proper transform of this divisor with multiplicity $\pm 1$.
Thus the curve $L_P$ is not nullhomologous on $\hat V$.

In the second case, let $S_2$ denote the irreducible component that is singular at $\Pi( P )$.
Then we have
\[
    ( \Pi^{-1} S_2  ) . L_p = \pm( 1 - 1 ) = 0.
\]
Of course this also holds for the other components that do not meet $\Pi( P )$.
Thus $L_P$ is nullhomologous on $\hat V$. \bigskip

It remains to examine the point $P_1$.

Let $S_i$ be an arbitrary component of $S$ of type $(a, b)$.
Because
\[
    ( \tilde \Pi^{-1} S_i ) . Q \cong S_i,
\]
we have
\begin{align*}
    ( \tilde \Pi^{-1} S_i ) . A_1 &= -a,
    \\
    ( \tilde \Pi^{-1} S_i ) . B_1 &= -b \quad \text{and}
    \\
    ( \tilde \Pi^{-1} S_i ) . L_1 &= \pm( a - b ).
\end{align*}
Hence $L_1$ is not nullhomologous on $\hat V$ precisely when a component of $S$ has type $(a, b)$ with $a \neq b$ (cf.\ \cite{finkelnberg87}, \cite{finkelnbergIV}). \bigskip

As $S$ is connected, one can read off from the intersection numbers that
\[
    \left( \sum_{i} \alpha_i\, \Pi^{-1} S_i \right) .\, L_p = 0
\]
for all $P \in \mathscr S$ if and only if $\alpha_i = \alpha_j$ for all $i, j$.
Hence only one relation exists between the divisors $\Pi^{-1} S_i$ and $H$, so the number of irreducible components of $S$ is exactly $\beta_4( V ) = d + 1$.
\end{proof}

\wernerpage{69}

Conversely we have the following

\begin{lem}
Every curve of type (3, 3) in $Q \subset \mathbb P^3$, $Q \cong \mathbb P^1 \times \mathbb P^1$ is the complete intersection of a cubic with the quadric $Q$, and therefore is the associated curve of a singular cubic hypersurface in $\mathbb P^4$.
\end{lem}

\begin{proof}
Without loss of generality let $Q = \{ x_0 x_3 - x_1 x_2 = 0 \}$.
Choose coordinates $( u_0 : u_1\:;\:v_0 : v_1 )$ in $\mathbb P^1 \times \mathbb P^1$ so that $x_0 = u_0 v_0$, $x_1 = u_0 v_1$, $x_2 = u_1 v_0$, and $x_3 = u_1 v_1$.
A (possibly reducible) curve $S$ of type (3, 3) in $Q$ is described by a bihomogeneous equation
\[
    F( u_0, u_1 ; v_0, v_1 ) = 0
\]
of bidegree (3, 3).
There exists a homogeneous polynomial $R( x_0, x_1, x_2, x_3 )$ of degree 3 such that
\[
    R( u_0 v_0, u_0 v_1, u_1 v_0, u_1 v_1) = F( u_0, u_1 ; v_0, v_1 ).
\]
Then $\{ R = 0 \}$ describes a cubic in $\mathbb P^3$ whose intersection with $Q$ is exactly $S$.
Hence
\[
    \big\{ x_4( x_0 x_3 - x_1 x_2 ) + R( x_0, x_1, x_2, x_3 ) = 0 \big\}
\]
is a cubic in $\mathbb P^4$ with $S$ as its associated curve.
\end{proof}

Using 
the various nodal curves of type (3, 3) in $\mathbb P^1 \times \mathbb P^1$, we can now classify all nodal cubics in $\mathbb P^4( \mathbb C )$ according to 
the homology groups of their small resolutions.
Consider how a curve of type (3, 3) in $\mathbb P^1 \times \mathbb P^1$ can be decomposed into irreducible components, and note that an irreducible component of type $(a, b)$ can have up to $(a - 1)(b - 1)$ ordinary double points.
\wernerpage{70}

One obtains the following complete table:
\[
    \begin{tabular}{c m{2.25cm}<{\centering} m{3cm}<{\centering} m{4.5cm}<{\centering}}
    Defect & Number of Nodes & Nodes resolved nullhomologously & Number of {projective algebraic} {small resolutions}  
    \\ \hline
    0 & 1 & 1 & 0
    \\
      & 2 & 2 & 0
    \\
      & 3 & 3 & 0
    \\
      & 4 & 4 & 0
    \\
      & 5 & 5 & 0
    \\ \hline
    1 & 4 & 0 & 2
    \\
      & 5 & 1 & 0
    \\
      & 6 & 2 & 0
    \\
      & 6 & 0 & 2
    \\ \hline
    2 & 6 & 0 & 6
    \\
      & 7 & 1 & 0
    \\
      & 7 & 0 & 6
    \\ \hline
    3 & 8 & 0 & 24
    \\
    4 & 9 & 0 & 102
    \\
    5 & 10 & 0 & 332
    \\ \hline
    \end{tabular}
\]

The (uniquely determined) cubic with 10 double points is examined in detail by Finkelnberg in \cite{finkelnberg87}; he also calculates the number of projective algebraic small resolutions to be 332.
Details of the other cubics will be included in \cite{finkelnbergIV}, where the case $s = 9$ will be discussed as an example. \bigskip

\wernerpage{71}

Let us assume the condition ``$d = 4$.''
Then $S$ decomposes into two components of type (1, 0), two of type (0, 1), and an irreducible (1, 1)-component, so $s = 9$.
The components of $S$ and the singularities are labeled as follows:
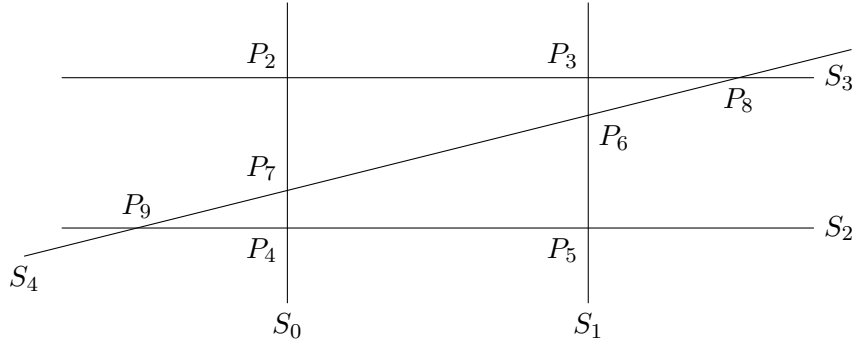
\[
\begin{tikzpicture}
\draw (-2, -2) -- (-2, 2);
\draw (2, -2) -- (2, 2);
\draw (-5, -1) -- (5, -1);
\draw (-5, 1) -- (5, 1);
\draw (-5.5, -1.375) -- (5.5, 1.375);
\draw (-2, -2) node[anchor=north] {$S_0$};
\draw (2, -2) node[anchor=north] {$S_1$};
\draw (5, -1) node[anchor=west] {$S_2$};
\draw (5, 1) node[anchor=west] {$S_3$};
\draw (-5.5, -1.375) node[anchor=north] {$S_4$};
\draw (-2, 1) node[anchor=south east] {$P_2$};
\draw (2, 1) node[anchor=south east] {$P_3$};
\draw (-2, -1) node[anchor=north east] {$P_4$};
\draw (2, -1) node[anchor=north east] {$P_5$};
\draw (2, 0.5) node[anchor=north west] {$P_6$};
\draw (-2, -0.5) node[anchor=south east] {$P_7$};
\draw (4, 1) node[anchor=north] {$P_8$};
\draw (-4, -1) node[anchor=south] {$P_9$};
\end{tikzpicture}
\]
The intersection numbers of the proper transforms of the divisors $\Pi^{-1} S_i$ with the exceptional curves on a small resolution are assembled 
in the following matrix: \[
    \begin{array}{c | r r r r r r r r r }
    & 1 & 2 & 3 & 4 & 5 & 6 & 7 & 8 & 9
    \\ \hline
    S_0 & 1 & 1 & 0 & 1 & 0 & 0 & 1 & 0 & 0
    \\
    S_1 & 1 & 0 & 1 & 0 & 1 & 1 & 0 & 0 & 0
    \\
    S_2 & -1 & 0 & 0 & -1 & -1 & 0 & 0 & 0 & -1
    \\
    S_3 & -1 & -1 & -1 & 0 & 0 & 0 & 0 & -1 & 0
    \\
    S_4 & 0 & 0 & 0 & 0 & 0 & -1 & -1 & 1 & 1
    \end{array}
\]

\wernerpage{72}

The row vectors of this matrix generate $\mathscr A$.
Therefore $\mathscr B$ is generated by the 9-tuples
\begin{align*}
    &e_1 - e_2 - e_5,
    \\
    &e_1 - e_3 - e_4,
    \\
    &e_2 - e_7 - e_8,
    \\
    &e_4 - e_7 - e_9,
    \\
    &e_3 - e_6 - e_8.
\end{align*}
These completely describe the linear relations between the exceptional curves.
The division into projective and non-projective small resolutions is obtained by specifying a divisor $D$ with $D . L_P > 0$ for all $P \in \mathscr S$ or an effective nullhomologous curve on the various small resolutions $\hat V$, respectively. 
In chapter \ref{chapter12} we will give a simple proof that exactly 102 small resolutions of $V$ are projective algebraic.
\chapter{Quartics in \texorpdfstring{$\mathbb P^4(\mathbb C)$}{P\^{}4(C)}} \label{chapter9}
\wernerpage{73}

Some of the methods developed for studying cubics can be carried over to special 
quartic hypersurfaces in $\mathbb P^4$. \bigskip

A nodal quartic in $\mathbb P^4$ can, by a change of coordinates, be put in the form
\[
    V = \{ x_4^2 Q + x_4 K + R = 0 \},
\]
where $Q, K, R \in \mathbb C[ x_0, \ldots x_3 ]$ are homogeneous polynomials of degrees 2, 3, and 4, respectively.
The point $P_1 = ( 0 : 0 : 0 : 0 : 1 )$ is then a singularity, and the surface $\{ Q = 0 \} \subset \mathbb P^3$ is smooth.
We only want to consider quartics for which the middle term vanishes, so
\[
    V = \{ x_4^2 Q + R = 0 \}.
\]
Let $P \in \mathscr S - \{ P_1 \}$.
There are two possible cases:
\begin{enumerate}[label=\arabic*.]
    \item $x_4( P ) = 0$.
    Then $Q( P ) \neq 0$, $R( P ) = 0$ and $\partial_i R( P ) = 0$ for all $i = 0, \ldots, 3$.  (Here $\partial_i := \frac{\partial}{\partial x_i}$.)
    
    \item $x_4( P ) \neq 0$.
    Since $\partial_4 V = 2 x_4 Q$, we have $Q( P ) = 0$, $R( P ) = 0$, and $x_4^2 \partial_i Q + \partial_i R( P ) = 0$ for all $i = 0, \ldots, 3$. 
\end{enumerate}

The curve $S := \{ Q = 0 \} \cap \{ R = 0 \}$ is again called the associated curve, and $\Pi$ again denotes the map
\[
    \begin{array}{c c c}
    \mathbb P^4 - \{ P_1 \} & \longrightarrow & \mathbb P^3 \\
    ( x_0 : \ldots : x_4 ) & \longmapsto & ( x_0 : \ldots : x_3 ) 
    \end{array}
\]
\wernerpage{74}
The point $( \alpha_0 : \alpha_1 : \alpha_2 : \alpha_3 )$ is a singularity of $S$ if and only if there exists a pair $( \lambda, \mu ) \neq (0, 0)$ such that
\[
    \lambda \partial_i Q( \alpha ) + \mu \partial_i R( \alpha ) = 0
\]
for $i = 0, \ldots, 3$.

In our case, we must also have $\lambda \ne 0$ and $\mu \neq 0$!

If $\mu = 0$, then $\partial_i Q( \alpha ) = 0$ for $i = 0, \ldots, 3$, so $Q$ is singular at $\alpha$.
Contradiction!

If $\lambda = 0$, then $\partial_i R( \alpha ) = 0$ for $i = 0, \ldots, 3$, so $R$ is singular at $\alpha$.
The point $( \alpha_0 : \alpha_1 : \alpha_2 : \alpha_3 : 0 )$ lies in $V$, and all partial derivatives vanish there.
The point $( \alpha_0 : \alpha_1 : \alpha_2 : \alpha_3 : 0 )$ would then be a singularity of $V$ with $x_4 = 0$ and $Q = 0$.
Contradiction! \bigskip

So for every $\alpha \in S_{\text{sing}}$ there is a $\lambda \ne 0$ with
\[
    \lambda \partial_i Q( \alpha ) + \partial_i R( \alpha ) = 0
\]
for all $i = 0, \ldots, 3$.
The equation
\[
    x_4^2 = \lambda
\]
has two solutions.
The two points $(\alpha_0 : \alpha_1 : \alpha_2 : \alpha_3 : \alpha_4 )$ with $\alpha_4^2 = \lambda$ belong to $V$ because $\Pi^{-1} S \subset V$.
Hence two singularities of $V$ lie over each $\alpha \in S_{\text{sing}}$. \bigskip

Unfortunately, on the other hand, $\Pi$ does not map the nodes with $x_4( P ) = 0$ into $S$, but only those with with $x_4( P ) \neq 0$.
For the latter points we have $\Pi( P ) \in S_{\text{sing}}$.

\wernerpage{75}

Since the singularities of $V$ with $x_4( P ) = 0$ are in one-to-one correspondence with the singularities of $R$, we have
\[
    s = 1 + 2 \# S_{\text{sing}} + \# R_{\text{sing}}.
\]
What information does the associated curve contain about the singularities of $V$ and their resolving curves?

Nothing at all can be said about the singularities with $Q \neq 0$; here one must try to project from another point.
Nothing can be said either about the $P \in \mathscr S - \{ P_1 \}$ such that $\Pi( P )$ is a singular point of an irreducible component of $S$, as the cones over the irreducible components of $S$ need not generate all of $H_4( V, \mathbb Z )$.

Now let $P \in \mathscr S - \{ P_1 \}$, and suppose that $\Pi( P )$ is at the intersection of two locally smooth components of $S$.
Then a locally smooth divisor on $V$ goes through $P$, so $L_P$ is not nullhomologous on $\hat V$.

Generally speaking, one only has complete information about $P_1$: the resolving curve is then not nullhomologous on $\hat V$ precisely when at least one component of $S$ is of type $(a, b)$ with $a \neq b$.  The associated curve $S$ is of type (4, 4) in $\{ Q = 0 \}$. \bigskip

\wernerpage{76}

As an example, we will investigate a quartic with 45 double points, the maximum number of nodes possible.
Let
\[
    V = \left\{ \sum_{i = 0}^{5} x_i = 0,\ \sigma_4( x_0, \ldots, x_5 ) = 0 \right\} \subset \mathbb P^5, 
\]
where $\sigma_i$ denotes the $i$-th elementary symmetric function.
Todd treated this variety in detail in \cite{todd47}.
The set of singularities $\mathscr S$ is the orbit of the two points $(1 : -1 : 0 : 0 : 0 : 0 )$ and $(1 : 1 : \epsilon : \epsilon : \epsilon^2 : \epsilon^2 )$, where $\epsilon = e^{2\pi i / 3}$, under the action $S_6$ on $\mathbb P^5$ that permutes the coordinates.
If we choose $P_1$ to be the point $(0 : 0 : 0 : 0 : 1 : -1 )$, a coordinate change
\begin{align*}
    y_4 &:= (x_4 - x_5)/2
    \\
    y_5 &:= (x_4 + x_5)/2
    \\
    y_i &:= x_i, \qquad i \neq 4, 5
\end{align*}
brings the quartic into the form
\[
    2 y_5 + \sigma_1 = 0 \qquad \text{and} \qquad ( y_5^2 - y_4^2 ) \sigma_2 + 2 y_5 \sigma_3 + \sigma_4 = 0,
\]
where $\sigma_i = \sigma_i( y_0, y_1, y_2, y_3 )$.
In  $\mathbb P^4$ with homogeneous coordinates \linebreak 
$( y_0 : y_1 : y_2 : y_3 : y_4 )$, we get
\[
    V = \left\{ y_4^2 \sigma_2 + \left( \sigma_1 \sigma_3 - \sigma_4 + \frac{\sigma_1^2 \sigma_2}{4} = 0 \right) \right\}.
\]
Then $P_1 = ( 0 : 0 : 0 : 0 : 1 )$ and $Q = \sigma_2$.

\wernerpage{77}

The singularities can be divided into classes as follows:
\[
    \begin{array}{c | l | l | c | c}
    \text{Type} & x\text{-coordinates in }\mathbb P^5 & y\text{-coordinates in }\mathbb P^4 & Q & \text{Number}
    \\ \hline
    0 & ( 0 : 0 : 0 : 0 : 1 : -1 ) & ( 0 : 0 : 0 : 0 : 1 ) & 0 & 1
    \\
    1 & ( .\ldots \ldots \ldots : 1 : 0 ) & ( .\ldots \ldots \ldots : \frac{1}{2} ) & 0 & 4
    \\
    2 & ( .\ldots \ldots \ldots : 0 : 1 ) & ( .\ldots \ldots \ldots : -\frac{1}{2} ) & 0 & 4
    \\
    3 & ( .\ldots \ldots \ldots : 0 : 0 ) & ( .\ldots \ldots \ldots : 0 ) & -1 & 6
    \\
    4 & ( .\ldots \ldots \ldots : 1 : 1 ) & ( .\ldots \ldots \ldots : 0 ) & 3 & 6
    \\
    5 & ( .\ldots \ldots \ldots : 1 : \epsilon ) & ( .\ldots \ldots \ldots : a ) & 0 & 12
    \\
    6 & ( .\ldots \ldots \ldots : 1 : \epsilon^2 ) & ( .\ldots \ldots \ldots : b ) & 0 & 12
    \\ \hline
    \end{array}
\]

Letting $a := \frac{1 - \epsilon}{2}$ and $b := \frac{1 - \epsilon^2}{2}$, the values of the coordinates can be completed using the
description of $\mathscr S$ above.
In $\{ Q = 0 \}$, the associated curve $S$ is described by the equation
\[
	\sigma_1 \sigma_3 - \sigma_4 = 0.
\]
There are 32 nodes lying over $S$, so $\# S_\text{sing} = 16$.
Therefore $S$ decomposes into $8$ lines. 
Consequently, all singularities of type 0, 1, 2, 5 and 6 are resolved by curves that are not nullhomologous.
This must also hold for the nodes of type 3 and 4, since they can be taken to singularities of the other types by permutation of coordinates.
Thus there exist projective algebraic small resolutions. \bigskip

\wernerpage{78}

The principle of the associated curve can be further generalized as follows:
consider the hypersurface
\[
    V = \{ x_4^{n - 2} Q + R = 0 \}
\]
of degree $n$ in $\mathbb P^4$;
we assume that $\{ Q( x_0, \ldots, x_3 ) = 0 \}$ and $\{ R(x_0, \ldots, x_3 ) = 0 \}$ are smooth surfaces in $\mathbb P^3$ of degrees 2 and $n$, respectively.\footnote{Note that in the $n=4$ discussion above, the surface cut out by $R$ was not assumed to be smooth.}
The curve
\[
    S := \{ Q = 0 \} \cap \{ R = 0 \} \subset \mathbb P^3
\]
is of type $(n, n)$ in $Q$, decomposing into $d' + 1$ irreducible components $S_0$ through $S_{d'}$.
$S_\text{sing}$ consists only of ordinary double points.
The map $\Pi$ is defined as above.  $P_1 := ( 0 : 0 : 0 : 0 : 1 ) \in \mathscr S$.
The map $\Pi$ can again be restricted to $V$ and then extended to $V( P_1 )$. \bigskip

Using the assumption that $R$ is smooth, one proves the analogue of the $n = 4$ case:

All singularities of $V$ other than $P_1$ are mapped by $\Pi$ to singularities of $S$;
over each singularity of $S$ lie exactly $(n - 2)$ singularities of $V$.
If $S$ has only ordinary double points, then so does $V$ --- there exists a $\lambda \neq 0$ such that
\[
    \lambda \partial_i Q + \partial R_i = 0
\]
for all $i = 0, \ldots, 3$; one solves the equation
\[
    x_4^{n-2} = \lambda.
\]
Therefore $s = 1 + (n - 2) \# S_{\text{sing}}$.

\wernerpage{79}

The cones $\Pi^{-1}( S_i )$ over the components of $S$ are divisors on $V$;
a point $P \in \mathscr S$ whose image under the map $\Pi$ is at the intersection of two (locally smooth) components of $S$ is resolved by a curve in $\hat V$ that is not nullhomologous in $H_2( \hat V, \mathbb Q )$.
$P_1$ is resolved by a nullhomologous curve if and only if none of the components of the associated curve $S$ is of type $(a, b)$ with $a \neq b$.

For the defect we have
\[ d \geq d'. \]

\addtocounter{thm}{-1}
\begin{conj}
$H_4( V, \mathbb Q )$ is generated by the divisors $\Pi^{-1} S_i$.\footnote{And, presumably, the hyperplane class.}
\end{conj}

This conjecture was proved for $n = 3$.
But in the general case we have not been able to prove that no divisor exists that passes smoothly through given singularities and is not mapped by $\Pi$ into $S$.\footnote{This statement is stronger than would be needed to prove the conjecture (which deals with homology classes, not actual divisors) and is not even true in the case $n=3$.  Consider cubic threefolds containing a plane, which will be studied at the end of chapter \ref{chapter10} and written there in the form $x_2 Q_1 + x_4 Q_2 = 0$.  Generically they have four nodes, all lying on the plane, and by the table in chapter \ref{chapter8}, they have two small resolutions.  One is obtained by blowing up along the plane $x_3=x_4=0$, the other by blowing up along the smooth quadric surface $x_3 = Q_2 = 0$.  But the projection $\Pi$ from a node does not map the latter divisor into the associated curve $S$, even though its homology class is the hyperplane class minus the class of the plane.}
The consequences of the conjecture would be:
\begin{enumerate}[label=\arabic*.]
    \item For a point $P \in \mathscr S - \{ P_1 \}$, if $\Pi( P )$ is a singularity of an irreducible component of the associated curve, then $P$ is resolved by a nullhomologous curve in a small resolution $\hat V$.
    
    \item $d = d'$.
\end{enumerate}

But even without a proof of this conjecture, the following characteristic examples are interesting;
they suggest various possibilities for combining the number of singularities, the defect, and nullhomologous curves.

\begin{exa}
\wernerpage{80}
The associated curve $S$ decomposes into two smooth components of type $(1, 0)$ and $(n - 1, n )$.
Then
\begin{align*}
    s = 1 + (n - 2) n &= (n - 1)^2
    \qquad \text{and}
    \\
    d' &= 1.
\end{align*}
No singularity is resolved nullhomologously.

As one can see, $( n - 1 )^2$ nodes in very special position are already enough to get a positive defect and projective algebraic small resolutions.
All double points lie in a plane contained in $V$.
Blowing up this divisor leads to a projective algebraic small resolution.

We will encounter this example again in the next chapter, where we will show that in fact the defect $d=1$.

The $( n - 1, n )$-component of $S$ can have up to $(n - 2)(n - 1)$ singularities; the associated $(n - 2)^2(n - 1)$ nodes of $V$ are presumably resolved nullhomologously.
\end{exa}

\begin{exa}
The associated curve $S$ decomposes into two smooth components of type $( 1, n - 1 )$ and $( n - 1, 1 )$.
We have
\begin{align*}
    s = 1 + (n - 2) \big( (n - 1)^2 + 1 \big) &= (n - 1)( n^2 - 3n + 3 )
    \qquad \text{and}
    \\
    d' &= 1.
\end{align*}
No singularity is resolved nullhomologously; there are projective algebraic small resolutions.

\wernerpage{81}

Again $d'$ is constant, but the number of singularities is larger than in the first case.
\end{exa}

\begin{exa}
The curve $S$ decomposes into $n$ components of type $(1, 0)$ and $n$ of type $(0, 1)$.
Then we have
\begin{align*}
    s = 1 + ( n - 2 ) n^2 &= (n - 1)( n^2 - n - 1 ) \qquad \text{and}
    \\
    d' &= 2n - 1.
\end{align*}
Again projective algebraic small resolutions exist.
\end{exa}

\begin{exa}
Now $S$ decomposes into $n - 1$ components of type $(1, 0)$ and one of type $(1, n)$.
We have
\begin{align*}
    s &= 1 + (n - 2)(n - 1) n \qquad \text{and}
    \\
    d' &= n - 1.
\end{align*}
The cones over the $(1, 0)$-components are smooth divisors on $V$, and every singularity lies in such a divisor.
Thus projective algebraic small resolutions exist.
\end{exa}

Both examples show that the defect can be unbounded with increasing $n$, without $s$ having to grow asymptotically like $n^4$ as in the case of Chmutov hypersurfaces.

\begin{exa}
\wernerpage{82}
Now $S$ decomposes into $n - 2$ components of type $(1, 0)$ and one $(2, n)$-component.
We have
\begin{align*}
    s &= 1 + (n - 2)^2 n \qquad \text{and}
    \\
    d' &= n - 2.
\end{align*}
As in the previous case, a smooth divisor passes through each singularity, so there exist projective algebraic small resolutions.
\end{exa}

However, it is possible that the $(2, n)$-component of the associated curve has up to $n - 1$ self-intersection points and hence $V$ possesses up to \linebreak 
$(n - 1)(n - 2)$ additional singularities.
These singularities are probably resolved nullhomologously.

\begin{exa}
This time $S$ decomposes into a $(1, 1)$-component and one $(n - 1, n - 1)$-component.
Then
\begin{align*}
    s &= 1 + 2(n - 2)(n - 1) \qquad \text{and}
    \\
    d' &= 1.
\end{align*}
The point $P_1$ is always resolved nullhomologously, so despite the defect being positive there are no projective algebraic small resolutions.
\end{exa}

This also holds for
\begin{exa}
The associated curve $S$ consists of $n$ components of type $(1, 1)$.
We have
\begin{align*}
    s &= 1 + (n - 2)(n - 1) n \qquad \text{and}
    \\
    d' &= n - 1.
\end{align*}
\wernerpage{83}
Furthermore, in this example the defect is unbounded.
Yet no projective algebraic small resolutions exist, since $P_1$ is resolved nullhomologously.
\end{exa} \bigskip

Supposing the conjecture is true, it is even possible to construct cases in which both the defect and the number of nullhomologous curves are unbounded:

The curve $S$ consists of $(n - m)$ curves of type $(1, 1)$ and one $(m, m)$-curve with $(m - 1)^2$ self-intersection points, where $m := \lfloor n / 2 \rfloor$.

As the example of Chmutov hypersurfaces shows, the defect can even grow asymptotically like $n^2 / 2$.
In chapter \ref{chapter12}, we will construct examples with $d \geq (n - 1)^2$ in other ways.
\chapter{Further examples} \label{chapter10}
\wernerpage{84}

The maximum number $\mu_4( 5 )$ of ordinary double points of a quintic $\mathbb P^4( \mathbb C )$ is not known; to date, we only have the bound
\[
    126 \leq \mu_4( 5 ) \leq 135.
\]
Schoen \cite{schoen86} and Hirzebruch \cite{hirzebruch85} construct examples of quintics with many nodes, which we briefly present here. \bigskip

Schoen's quintic
\[
    \left\{ \sum_{i = 0}^{4} x_i^5 - 5 \prod_{i = 0}^{4} x_i = 0 \right\}
\]
has 125 nodes, namely the orbit of the point $( 1 : 1 : 1 : 1 : 1 )$ under the action of the group generated by the coordinate transformations
\[
    ( x_0 : \ldots : x_4 ) \longmapsto ( x_0 : x_1 \xi_5^{a_1} : \ldots : x_4 \xi_5^{a_4} )
\]
with $\sum_{i = 1}^{4} a_i \equiv 0 \pmod 5$ and $\xi_5 := e^{2\pi i / 5}$.
This group acts transitively on $\mathscr S$.
Schoen calculates the defect as $d = 24$; there are projective algebraic small resolutions. \bigskip

Hirzebruch's quintic has 126 nodes.
For its construction, consider the curve $\{ f( x, y ) = 0 \}$ in the affine $(x, y)$-plane, given by the
\wernerpage{85}
product $f = \prod_{i = 1}^{5} f_i$ of the 5 lines of a regular pentagon:

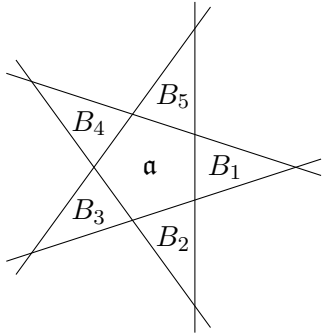
\[
\begin{tikzpicture}
\draw (0,0) node {$\mathfrak a$};
\draw (.6, 2.2) -- (.6, -2.2);
\draw (1, -0) node {$B_1$};
\draw (2.27773, .109203) -- (-1.90691, -1.25047);
\draw (.309017, -.951057) node {$B_2$};
\draw (.807717, -2.13251) -- (-1.77854, 1.42717);
\draw (-.809017, -.587785) node {$B_3$};
\draw (-1.77854, -1.42717) -- (.807717, 2.13251);
\draw (-.809017, .587785) node {$B_4$};
\draw (-1.90691, 1.25047) -- (2.27773, -.109203);
\draw (.309017, .951057) node {$B_5$};
\end{tikzpicture}
\]

As a function of two real variables, $f$ has relative extrema at the center of the pentagon as well as at an interior point $b_j$ of each triangle $B_j$.
Thus both partial derivatives vanish at these 6 points and at the 10 intersection points of the lines.
By symmetry, we have $f( b_i ) = f( b_j )$ for all $i$ and $j$.
The equation
\[
    f( u, v ) - f( z, w ) = 0
\]
defines a quintic in $\mathbb C^4$, and after homogenizing, a quintic in $\mathbb P^4$ with 126 ordinary double points.
Let $P \in \mathscr S$ be one of the 100 nodes whose coordinates arise as intersections of two lines. 
At least one of the meromorphic functions
\[
    \frac{ f_i( u, v ) }{ f_j( z, w ) }
\]
has a point of indeterminacy at $P$.

\wernerpage{86}

The divisors
\[
    \{ z  = \xi_5^i u ,\, w = \xi_5^i v \} \qquad ( i = 0, \ldots, 4 )
\]
on $V$ are smooth and contain the remaining 26 nodes.
Thus no singularity is resolved nullhomologously; there exist projective algebraic small resolutions, as Schoen first established. \bigskip

We switch back to double solids with branching surfaces with many nodes.
Unfortunately I neither know an explicit formula for the sextic in $\mathbb P^3$ with $s = 66$,\footnote{See footnote \ref{65_vs_66_footnote} on page \pageref{65_vs_66_footnote}: it turns out that $\mu_3(6) = 65$.} nor the Kreiss octic with 160 double points.
Thus, in general, we can only say that in both cases $2 s > \beta_3( V_t )$, so $d > 0$. \bigskip

Instead we consider examples of quartics which Kummer already examined in \cite{kummer72}.
Let $\{ F = 0 \} \subset \mathbb P^3$ be a smooth quadric that intersects the edges of the tetrahedron $\{ x_0 x_1 x_2 x_3 = 0 \}$ in 12 distinct points.
Then the quartic
\[
    B = \{ x_0 x_1 x_2 x_3 + F^2 = 0 \}
\]
has at least these 12 points as nodes.
The equation of the associated double solid can be written as
\[
    \omega^2 - F^2 = x_0 x_1 x_2 x_3.
\]
Successive blowing up of the divisors $\{ \omega - F = 0, x_i = 0 \}$ for $i = 0, 1, 2$ yields a partial resolution with $\beta_2( \check V ) \geq 4$ in which these 12 singularities are resolved.

\wernerpage{87}

In the case where there are no further nodes --- that is, $s = 12$ --- we have $d = 3$, and the partial resolution above is already a projective algebraic small resolution.
The following relations hold between the exceptional curves:

The two intersection points of a tetrahedron edge $\{ x_i = 0,\, x_j = 0 \}$, $i \neq j$, with the quadric are resolved by curves that are homologous to one another; 
let $L_{i j}$ denote a representative of the homology class of the resolving curve.
There are 3 further relations in $H_2( \hat V, \mathbb Z )$ between the 6 curves $L_{i j}$:
\begin{align*}
    L_{01} - L_{02} + L_{12} &\simeq 0
    \\
    L_{02} - L_{03} + L_{23} &\simeq 0
    \\
    L_{01} - L_{03} + L_{13} &\simeq 0.
\end{align*}
There are exactly 24 projective algebraic small resolutions. \bigskip

Choose a specific 
\[
    F = \frac{1}{4} \sum_{i = 0}^{3} x_i^2,
\]
so $B$ has 16 singularities altogether, and $d = 6$.
The 4 new singularities are $(1 : 1 : 1 : 1 ) - 2 e_i$ for $i = 0, 1, 2, 3$.
There are additional smooth divisors that also contain these new singularities; for example, write the equation of the branching surface as
\[
    -2 ( F - x_0 x_1 )^2 = x_0 x_1 ( x_0 - x_1 + i x_2 + i x_3 ) ( x_0 - x_1 - i x_2 - i x_3 ).
\]
There are projective algebraic small resolutions, but the relations between the exceptional curves are more complicated than in the case with $s = 12$.
Interestingly, the new divisors allow the points $\{ x_i = 0, x_j = 0 \} \cap \{ F = 0 \}$, $i \neq j$, to be separated; the resolving curves are no longer homologous. \bigskip

\wernerpage{88}

In the case $s = 13$ we still have $d = 3$, because the new singularities do not lie on the quadric $\{ F = 0 \}$.
The partial resolution described above has $\beta_2( \check V ) = 4$, so the 13th double point is resolved nullhomologously.
There are no projective algebraic small resolutions.
Such an example plays an important role in Poon \cite{poon92}. \bigskip

With
\[
    F = \frac{1}{4} \sum_{i = 0}^{3} x_i^2 + x_0 x_1
\]
one obtains $s = 14$ and $d = 4$.
The new singularities lie at $( 1 : -1 : i : -i )$ and $( 1 : -1 : -i : i )$.
Write the equation of the double solid as
\[
    \omega^2 - ( F - x_0 x_1 )^2 = \frac{x_0 x_1 ( x_0 + x_1 + i x_2 + i x_3 ) ( x_0 + x_1 - i x_2 - i x_3 )}{2}
\]
The divisor
\[
    \{ \omega = F - x_0 x_1, x_0 + x_1 = i ( x_2 + x_3 ) \} \subset V
\]
passes 
smoothly through the new singularities, which therefore are not resolved nullhomologously.
There exist projective algebraic small resolutions.
It is worth mentioning that the points $\{ x_i = 0, x_j = 0 \} \cap \{ F = 0 \}$ for $(i, j) \in \{ (0, 1), (2, 3) \}$ are still resolved homologously, as are the two new singularities;
all others can be separated by smooth divisors on $V$.

\wernerpage{89}

These examples can be further generalized by replacing the coordinate functions $x_0, \ldots, x_3$ by 4 arbitrary linear forms $\ell_0, \ldots, \ell_3$, any three of which are linearly independent.
Two of the Chmutov quartics that we encountered in chapter \ref{chapter7} are of this kind.

\bigskip
    
According to Kummer \cite{kummer64}, every quartic in $\mathbb P^3$ with 16 ordinary double points can be expressed as
\[
    B = \{ \phi^2 = 16 K x_0 x_1 x_2 x_3 \}, 
\]
where
\[
    \phi = x_0^2 + x_1^2 + x_2^2 + x_3^2 + 2a ( x_0 x_1 + x_2 x_3 ) + 2b ( x_0 x_2 + x_1 x_3 ) + 2c ( x_0 x_3 + x_1 x_2 )
\]
and
\[
    K = a^2 + b^2 + c^2 - 2 a b c - 1,
\]
$a$, $b$, and $c$ being parameters.

The planes $\{ x_i = 0 \}$ each contain 6 nodes, which lie on the conic that the quadric $\{ \phi = 0 \}$ cuts out from the planes.

Overall there are 16 of these singular tangent planes,\footnote{Werner writes ``singul\"are Tangentialebene'' which is Kummer's original term from \cite{kummer64}.  In English these planes are sometimes called ``tropes,'' following Cayley, or ``double planes.''  They correspond to singular points of the dual surface, which happens to be isomorphic to the original surface $B$.  The intersection of any of these plane with $B$ is a double conic.}
6 of them passing through each singularity of the Kummer surface, and 6 singular points lying on each plane.

Let $\{ p = 0 \}$ and $\{ q = 0 \}$ be the equations of any two of these tangent planes.
Then there exist homogeneous polynomials $\phi, \psi \in \mathbb C[ x_0, \ldots, x_3 ]$ of degree 2 such that
\[
    B = \{ \phi^2 = p q \psi \}.
\]
If $B$ is a Kummer surface with 16 nodes, then one can choose $\phi$ so that $\psi$ decomposes into two linear factors, which in turn describe two of the singular tangent planes.

\wernerpage{90}

Now let $B = \{ F = 0 \}$ be a Kummer surface with 16 ordinary double points, and let the defining polynomial $F$ be fixed.
Furthermore, let $P_1 \in \mathscr S$ and let $\{ q_i = 0 \}$ $(i = 1, \ldots, 6)$ be the singular tangent planes through $P_1$.
Then there are linear polynomials $r_{i j}$ and $s_{i j}$ as well as homogeneous polynomials $\phi_{i j}$ of degree 2 (which are uniquely determined up to sign) such that
\[
    F = \phi_{i j}^2 + q_i q_j r_{i j} s_{i j}.
\]
The double solid $V$ is given by the equation
\[
    \omega^2 - \phi_{i j}^2 = q_i q_j r_{i j} s_{i j}.
\]
For all $i, j, k \in \{1, \ldots, 6\}$, the difference $\phi_{i j}^2 - \phi_{i k}^2$ is divisible by $q_i$, so
\[
    \phi_{i j} = \pm \phi_{i k} + q_i \psi_i.
\]
If we now set
\[
    D_i := \{ q_i = 0, \omega - \phi_{i j} - 0 \},
\]
this definition is independent of $j$.  The divisors $D_i$ on $V$ are well-defined, and their homology classes in $H_4( V, \mathbb Z )$ are uniquely determined up to sign.

Any two of the divisors $D_i$ and $D_j$ $(i \neq j)$ have two singularities in common;
one is $P_1$, and the other we call $P_{i j}$.
Blowing up the divisors
\[
    \{ \omega - \phi_{i j} = 0,\ q_i = 0 \} \qquad\text{and}\qquad \{ \omega - \phi_{i j},\ q_j = 0 \}
\]
yields partial resolutions that resolve both $P_1$ and $P_{i j}$ differently.\footnote{That is, the birational map from the blow-up along $D_i$ to the blow-up along $D_j$ flops both the curve over $P_1$ and the curve over $P_{ij}$.}

\wernerpage{91}

Thus for the intersection numbers of the proper transforms of $D_i$ and $D_j$ with the exceptional curves $L_1$ and $L_{i j}$ on a small resolution, we have
\[
    \hat D_i . L_1 = \hat D_i . L_{i j} \iff \hat D_j . L_1 = \hat D_j . L_{i j}.
\]
Given a particular 
small resolution and in terms of a particular choice of divisors $D_i$, the intersection numbers are given by the following table:
\begin{center}
\resizebox{\textwidth}{!}{
$\begin{array}{c | c c c c c c c c c c c c c c c c c }
     & L_1 & L_{1 2} & L_{1 3} & L_{1 4} & L_{1 5} & L_{1 6} & L_{2 3} & L_{2 4} & L_{2 5} & L_{2 6} & L_{3 4} & L_{3 5} & L_{3 6} & L_{4 5} & L_{4 6} & L_{5 6} \\
     \hline
     D_1 & 1 & 1 & 1 & 1 & 1 & 1 & 0 & 0 & 0 & 0 & 0 & 0 & 0 & 0 & 0 & 0 \\
     D_2 & -1 & -1 & 0 & 0 & 0 & 0 & 1 & 1 & 1 & 1 & 0 & 0 & 0 & 0 & 0 & 0 \\
     D_3 & 1 & 0 & 1 & 0 & 0 & 0 & -1 & 0 & 0 & 0 & 1 & 1 & 1 & 0 & 0 & 0 \\
     D_4 & -1 & 0 & 0 & -1 & 0 & 0 & 0 & 1 & 0 & 0 & -1 & 0 & 0 & 1 & 1 & 0 \\
     D_5 & 1 & 0 & 0 & 0 & 1 & 0 & 0 & 0 & -1 & 0 & 0 & 1 & 0 & -1 & 0 & 1 \\
     D_6 & -1 & 0 & 0 & 0 & 0 & -1 & 0 & 0 & 0 & 1 & 0 & 0 & -1 & 0 & 1 & 1 \\
\end{array}$
}
\end{center}

\paragraph{Consequences:}
\begin{enumerate}[label=\arabic*.]
    \item The rank of the matrix is 6; since the defect is also 6, the divisors $D_i$ together with the image of a generator of $H_4( V_t, \mathbb Z )$ form a basis for $H_4( V, \mathbb Z )$.
    
    \item There are projective algebraic small resolutions; the exact number can be determined using the table above.
    
    \item The number of projective algebraic small resolutions for any double solid with $\deg B = 4$ and $s = 16$ is the same.
\end{enumerate}
\bigskip

\wernerpage{92}

Until now, we have mainly occupied ourselves with varieties with many singularities and their small resolutions.
Now we ask, how many nodes must a hypersurface in $\mathbb P^4$ or a double solid have in order for a projective algebraic small resolution to exist? \bigskip

Consider
\[
    V = \{ x_3 Q_1 + x_4 Q_2 = 0 \}
\]
with $Q_1, Q_2 \in \mathbb C[ x_0, \ldots, x_4 ]$ homogeneous of degree $(n - 1)$, such that $V$ only has ordinary double points at $\{ x_3 = 0,\,x_4 = 0,\,Q_1 = 0,\,Q_2 = 0 \}$.
The curves $Q_i( x_0 : x_1 : x_2 : 0 : 0 )$ intersect the plane $\{ x_3 = 0, x_4 = 0 \}$ transversely in $(n - 1)^2$ points, so $s = (n - 1)^2$.
Blowing up the plane $\{ x_3 = 0,\, x_4 = 0 \}$ yields a projective small resolution.

\begin{prop}
For all $n \geq 2$, the defect is equal to 1.
\end{prop}

\begin{proof}\!\!\!\footnote{Alternatively one could use a Koszul resolution of the ideal sheaf $\mathscr I_{\mathscr S}$ to compute $h^1(\mathscr I_{\mathscr S}(2n-5)) = 1$.}
The case $n = 2$, $s = 1$, and $d = 1$ has already been dealt witht in chapter \ref{chapter4}.
Now let $n \geq 3$.
Consider homogeneous polynomials $f$ of degree $(2n - 5)$ in $\mathbb P^4$, and write them in the form 
\[
    f = f_1( x_0, x_1, x_2 ) + x_3 g + x_4 h.
\]
Polynomials that are divisible by $x_3$ or $x_4$ vanish on $\mathscr S$.
Polynomials of degree $(2n - 5)$ in $\mathbb P^2$ with homogeneous coordinates $( x_0 : x_1 : x_2 )$ form a space of dimension $\binom{2n - 3}{2} = 2n^2 - 7n + 6$.

\wernerpage{93}

By Max Noether's fundamental theorem (\cite{fulton78}, p.~120), for every $f_1( x_0, x_1, x_2 )$ of degree $(2n - 5)$ with $f_1 |_{\mathscr S} = 0$, we have
\[
    f_1 = A Q_1( x_0, x_1, x_2, 0, 0 ) + B Q_2( x_0, x_1, x_2, 0, 0 )
\]
with $A, B \in \mathbb C[ x_0, x_1, x_2 ]$ homogeneous of degree $(n - 4)$.
Thus the dimension of the space of these polynomials $f_1$ equals $2\binom{n - 2}{2} = n^2 - 5n + 6$.
For the defect we get
\[
    d = (n^2 - 2n + 1) + (n^2 - 5n + 6) - ( 2n^2 - 7n + 6 ) = 1,
\]
and the proposition is proven in general.
There are exactly two projective algebraic small resolutions.
\end{proof}

Every hypersurface in $\mathbb P^4$ of degree $n$ with $(n - 1)^2$ nodes that lie in a plane can be brought into the form described above.

\begin{conj}
If $\deg V = n$ and $s < (n - 1)^2$, then $d = 0$.
\end{conj}

Analogous observations can also be made for double solids.
Let
\[
    B = \{ x_3 Q_1 + Q_2^2 \} \subset \mathbb P^3
\]
with $Q_1, Q_2 \in \mathbb C[ x_0, x_1, x_2]$ homogeneous of degrees $(n - 1)$ and $n / 2$, respectively, such that $B$ only has ordinary double points at $\{ x_3 = 0,\,Q_1 = 0,\,Q_2 = 0 \}$.
The curves $Q_i( x_0, x_1, x_2, 0 )$ intersect the plane $\{ x_3 = 0 \}$ transversely in $(n-1)\,n/2$ points,
\wernerpage{94}
so the double solid
\[
    \omega^2 - Q_2^2 = x_3 Q_1
\]
has exactly $(n-1)\,n/2$ nodes.
As in the preceding example, one shows that $V$ has defect 1 for all $n \geq 2$.

\begin{conj}
Let $B$ of degree $n$ be the branching surface of a double solid with $s < (n - 1)^2$.  Then $d = 0$.
\end{conj}
\chapter{Generalizations} \label{chapter11}
\wernerpage{95}

In the last two chapters we shall attempt to generalize some of the previous results: rather than double solids and hypersurfaces in $\mathbb P^4$, we will investigate general projective threefolds with isolated singularities, and rather than restricting our attention to ordinary double points, we will also allow higher singularities. \bigskip

We begin with the most general case.
Let $V$ be a projective threefold with isolated singularities.
All singularities are allowed that admit a small resolution, meaning one in which the singularity is resolved by a curve.
The set of singularities is again denoted by $\mathscr S$, and for $P \in \mathscr S$, let $L_P^i$ be the various irreducible curves of a small resolution of $P$; the index $i$ runs from 1 to $s_P$.
It is again possible to characterize the projective algebraic small resolutions.

\begin{thm}
For a small resolution $\hat V$, the following statements are equivalent:
\begin{enumerate}[label=\arabic*.]
    \item $\hat V$ is projective algebraic.
    
    \item There exists a divisor $D$ on $\hat V$ with
    \[
        D . L_P^i > 0
    \]
    for all $i = 1, \ldots, s_P$ and all $P \in \mathscr S$.

    \item There is no nontrivial relation between the irreducible exceptional curves
\wernerpage{96}
    \[
        \sum_{P \in \mathscr S} \sum_{i = 1}^{s_P} \beta_P^i L_P^i \simeq 0,
    \]
    such that $\beta_P^i \geq 0$ for all $i = 1, \ldots, s_P$ and all $P \in \mathscr S$.
\end{enumerate}
\end{thm}

\begin{proof}
The implications ``1 $\Rightarrow$ 2 $\Rightarrow$ 3'' are clear.
The proof of the direction ``3 $\Rightarrow$ 1'' based on Peternell's results [given in Theorem \ref{thm1}] can be carried over 
to the general case verbatim.

In order to adapt 
the alternative proof [given after Theorem \ref{thm1}] to some specific cases, a few preliminary considerations are necessary. \bigskip

Define
\[
    s := \sum_{P \in \mathscr S} s_P
\]
and define $\mathscr B$ to be the $\mathbb Q$-vector space generated by the $s$-tuples $(\beta_1, \ldots, \beta_s)$ of coefficients of the linear relation between the irreducible exceptional curves $L_i$.
Note that we are simply enumerating the exceptional curves here without regard for which curve resolves which singular point.
The $\mathbb Q$-vector space $\mathscr A$ is generated by the $s$-tuples $(\alpha_1, \ldots, \alpha_s)$ of intersection numbers of the irreducible exceptional curves with the divisors on $\hat V$.
It is clear that $\mathscr A$ and $\mathscr B$ are once again orthogonal to each other in $\mathbb Q^s$.
Less clear is how $\mathscr A$ and $\mathscr B$ behave when the small resolution $\hat V$ changes.

\wernerpage{97}

One can also define a defect $d := \beta_4( V ) - \beta_2( V )$ so that
\begin{align*}
    \beta_4( \hat V ) &= \beta_4( V ) \qquad\text{and}
    \\
    \beta_2( \hat V ) &= \beta_2( V ) + d.
\end{align*}
The vector space $\mathscr B$ has dimension $(s - d)$. \bigskip

If $P \in \mathscr S$ is of type $(2, 2, n + 1, n + 1)$, meaning that $V$ is given in a neighborhood of $P$ in affine coordinates by
\[
    z_1^2 + z_2^2 + z_3^{n + 1} + z_4^{n + 1} = 0,
\]
then one can again pass to a deformation $V_t$, locally described by an equation
\[
    z_1^2 + z_2^2 + z_3^{n + 1} + z_4^{n + 1} = t.
\]
We have
\[
    \beta_4( V_t ) = \beta_2( V_t ) = \beta_2( V ) = \beta_4( V ) - d.
\]
Thus the dimension of $\mathscr A$ equals $d$, and we have $\mathbb Q^s = \mathscr A \perp \mathscr B$. \bigskip

Between the vanishing cycles on $V_t$ there are $d$ relations, which are described by $\mathscr A$.

Notice that of the $n^2$ vanishing cycles, only $n$ at a time are involved in global relations. 
A precise description 
would take us too far afield.

Using these considerations, the alternative proof can be extended at least to the case of singularities of type $(2, 2, n + 1, n + 1)$.
\end{proof}

\wernerpage{98}

In the case of ordinary double points, a change of small resolution still entails a change in the signs in the relations between the exceptional curves.
Hence we have

\begin{thm}
Let $V$ be a nodal projective threefold. 
A projective algebraic small resolution exists exactly when all irreducible exceptional curves on $\hat V$ are not nullhomologous.
\end{thm}
\bigskip

Our first new examples are complete intersections of two quadrics in $\mathbb P^5$.

In \cite{poon86}, Poon examines the three dimensional variety
\[
    \left\{ 2(z_0^2 + z_1^2) + \lambda z_2^2 + \tfrac{3}{2} z_3^2 + z_4^2 + z_5^2 = 0,\ \sum_{i = 0}^5 z_i^2 = 0 \right\} \subset \mathbb P^5,
\]
where $\frac{3}{2} < \lambda < 2$, which has 4 nodes at $(1 : \pm i : 0 : 0 : 0 : 0)$ and $(0 : 0 : 0 : 0 : 1 : \pm i)$.

Poon found that the defect is equal to 2, and that the two pairs of points are each resolved homologously (modulo signs).
Thus there are 4 projective algebraic small resolutions. \bigskip

If $\lambda = \frac{3}{2}$ then $s = 6$, because there are additional singularities at \linebreak 
$(0 : 0 : 1 : \pm i : 0 : 0)$.
The divisor
\[
    \{ z_0 + i z_1 = 0,\, z_2 + i z_3 = 0,\, z_4 + i z_5 = 0 \}
\]
is smooth on $V$ and contains the points $(1 : i : 0 : 0 : 0 : 0)$, \linebreak 
$(0 : 0 : 1 : i : 0 : 0)$, and $(0 : 0 : 0 : 0 : 1 : i)$.
Similarly, the divisor
\[
    \{ z_0 - i z_1 = 0,\, z_2 + i z_3 = 0,\, z_4 + i z_5 = 0 \}
\]
\wernerpage{99}
goes through $(1 : -i : 0 : 0 : 0 : 0)$, $(0 : 0 : 1 : i : 0 : 0)$, and \linebreak 
$(0 : 0 : 0 : 0 : 1 : i)$, and so on.
The defect equals 4, and the intersection behavior of the divisors with the curves $L_P$ on a particular small resolution of $V$ are described by the matrix
\[
    \begin{pmatrix}
    1 & 1 & 0 & 0 & 0 & 0 \\
    0 & 0 & 1 & 1 & 0 & 0 \\
    0 & 0 & 0 & 0 & 1 & 1 \\
    1 & 0 & 1 & 0 & 1 & 0
    \end{pmatrix}.
\]
The vector space $\mathscr B$ is spanned by the vectors $(1, -1, -1, 1, 0, 0)$ and \linebreak 
$(1, -1, 0, 0, -1, 1)$.
Of the 64 small resolutions, exactly 46 are projective algebraic. \bigskip

These examples suggest the conjecture that with complete intersections in higher-dimensional projective spaces, the defect and thus the number of projective algebraic small resolutions $\hat V$ is higher in comparison to the number of nodes than with hypersurfaces in $\mathbb P^4$.
So far, however, there is no general formula for calculating the defect. \bigskip

We now turn to higher singularities.  Let $V$ be a singular threefold with $P \in \mathscr S$ given in local affine coordinates by the equation
\[
    z_1^2 + z_2^2 + z_3^{n + 1} + z_4^{n + 1} = 0.
\]
\wernerpage{100}
A suitable change of coordinates brings this equation into the form
\[
    z_1 z_2 = \prod_{j = 0}^n ( z_3 + \xi_{n + 1}^j z_4 ),
\]
where $\xi_{n + 1} := e^{2 \pi i / (n + 1)}$.
Set $f_j := z_3 + \xi_{n + 1}^j z_4$.
Then $\{ f_j = 0 \}$ and $\{ z_i = 0 \}$ describe local Cartier divisors, which decompose into Weil divisors.
According to Brieskorn \cite{brieskorn66}, one 
small resolution of $P$ is given by the closure of the graph of the local meromorphic function
\[
    \left( \frac{z_1}{f_0},\, \frac{z_1}{f_0 f_1},\, \ldots,\, \frac{z_1}{f_0 f_1 \cdots f_{n - 1}} \right)
\]
in $\mathbb P^N \times (\mathbb P^1)^n$.
The exceptional curves are
\[
    L_{P}^i := \{ P \} \times ( 0 : 1 ) \times \cdots \times \mathbb P^1 \times ( 1 : 0 ) \times \cdots \times ( 1 : 0 ),
\]
with $i = 1, \ldots, n$ and $\mathbb P^i$ in the $i$-th position.
Thus $s_P = n$.

\begin{lem}
The proper transforms of the Weil divisors $\{ z_i = 0, f_j = 0 \}$, $(j = 0, \ldots, n)$ have the following intersection numbers with the curves $L_P^i$:
\[
    \begin{array}{l l}
    \{ z_1 = 0, f_0 = 0 \} : & (-1, 0, 0, \ldots, 0, 0)
    \\
    \{ z_1 = 0, f_1 = 0 \} : & (1, -1, 0, \ldots, 0, 0)
    \\
    \{ z_1 = 0, f_2 = 0 \} : & (0, 1, -1, 0 \ldots, 0)
    \\
    \hspace{3.5em} \vdots
    \\
    \{ z_1 = 0, f_{n - 1} = 0 \} : & (0, \ldots, 0, 0, 1, -1)
    \\
    \{ z_1 = 0, f_n = 0 \} : & (0, \ldots, 0, 0, 0, 1).
    \end{array}
\]
The $i$-th component denotes the intersection number with the curve $L_P^i$.
\end{lem}

\begin{proof}
\wernerpage{101}
Consider the case $n = 2$.
The local equation of $V$ in an affine neighborhood of $P$ is
\[
    z_1 z_2 = \prod_{j = 0}^2 f_j.
\]
The closure of the graph of the meromorphic function
\[
    \left( \frac{z_1}{f_0},\, \frac{z_1}{f_0 f_1} \right)
\]
describes a small resolution with exceptional curves $\mathbb P^1 \times (1 : 0)$ and \linebreak 
$(0 : 1) \times \mathbb P^1$.
The proper transform of the Weil divisor $\{ z_2 = 0, f_0 = 0 \}$ intersects these two curves only at the point $(1 : 0) \times (1 : 0)$; the intersection numbers are therefore $(1, 0)$.
As the Cartier divisor $\{ f_0 = 0 \}$ must have intersection number 0 with all exceptional curves, the intersection numbers of these curves with the Weil divisor $\{ z_1 = 0, f_0 = 0 \}$ are $(-1, 0)$.

The proper transforms of $\{ z_1 = 0, f_1 = 0 \}$ and $\{ z_2 = 0 , f_1 = 0 \}$ contain $(0 : 1) \times \mathbb P^1$ and $\mathbb P^1 \times (0 : 1)$, respectively.
This gives the intersection numbers with the exceptional curves as $(1, -1)$ and $(-1, 1)$, respectively.
This procedure also allows one to calculate all intersection numbers in the general case.
\end{proof}

\wernerpage{102}

According to Brieskorn, the other small resolutions are described by the graphs of the functions
\[
    \left( \frac{z_1}{f_{j_0}}, \,\frac{z_1}{f_{j_0} f_{j_1}},\, \ldots,\, \frac{z_1}{f_{j_0} f_{j_1} \cdots f_{j_n}} \right),
\]
and thus correspond to permutations of the above tuples of intersection numbers.
There are $(n + 1)!$ different possibilities for the small resolution. \bigskip

We will now attempt to characterize global properties of small resolutions using this local information.
Let $V \subset \mathbb P^N$ be a singular threefold where at each $P \in \mathscr S$ there exists a $n \in \mathbb N$ such that $V$ is given in an affine neighborhood of $P$ by an equation of the form
\[
    z_1^2 + z_2^2 + z_3^{n + 1} + z_4^{n + 1} = 0.
\]
To address the existence of projective algebraic small resolutions, one should search for global divisors on $V$ whose proper transforms on a small resolution $\hat V$ have positive intersection with all exceptional curves.
Whether this search is successful depends on how the local Weil divisors
\[
    W_P^j := \{ z_1 = 0, f_j = 0 \},
\]
$P \in \mathscr S$, can be extended to global divisors.

\wernerpage{103}

Let us assume that all the local Weil divisors $W_P^j$ through a point $P \in \mathscr S$ can be extended to distinct global divisors, each of which is smooth at $P$.
Consider the linear combination
\[
    D := \sum_j \alpha_P^j D_P^j,
\]
$\alpha_P^j \in \mathbb N_0$.
If the $\alpha_P^j$ are all distinct --- without loss of generality we can assume that $\alpha_P^0 < \ldots < \alpha_P^n$ --- then there is exactly one small resolution of the point $P$ such that the proper transform of $D$ intersects all the curves $L_P^i$ positively: one must resolve in such a way that $W_P^0$ is assigned the tuple of intersection numbers $(-1, 0, 0, \ldots, 0, 0)$, $W_P^1$ the tuple $(1, -1, 0, \ldots, 0, 0)$, and so on.

If two of the numbers $\alpha_P^i$ are equal, then it is not possible to find a small resolution with
\[
    \hat D . L_P^i > 0
\]
for all $i = 1, \ldots, n$, as one can see from the intersection numbers.
This also holds 
if several $W_P^i$ can only be completed to a single global divisor, which is then singular at $P$.

\begin{thm}
Let $V$ be a threefold with the above properties.
There exists a projective algebraic small resolution if and only if for all singularities $P \in \mathscr S$, all local Weil divisors $W_P^i$ can be extended to global divisors $D_P^i$ on $V$ which are smooth at $P$.
\end{thm}
\begin{proof}
\wernerpage{104}
If this condition is satisfied, consider all linear combinations
\[
    D := \sum_{P \in \mathscr S} \sum_{j = 0}^{s_P} \alpha_P^j D_P^j,
\]
$\alpha_P^j \in \mathbb N_0$, as well as their restrictions $D_P$ to affine neighborhoods of all $P \in \mathscr S$
\[
    D_P =: \sum_{j = 0}^{s_P} \beta_P^j W_P^j.
\]
For a point $P \in \mathscr S$, the set of divisors $D$ with $\beta_P^i \neq \beta_P^j$ for all $i \neq j$ is nonempty.
Then there also exists a divisor $D$ such that $\beta_P^i \neq \beta_P^j$ for all $i \neq j$ and all $P \in \mathscr S$.
By our earlier considerations, for each such divisor $D$ there is a small resolution $\hat V$ such that 
\[
    \hat D . L_P^i > 0
\]
for all $i = 1, \ldots, s_P$ and all $P \in \mathscr S$.
The different projective algebraic small resolutions are then in one-to-one correspondence with the possible orderings of the $\beta_P^i$ at each singularity.
\end{proof}

\begin{cor}
Let $V$ be as in the previous theorem. If there exist projective algebraic small resolutions, then $V$ has only $(2, 2, n + 1, n + 1)$-singularities with $n \leq d$.
\end{cor}
\chapter{Examples with higher singularities} \label{chapter12}
\wernerpage{105}

As in the case of ordinary double points, singular cubics in $\mathbb P^4$  provide nice examples to which we can apply the general theory. \bigskip

In this chapter, let $K \subset \mathbb P^4$ be a cubic with singularities only of type $(2, 2, n + 1, n + 1)$.
Furthermore, we suppose that $K$ has at least one ordinary double point, which lies at the point $(0 : 0 : 0 : 0 : 1)$.
Then we can write the cubic as
\[
    \{ x_4 Q + R = 0 \},
\]
where $Q$ and $R$ are homogeneous polynomials of degrees 2 and 3, respectively, 
in $(x_0 : \ldots : x_3)$.
We again call
\[
    S := \{Q = 0\} \cap \{R = 0\} \subset \mathbb P^3
\]
the associated curve.  The map $\Pi \colon K - \{ P_1 \} \to \mathbb P^3$ sends a point $P \in \mathscr S - \{ P_1 \}$ of type $(2, 2, n + 1, n + 1)$ to a singularity of $S$ at which $(n + 1)$ local components intersect, pairwise transversely.
Since $S$ is of type $(3, 3)$, only the cases $n = 1, 2, 3$ occur.
The cones over the irreducible components $S_i$ of $S$ form a basis for $H_4( K, \mathbb Z )$,\footnote{This is stronger than the statement in chapter \ref{chapter8} that these divisors and the class of the generic hyperplane section generate $H_4( K, \mathbb Z )$.  It may even be too strong, although the meaning of ``basis'' is ambiguous since $H_4( K, \mathbb Z )$ typically has torsion.} and we have the following

\begin{thm}
Let $K$ be a cubic with the properties given above.
A projective algebraic small resolution exists if and only if all irreducible components of $S$ are smooth and at least one is of type $(a, b)$ with $a \neq b$.
\end{thm}

\wernerpage{106}

The ample divisors on the various projective algebraic small resolutions $\hat V$ arise as proper transforms of
\[
    \sum_{i = 0}^d \alpha_i\, \Pi^{-1} S_i,
\]
where the $\alpha_i$ are chosen so that $\alpha_i \neq \alpha_j$ if $S_i \cap S_j \neq \varnothing$.
Note also that we must have $\sum_{i = 0}^d \alpha_i (a_i - b_i) \neq 0$, where $S_i$ is of type $(a_i, b_i)$.

The different projective algebraic small resolutions correspond to the different possibilities, at each singularity $P \in \mathscr S - \{ P_1 \}$, for ordering the $\alpha_i$ with $\Pi( P ) \in S_i$. 
The specific 
resolution of $P_1$ is determined by the sign of $\sum_{i = 0}^d \alpha_i (a_i - b_i) \neq 0$. \bigskip

These considerations immediately yield the following

\begin{thm}
Let $K^1$ be another cubic, whose associated curve $S^1$ arises from a deformation of $S$ in which a $(2, 2, n + 1, n + 1)$-singularity splits into $\binom{n + 1}{2}$ nodes.
Then $K^1$ and $K$ have the same number of projective algebraic small resolutions.
\end{thm}

This is to be understood as the quadric $Q$ remaining unchanged while the associated curve in $Q$ and the cubic $R$ are deformed.

\begin{remark}
\wernerpage{107}
An analogous theorem can be formulated for arbitrary singular threefolds $V$ with only isolated $(2, 2, n + 1, n + 1)$-singularities.  Let $\mathscr W$ be the total space of a deformation of $V$, let
\[
    \Pi : \mathscr W \to \mathbb C
\]
be the projection, let $V = V^0 = \Pi^{-1}( 0 )$ and let $V^1 = \Pi^{-1}( 1 )$.
Suppose the deformation has 
the following properties:
\begin{enumerate}[label=\arabic*.]
    \item A $(2, 2, n + 1, n + 1)$-singularity in $V$ splits into $\binom{n + 1}{2}$ nodes in $V^1$.
    
    \item There exist divisors $D_1, \ldots, D_N$ on $\mathscr W$ such that $D_1 |_{V^i}, \ldots, D_N |_{V^i}$ are the global divisors\footnote{Presumably meaning the ones that extend the local divisors $W_P^i$ as in the last part of chapter \ref{chapter11}.} on $V^i$, $i = 0, 1$.
\end{enumerate}
Then $V$ and $V^1$ have the same number of projective algebraic small resolutions. \bigskip
\end{remark}

We return to the cubic.
Examples with higher singularities can be found as deformations of nodal cubics given in chapter \ref{chapter8};
in this way it is possible to classify all cubics with only isolated singularities of type $(2, 2, n + 1, n + 1)$ and at least one ordinary double point.

The numbers of projective algebraic small resolutions calculated in chapter \ref{chapter8} can be verified in the cases $d \leq 3$ using elementary methods.
Let us take a closer look at the case $d = 4$.

\wernerpage{108}

The associated curve
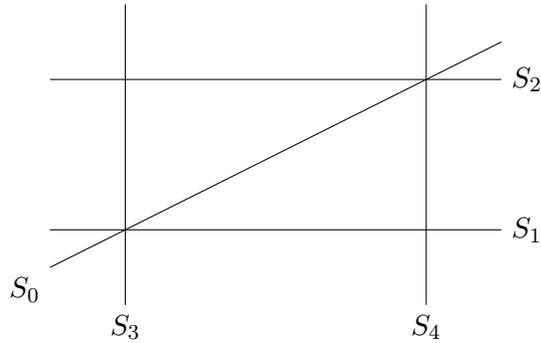
\[
\begin{tikzpicture}
\draw (-2, -2) -- (-2, 2);
\draw (2, -2) -- (2, 2);
\draw (-3, -1) -- (3, -1);
\draw (-3, 1) -- (3, 1);
\draw (-3, -1.5) -- (3, 1.5);
\draw (-3, -1.5) node[anchor=north east] {$S_0$};
\draw (-2, -2) node[anchor=north] {$S_3$};
\draw (2, -2) node[anchor=north] {$S_4$};
\draw (3, -1) node[anchor=west] {$S_1$};
\draw (3, 1) node[anchor=west] {$S_2$};
\end{tikzpicture}
\]
belongs to a cubic with three nodes and two singularities of type $(2, 2, 3, 3)$.
The divisors on this cubic are given by
\[
    \sum_{i = 0}^4 \alpha_i \, \Pi^{-1} S_i.
\]
Again we assign to a divisor $D$ the small resolution with $\hat D . L_P^i > 0$ for all $i = 1, \ldots, s_P$ and all $P \in \mathscr S$.
We set
\[
    \alpha_0 < \alpha_1 < \alpha_3
\]
and thus fix one of the 6 possible small resolutions of the singularity corresponding to $S_0 \cap S_1 \cap S_3$.
We wonder: how many projective algebraic small resolutions are possible given this constraint?

If $\alpha_4 < \alpha_0$, there are 6 cases:
\[
    \begin{array}{lllll}
    \alpha_2 < \alpha_4, & & & & 
    \\
    \alpha_4 < \alpha_2 < \alpha_0 & & \text{and} & & \alpha_1 + \alpha_2 < \alpha_3 + \alpha_4
    \\
    & & \text{or} & & \alpha_1 + \alpha_2 > \alpha_3 + \alpha_4,
    \\
    \alpha_0 < \alpha_2 < \alpha_3 & & \text{and} & & \alpha_1 + \alpha_2 < \alpha_3 + \alpha_4
    \\
    & & \text{or} & & \alpha_1 + \alpha_2 > \alpha_3 + \alpha_4,
    \\
    \alpha_2 > \alpha_3. & & & & \\
    \end{array}
\]

\wernerpage{109}

If $\alpha_0 < \alpha_4 < \alpha_1$, there are 5 distinct possibilities:
\[
    \begin{array}{lllll}
    \alpha_2 < \alpha_0, & & & & 
    \\
    \alpha_0 < \alpha_2 < \alpha_4, & & & &
    \\
    \alpha_4 < \alpha_2 < \alpha_3 & & \text{and} & & \alpha_1 + \alpha_2 < \alpha_3 + \alpha_4
    \\
    & & \text{or} & & \alpha_1 + \alpha_2 > \alpha_3 + \alpha_4,
    \\
    \alpha_2 > \alpha_3. & & & & \\
    \end{array}
\]
If $\alpha_4 > \alpha_1$, the following 6 cases are possible:
\[
    \begin{array}{lllll}
    \alpha_2 < \alpha_0, & & & & 
    \\
    \alpha_0 < \alpha_2 < \min\{ \alpha_3, \alpha_4 \}, & & & &
    \\
    \alpha_4 < \alpha_2 < \alpha_3, & & & &
    \\
    \alpha_2 > \max\{ \alpha_3, \alpha_4 \}& & \text{and} & & \alpha_1 + \alpha_2 < \alpha_3 + \alpha_4
    \\
    & & \text{or} & & \alpha_1 + \alpha_2 > \alpha_3 + \alpha_4,
    \\
    \alpha_3 < \alpha_2 < \alpha_4. & & & & \\
    \end{array}
\]
Changing the sign of $(\alpha_1 + \alpha_2) - (\alpha_3 + \alpha_4)$ corresponds to changing the small resolution of $P_1$.
Altogether there are 102 projective algebraic small resolutions, which is also true for a cubic with 9 nodes. \bigskip

What can be said about the relations between the exceptional curves on the different small resolutions?

In the case of higher singularities, the change of relations resulting from a change in small resolution is not as simple to describe as in the
\wernerpage{110}
case of nodes, where it caused a change in the signs of the corresponding coefficients.
In general, a change in the small resolution of a point $P \in \mathscr S$ entails a permutation of the intersection numbers of the proper transforms of divisors that contains $P$ with $L_P^i$.
This alteration of the vector space $\mathscr A$ causes an alteration of the vector space $\mathscr B$ that describes the linear relations between the exceptional curves $L_P^i$. \bigskip

Here is a typical example.
Let $K$ be a cubic whose associated curve decomposes into 3 components of type $(2, 2)$, $(0, 1)$, and $(1, 0)$ in the following way:
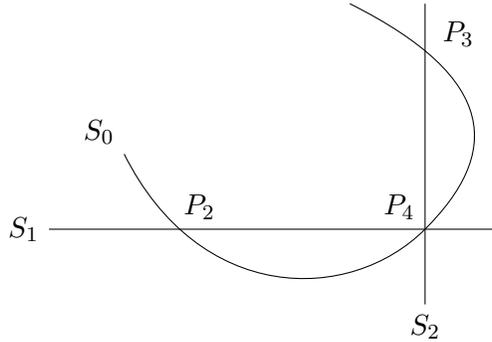
\[
\begin{tikzpicture}
\draw (2, -2) -- (2, 2);
\draw (-3, -1) -- (3, -1);
\draw (-2,0)
      .. controls (-1,-2) and (1,-2)
      .. (2,-1)
      .. controls (3,0) and (3,1)
      .. (1,2);
\draw (-2,0) node[anchor=south east] {$S_0$};
\draw (-1,-1) node[anchor=south] {$P_2$};
\draw (2,-1) node[anchor=south east] {$P_4$};
\draw (2.1,1.6) node[anchor=west] {$P_3$};
\draw (2, -2) node[anchor=north] {$S_2$};
\draw (-3, -1) node[anchor=east] {$S_1$};
\end{tikzpicture}
\]
One possible (projective) small resolution comes from assigning the divisors $\Pi^{-1} S_1$ and $\Pi^{-1} S_2$ the intersection numbers $(-1, 1, 0, 1, -1)$ and $(1, 0, 1, 0, 1)$, respectively, with $L_1$, $L_2$, $L_3$, $L_4^1$, and $L_4^2$.
An ample divisor
\[
    \sum_{i = 0}^2 \alpha_i\, \widehat{\Pi^{-1} S_i }
\]
satisfies the condition
\[
    \alpha_0 < \alpha_1 < \alpha_2.
\]

\wernerpage{111}

Passing to a small resolution whose ample divisor satisfies
\[
    \alpha_1 < \alpha_0 < \alpha_2,
\]
entails a change in the small resolution of the points $P_2$ and $P_4$:
$\Pi^{-1} S_1$ and $\Pi^{-1} S_2$ are now assigned the intersection numbers $(-1, -1, 0, -1, 0)$ and $(1, 0, 1, 0, 1)$, respectively.
The linear relations in the first case are
\[
    L_1 \simeq L_4^2,\, L_2 \simeq L_4^1,\, L_3 \simeq L_4^1 + L_4^2,
\]
and in the second,
\[
    L_3 \simeq L_4^2, L_2 \simeq L_4^1, L_1 \simeq L_4^1 + L_4^2.
\]\smallskip

Also interesting is the following example with $d = 1$.  The associated curve $S$ decomposes into two components of types $(3, 2)$ and $(0, 1)$, 
the $(3, 2)$-component having a double point which coincides with one of the two intersection points of the two components.
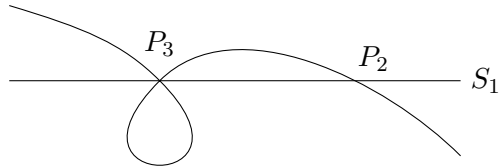
\[
\begin{tikzpicture}
\draw (-3, 0) -- (3, 0);
\draw (-3,1)
      .. controls (-2,0.7) and (-1.5,0.5)
      .. (-1, 0)
      .. controls (0.5,-1.5) and (-2.5,-1.5)
      .. (-1,0)
      .. controls (0,1) and (2,0)
      .. (3, -1);
\draw (3, 0) node[anchor=west] {$S_1$};
\draw (-1, 0.2) node[anchor=south] {$P_3$};
\draw (1.5, 0) node[anchor=south west] {$P_2$};
\end{tikzpicture}
\]
The intersection numbers of the divisor $\widehat{\Pi^{-1} S_1}$ with $L_1$, $L_2$, $L_3^1$, and $L_3^2$ are $(\pm 1, \pm 1, 0, 1)$, $(\pm 1, \pm 1, -1, 0)$, or $(\pm 1, \pm 1, 1, -1)$.
All sign combinations are independent from one another.
In the first case $L_3^1$ is nullhomologous; in the second case $L_3^2$ is nullhomologous;
in the third case, neither irreducible exceptional curve is nullhomologous, but instead the sum $L_3^1 + L_3^2$ is. \bigskip

\wernerpage{112}

Thus with higher singularities, whether an irreducible exceptional curve is nullhomologous or not may depend on the particular 
small resolution.
In contrast, the following condition is independent of the 
small resolution:
\begin{enumerate}[label=$(*)$]
    \item A singularity $P \in \mathscr S$ is resolved by curves $L^i$, between which there is a nontrivial relation
    \[
        \sum_{i = 1}^{s_P} \alpha^i L^i \simeq 0,
    \]
    where $\alpha^i \geq 0$ for all $i = 1, \ldots, s_P$
\end{enumerate}

This 
holds not only for cubics in $\mathbb P^4$, but in general for threefolds with isolated $(2, 2, n + 1, n + 1)$-singularities.
A projective algebraic small resolution exists if and only if the condition $(*)$ is not satisfied for any $P \in \mathscr S$.\bigskip

With that we leave cubics and turn to special 
double solids whose branching surfaces (generalizing Clemens' original definition) can now also have higher singularities.

Recall an example from chapter \ref{chapter10}; there we investigated a double solid with branching surfaces of the form
\[
    B = \{ Q^2 - x_0 x_1 x_2 x_3 = 0 \} \subset \mathbb P^3
\]
where $\deg Q = 2$.
The general such surface has 12 nodes.
But if one chooses $Q$ such that the quadric $\{ Q = 0 \}$
\wernerpage{113}
goes through the four vertices of the tetrahedron $\{ x_0 x_1 x_2 x_3 = 0 \}$, then $B$ has four $D_4$-singularities.
According to Kummer \cite{kummer66}, in this case we have
\[
    Q = \sum_{0 \leq i < j \leq 3} a_{i j} x_i x_j.
\]
The singularities of the double solid are of type $(2, 2, 3, 3)$, so they correspond to the triple cover of an $A_2$-surface singularity
and thus allow a small resolution.
The double solid is again described by the equation
\[
    ( \omega - Q ) ( \omega + Q ) = x_0 x_1 x_2 x_3,
\]
and can be deformed to a general double solid with $\deg B = 4$ and $s = 12$.
Hence $d = 3$, and the divisors
\[
    D_i = \{ \omega - Q = 0, x_i = 0 \}
\]
$i = 0, \ldots, 3$, generate $H_4( V, \mathbb Z )$.
Three of these global divisors pass through each singularity, where they represent the 3 different local divisors.
Thus projective algebraic small resolutions exist.
To calculate how many, we consider all divisors
\[
    D = \sum_{i = 0}^3 \alpha_i D_i,
\]
with $\alpha_i \geq 0$.
There are 24 possible orderings the $\alpha_i$ at the singular points.
Since any two divisors intersect in two singularities and represent distinct local divisors there, this leads to 24 different projective algebraic small resolutions.
This agrees with the number we calculated in chapter \ref{chapter10} in the case of 12 nodes.

\wernerpage{114}

If we choose all the constants $a_{i j}$ in the equation for $Q$ equal to $\frac{1}{2}$, then $B$ has three more singularities, namely nodes at the points $(1 : 1 : -1 : -1)$, $(1 : -1 : 1 :-1)$, and $(1 : -1 : -1 : 1)$.
The first two nodes lie on the smooth divisor
\[
    \{ x_0 + x_3 = 0,\, \omega = x_1 x_2 - x_0 x_3 \},
\]
which in turn lies in $V$.
One obtains analogous divisors by transposing $x_3$ with $x_1$ or $x_2$.
Thus there are global smooth divisors through all 3 nodes, ensuring the existence of a projective algebraic small resolution in this example as well. \bigskip

Finally, the following beautiful example of a family of hypersurfaces in $\mathbb P^4$ should not go unmentioned:

Let
\[
    V = \{ F( z_0, z_1, z_2 ) + z_3^{n + 1} + z_4^{n + 1} = 0 \}, 
\]
where $\{ F = 0 \}$ describes a nodal curve of degree $(n + 1)$ in $\mathbb P^2$.
The singularities of this curve are in one-to-one correspondence with the singularities of $V$, which are all of type $(2, 2, n + 1, n + 1)$.

Write $F = \prod_{i = 1}^m f_i$ ;
we assume furthermore that all irreducible components of the curve $\{ F = 0 \}$ are smooth.
Then the only singularities are the intersection points of the irreducible components in question,
\wernerpage{115}
hence one gets a projective algebraic small resolution of $V$ by blowing up the divisors
\[
    \{ f_i = 0, z_3 + \xi_{n + 1}^j z_4 = 0 \},
\]
where $1 \leq i \leq m - 1$ and $1 \leq j \leq n$.
The number of singularities depends on the number of components of the curve $\{ F = 0 \}$ and their degrees.
The two extreme cases are:
\begin{enumerate}[label=\arabic*.]
    \item $F$ factors into two polynomials, of degrees $n$ and $1$.
    Then $V$ has $n$ singularities, and $d \geq n$.
    
    \item $F$ factors into $(n + 1)$ linear terms.
    In this case, $V$ has $\binom{n + 1}{2}$ singularities, and $d \geq n^2$.
\end{enumerate}

If we now deform $V$ so that each singularity splits into $\binom{n + 1}{2}$ nodes, then the new variety $V^1$ has exactly as many projective algebraic small resolutions as $V$.
The nodes lie in very special position.
In the first case we have $s = n^2( n + 1 ) / 2$ and $d \geq n$, and in the second, $s = n^2( n + 1 )^2 / 4$ and $d \geq n^2$.
Note that $\deg V^1 = \deg V = n + 1$.

\begin{conj}
Let $V$ be a hypersurface of degree $(n + 1)$ in $\mathbb P^4$, and let all its singularities be of type $(2, 2, n + 1, n + 1)$, with the same $n$ for each singularity.
Then the number of singularities is less than or equal to $\binom{n + 1}{2}$.  It is greater than or equal to $n$ if $V$ admits a projective algebraic small resolution.
\end{conj}

\clearpage
\phantomsection
\addcontentsline{toc}{chapter}{Bibliography}
\nocite{*}
\bibliographystyle{plain}
\bibliography{ref}
\end{document}